\documentclass[10pt]{article}
\usepackage{microtype}
\usepackage{amsmath,amssymb,amsthm,mathtools}
\usepackage[a4paper,margin=29mm]{geometry}
\usepackage{tikz-cd}
\usetikzlibrary{arrows.meta}
\usepackage{float}
\usepackage{algorithm,algpseudocode}
\usepackage{enumitem,booktabs,graphicx}
\usepackage[hidelinks,pdftitle={Discrete Helicity: Linking on Triangulated Skeletons and Algorithms},pdfauthor={Jingqi Li}]{hyperref}
\setlist[enumerate]{leftmargin=2em,itemsep=2pt,topsep=4pt}
\setlist[itemize]{leftmargin=1.7em,itemsep=2pt,topsep=4pt}

\newif\ifshowalgorithms
\showalgorithmstrue

\makeatletter
\renewcommand{\@maketitle}{%
  \null\vskip 1em
  \begin{center}
    {\normalfont\normalsize\bfseries\MakeUppercase{\@title}\par}
    \vskip 1.2em
    {\normalfont\footnotesize\MakeUppercase{\@author}\par}
  \end{center}
  \vskip .6em
}
\renewcommand{\@seccntformat}[1]{\textup{\csname the#1\endcsname.}\quad}
\renewcommand{\section}{\@startsection{section}{1}{\z@}%
  {-3.5ex \@plus -1ex \@minus -.2ex}{1.5ex \@plus .2ex}%
  {\normalfont\normalsize\scshape\centering}}
\renewcommand{\subsection}{\@startsection{subsection}{2}{\z@}%
  {-3.25ex \@plus -1ex \@minus -.2ex}{1ex \@plus .2ex}%
  {\normalfont\normalsize\bfseries}}
\renewcommand{\subsubsection}{\@startsection{subsubsection}{3}{\z@}%
  {-3.25ex \@plus -1ex \@minus -.2ex}{1ex \@plus .2ex}%
  {\normalfont\normalsize\itshape}}
\makeatother
\renewenvironment{abstract}
  {\footnotesize
   \list{}{\leftmargin=2.5em\rightmargin=2.5em
     \topsep=0pt\partopsep=0pt\parsep=0pt\itemsep=0pt}
   \item\relax\noindent\textsc{\abstractname.}\enspace\ignorespaces}
  {\endlist\vspace{1em}}

\allowdisplaybreaks[2]
\numberwithin{equation}{section}
\theoremstyle{plain}
\newtheorem{theorem}{Theorem}[section]
\newtheorem{main}[theorem]{Theorem}
\newtheorem{introassumption}{Assumption}
\newtheorem{proposition}[theorem]{Proposition}
\newtheorem{lemma}[theorem]{Lemma}
\newtheorem{corollary}[theorem]{Corollary}
\theoremstyle{definition}
\newtheorem{definition}[theorem]{Definition}
\newtheorem*{definition*}{Definition}
\newtheorem{assumption}[introassumption]{Assumption}
\newtheorem{example}[theorem]{Example}

\theoremstyle{definition}
\newtheorem{remark}[theorem]{Remark}
\newtheorem*{remark*}{Remark}
\newcommand{\R}{\mathbb{R}}
\newcommand{\Hc}{H_{\mathrm c}}
\newcommand{\Hw}{H_{\mathrm w}}
\newcommand{\Hsm}{\mathcal{H}}
\DeclareMathOperator{\Lk}{Lk}
\newcommand{\cpr}{c_{\mathrm{pr}}}
\newcommand{\cdu}{c_{\mathrm{du}}}

\DeclareMathOperator{\diag}{diag}

\newcommand{\Cross}{\mathcal X}

\DeclareMathOperator{\sgn}{sgn}
\newcommand{\E}{\mathbb E}
\newcommand{\Prob}{\mathbb P}

\ifshowalgorithms
\title{Discrete Helicity: A Linking explanation on\\Triangulations and Algorithms}
\else
\title{Discrete Helicity: A Linking explanation on\\Triangulations}
\hypersetup{pdftitle={Discrete Helicity: A Linking Explanation on Triangulations}}
\fi
\author{Jingqi Li}
\date{}
\begin{document}
\maketitle
\begin{abstract}
We interpret discrete helicity in terms of linking numbers on a closed oriented simply connected 3-manifold.
We first introduce cup helicity, defined with the cup product on a triangulation of the manifold. Combining cup--cap duality with the local antisymmetrization formula for the Whitney product~\cite{SchubelEtAl2024}, we express Kotiuga's discrete helicity, defined using Whitney interpolation, in terms of linking numbers~\cite{Kotiuga1990}.
We further construct continuous-time Markov processes on the primal and dual vertex sets and use their trajectories to define discrete  flow lines. These yield a discrete helicity theorem of the same expression as its classical counterpart, expressing Whitney helicity as a weighted average of asymptotic linking numbers.
Also, inspired by the extrusion construction idea in \cite{Hirani2003}, we construct a new discrete Lie derivative operator and use its relation to the Markov generators to provide a geometric justification for our stochastic dynamics. \ifshowalgorithms
Finally, we use the linking interpretation to develop deterministic and randomized algorithms for computing discrete helicity.
\fi 
\end{abstract}
\section{Introduction}\label{sec:introduction}

Helicity connects the local description of a divergence-free field with
the global topology of its trajectories~\cite{Moffatt1969,Vogel2003,ArnoldKhesin2021}.  In fluid mechanics and
magnetohydrodynamics, it measures the mutual winding of vortex or magnetic
flux and constrains its evolution under ideal transport.  Moffatt's
interpretation in terms of linked vortex filaments made this connection
explicit~\cite{Moffatt1969}.  Let $X$ be a smooth divergence-free field on a
closed oriented 3-manifold with volume form $\mu$.  Assume $\iota_X\mu$ is exact, where $\iota$ denotes the contraction operator, and choose a globally defined real 1-form $\alpha$
with $\iota_X\mu=d\alpha$.  Choose any vector field  $A$ such that $\operatorname{curl}A=X$.
The helicity of $X$ is defined by
\[
 \Hsm(X):=\int_M A\cdot X\,\mu=\int_M\alpha\wedge d\alpha.
\]
The expression in terms of differential forms is independent of the choice of
potential and of the auxiliary metric, and is invariant under
volume-preserving diffeomorphisms~\cite{ArnoldKhesin2021}.

Arnold also proved that helicity bounds the energy from below~\cite{ArnoldKhesin2021}:
\[
\mathcal{E}(X):=\frac12\int_M |X|^2\,\mu
\ge C_M\,|\Hsm(X)|,
\]
where $C_M>0$ depends only on $M$.
Thus a field with nonzero helicity cannot have its energy reduced
arbitrarily close to zero by volume-preserving deformations.

For notational simplicity, unless otherwise stated, we make the following assumption throughout the article:

\begin{introassumption}
Let $M$ be a smooth, closed, oriented, simply connected $3$-manifold
with a fixed Riemannian metric $g$ and its induced volume form
$\mu=\operatorname{vol}_g$.
\end{introassumption}

In the following passage, we denote smooth helicity by $\mathcal{H}$
and its discrete counterparts by $H_c$ and $H_w$.

Arnold's theorem on asymptotic linking gives this integral a dynamical meaning.
Let $\varphi_t$ denote the flow of $X$. For $T,S>0$, append suitable short return paths to the trajectory segments $\varphi_{[0,T]}(x)$ and $\varphi_{[0,S]}(y)$, obtaining closed curves $\Gamma_T(x)$ and $\Gamma_S(y)$, respectively. Their normalized linking numbers have an $L^1$ limit:
\[
 \lim_{T,S\to\infty}
 \frac{\Lk\bigl(\Gamma_T(x),\Gamma_S(y)\bigr)}{TS}
 =\lambda_X(x,y)
 \quad\text{in }L^1(M\times M,\mu\otimes\mu).
\]
Arnold's theorem provides
\[
 \Hsm(X)=\int_M\int_M\lambda_X(x,y)\,\mu(x)\mu(y).
\]
Thus an integral defined through a potential can be recovered
from the long-time entanglement of the field's
trajectories~\cite[\S~III.4]{ArnoldKhesin2021}.  The choice and
control of closing paths are substantive parts of this statement;
Vogel supplies a construction and an $L^1$ formulation applicable in our
simply connected setting~\cite{Vogel2003}.  This correspondence between
helicity, linking, and long-time dynamics is the model
for our main result, the discrete helicity theorem stated in
Theorem~\ref{thm:main-discrete-helicity}.

The differential form expression also places magnetic helicity in direct
contact with Abelian Chern--Simons theory. For the globally defined real
potential $\alpha$, the Abelian Chern--Simons functional is
$\mathcal{S}_{\mathrm{CS}}(\alpha)=c_{\mathrm{CS}}\Hsm(X)$, where
$c_{\mathrm{CS}}$ specifies the chosen coupling and normalization
convention~\cite{Adams1996,Adams1997}. Discrete linking pairings have
also been studied in simplicial and lattice formulations of Abelian
Chern--Simons theory and its doubled (BF) version. This connection
provided one of the initial motivations for our present work.
For readers interested in this viewpoint, we refer to Adams's
primal--dual simplicial constructions, which recover linking numbers
through Wilson loop observables~\cite{Adams1996,Adams1997}, and to the
explicit lattice Gauss linking formulas of Di Bartolo, Grau, and
Leal~\cite{DiBartoloEtAl2013}.

  We develop a topological interpretation of discrete helicity in discrete exterior calculus (DEC) and the Whitney finite element setting of finite element exterior calculus (FEEC). More explicitly, we express Whitney helicity as the linking of flux-weighted primal and dual chains and recover it from the asymptotic linking of discrete flow trajectories. \ifshowalgorithms
This interpretation leads to practical deterministic and randomized algorithms for computing helicity from discrete flux data.
\fi

Whitney forms provide the finite element setting for this question.
Kotiuga's finite element approach to helicity~\cite{Kotiuga1990,GrossKotiuga2001}
uses Whitney interpolation  to reconstruct a real $1$-cochain $a$ as a
Whitney $1$-form $Wa$. In our notation, the resulting Whitney helicity is
\begin{equation}
\Hw(a):=\int_M Wa\wedge dWa.
\label{eq:intro-whitney-helicity}
\end{equation}
Although written as an integral, $\Hw(a)$ is a discrete quantity, since $Wa$ lies
in a finite-dimensional Whitney space and $\Hw(a)$ is a quadratic form
in the finitely many edge coefficients of $a$. We explain this construction
in more detail in \S\ref{sec:preliminaries}.
Kotiuga studied the matrices arising from this discretization~\cite{Kotiuga1990},
and computational work on discrete helicity includes the sparsity and
Lanczos methods of Chammas and Kotiuga~\cite{ChammasKotiuga1996}.
Gross and Kotiuga describe data structures based on simplicial chains
and cochains for tetrahedral finite element meshes~\cite{GrossKotiuga2001}.
The relevant cochain calculus draws on Whitney interpolation~\cite{Whitney1957},
product convergence under refinement~\cite{Wilson2007}, and discrete
exterior calculus~\cite{Hirani2003,DesbrunEtAl2005}; the averaging and
naturality results of Schubel, Berwick-Evans, and Hirani clarify the
local structure of the wedge product~\cite{SchubelEtAl2024}. Accordingly, the
contribution sought here is a specific connection between discrete 
helicity and  its flux linking explanation, building on these
existing algebraic, topological, and computational developments.

The main difficulty is that a cochain supplies neither a distinguished
vector field nor a time parametrization on the mesh. Hirani represents discrete
vector fields by vectors attached to primal or dual vertices~\cite[\S~5.1]{Hirani2003}.
The discrete vector bundles of Berwick-Evans, Hirani, and Schubel provide a
more general setting, in which sections are bundle-valued $0$-cochains~\cite{BerwickEvansEtAl2026}.
A different approach is Forman's combinatorial vector field: a pairing of
incident cells represented by an operator $V\colon C_p\to C_{p+1}$ that raises
chain degree by one~\cite{Forman1998Morse}. This construction also admits
combinatorial paths and a discrete flow~\cite{Forman1998Dynamics}.
However, these geometric and combinatorial descriptions do not, by themselves,
prescribe continuous-time trajectories compatible with the given flux
cochain and cell volumes required for a helicity theorem.

Our approach is to give the flux a stochastic dynamical realization.
Viewing each vertex as a contraction core of its dual cell, we first construct continuous-time Markov processes on the primal and dual
vertex sets, with jump rates determined by positive outgoing flux divided
by the volume of the departure cell. Each sample path produces a recording chain
by counting its signed edge traversals. These recording chains serve as
our discrete trajectories and describe the long-time transport encoded
by the flux. After closing them by bounded return paths, their normalized
linking numbers admit an asymptotic limit whose volume-weighted average
is Whitney helicity. Thus the Markov construction supplies both the time
parametrization and the long-time behavior needed for a discrete analogue
of Arnold's helicity theorem.

Motivated by cup--cap duality and the antisymmetrization formula relating
cup and wedge products in the work of Hirani et al.~\cite{Hirani2003,SchubelEtAl2024},
we first introduce a new cup helicity,
\[
 \Hc(a):=\langle a\smile\delta a,[K]\rangle.
\]
We give a linking realization of cup helicity and, through the
cup--wedge relation, also obtain a linking realization of Whitney helicity $H_w$.
We further give a heuristic flux interpretation of the chain
coefficients, providing a geometric rationale for this discretization.
We also establish gauge invariance of both helicities.
On the Whitney finite element space, we prove a discrete analogue of
Arnold's energy inequality.

Our main result is the following discrete helicity theorem. It realizes
Whitney helicity as the volume-weighted average of the asymptotic linking
of independent primal and dual Markov trajectories.

\begin{main}
\label{thm:main-discrete-helicity}
There exist trajectories on the primal and dual $1$-skeletons whose closures $\Gamma_t(x)$ and $\Gamma_s(y)$ satisfy,
with $\lambda$ extended piecewise constantly to $M\times M$,
\begin{equation}
 \lim_{t,s\to\infty}\frac{\Lk(\Gamma_t(x),\Gamma_s(y))}{ts}
 =\lambda(x,y),\qquad
 \Hw(a)=\int_{M\times M}\lambda(x,y)\,d\mu(x)\,d\mu(y).
 \label{eq:intro-discrete-helicity}
\end{equation}
\end{main}

This is the discrete analogue of Arnold's helicity theorem.
The trajectory construction and proof are given in
\S\ref{sec:markov}.

We next relate the stochastic dynamics to discrete Lie transport.
Inspired by the extrusion approach discussed by Hirani~\cite[\S~8.3]{Hirani2003},
we define contraction by taking a weighted average of the sweeps
associated with individual tetrahedra and obtain a
discrete Lie derivative through Cartan's formula.  A compatible
primal--dual pairing supplies the dual operators.  Under a particular geometric assumption, the scalar Lie derivatives on both the primal and dual skeletons coincide with the transport parts in the transport--diffusion decomposition of the corresponding Markov generators with respect to the inner products weighted by cell volumes.  This provides
a geometric account of the directed transport contained in the
stochastic model.

Also, for convergence, \S\ref{sec:convergence} establishes
first-order convergence of both cup and Whitney helicity to smooth
helicity under controlled shape-regular refinement. Under additional
volume compatibility, we obtain first-order convergence of the primal
scalar Lie derivative to the smooth one in the $L^\infty$ norm.

\ifshowalgorithms
Finally, on $S^3$, the linking identity yields deterministic and
randomized algorithms for evaluating Whitney helicity directly from
an exact $2$-cochain, without recovering a potential or solving a linear
system. Using polygonal realizations in a punctured chart that preserve
the linking of the two flux chains, we reduce the computation to a
weighted linking pairing. After geometric preprocessing, importance
sampling replaces the full crossing sum by a prescribed number of
sampled contributions. This can substantially reduce the summation cost
when the required sample size is small relative to the crossing count,
and an explicit probability bound controls the absolute error of the
resulting unbiased estimator. Methods for computing linking numbers in
computational geometry include formulas based on solid angles~\cite{Arai2013}, certified angle
summation~\cite{BertolazziEtAl2019}, and direct, Barnes--Hut, or fast
multipole evaluation of the Gauss linking integral~\cite{qu2021}.
We use the classical formula based on signed crossings under a generic planar
projection; these references provide interested readers with further
computational background and alternative methods.
\fi

The paper is organized as follows.
\S\ref{sec:preliminaries} reviews the geometric, algebraic, and probabilistic preliminaries.
\S\ref{sec:discrete-helicity} introduces cup and Whitney helicity,
establishes gauge invariance, and proves the linking formula for cup helicity. \S\ref{sec:whitney-flux} identifies Whitney helicity
with the linking of flux-weighted primal and dual chains.
\S\ref{sec:markov} constructs the Markov trajectories and proves the
discrete helicity theorem. \S\ref{sec:transport} develops the
discrete Lie derivative via extrusion and relates it to the Markov
generators. \S\ref{sec:convergence} proves convergence of the two
helicities and of the primal scalar Lie derivative. \ifshowalgorithms
Finally,
\S\ref{sec:algorithms} presents deterministic and randomized linking
algorithms together with numerical tests.
\fi

\section{Preliminaries}
\label{sec:preliminaries}

Throughout, $X$ denotes a smooth vector field on $M$ that is divergence-free. Write $V:=\int_M\mu$ for the
total volume of $M$.

We assume familiarity with the basic notions of simplicial complexes
and simplicial homology. For background, we refer the reader to
 Hatcher~\cite{Hatcher2002} and
 Weibel~\cite{Weibel1994}.
We also assume familiarity with calculus on smooth manifolds, including
 differential forms and integration; standard references are
Lee~\cite{Lee2012} and Bott and Tu~\cite{BottTu1982}.

\subsection{Smooth triangulations and orientation }
\label{sec:smooth-triangulations}
A smooth triangulation of $M$ consists of a simplicial complex $K$
and a homeomorphism $\Phi:|K|\to M$ whose restriction to each simplex
extends to a smooth map on a neighborhood of that simplex in its affine
span, and whose Jacobian  is nondegenerate at every point  including its boundary~\cite{Lurie2009}.
Here $| \cdot|$ denotes the geometric realization .
Whitehead's triangulation theorem guarantees the existence of such a
triangulation; since $M$ is compact, $K$ can be chosen
finite~\cite{Whitehead1940}. Throughout the following, we assume that $K$ is a finite smooth triangulation of $M$.

Let
\[
\widehat T:=\operatorname{conv}\{0,e_1,e_2,e_3\}\subset\R^3
\]
be the fixed  reference tetrahedron, where $e_1,e_2,e_3$ are
the standard basis vectors in $\mathbb R^3$. For each tetrahedron $T=\Phi(|\sigma_T|)$, where $\sigma_T\in K$, let
$A_T:\widehat T\to|\sigma_T|$ be the affine isomorphism determined by
the chosen vertex correspondence. Define the cell map
$\Phi_T:=\Phi|_{|\sigma_T|}\circ A_T$ by the commutative diagram
\[
\begin{tikzcd}[row sep=large, column sep=large]
\widehat T \arrow[r, "A_T"] \arrow[dr, "\Phi_T"']
  & \lvert\sigma_T\rvert \arrow[d, "\Phi|_{|\sigma_T|}"] \\
  & T\subset M.
\end{tikzcd}
\]
We choose the vertex correspondence so that $\Phi_T$ preserves the
standard orientation of $\widehat T$ and the orientation of $T$ induced
by $M$. These maps are smooth embeddings up to
the boundary and agree on shared faces after affine
identifications.

We now introduce an assumption used for the transport construction in
\S\ref{sec:transport} and its convergence analysis in
\S\ref{sec:scalar-convergence}:
\begin{assumption}\label{ass:volume-compatible}
For each tetrahedron $T$, set $V_T:=\int_T\mu$ and equip the reference
tetrahedron with the normalized volume form
$\widehat\mu_T:=6V_T\,dV_{\widehat T}$. The cell maps extend to $C^2$
diffeomorphisms between neighborhoods and satisfy
\[
\Phi_T^*\mu=\widehat\mu_T.
\]
Thus the local volume density relative to $\widehat\mu_T$ is one.

\end{assumption}
Here $dV_{\widehat T}$ is the Euclidean volume form; we choose this normalization so that the volume integrals over the
reference and physical tetrahedra agree.

Bertelson and Distexhe prove that any smooth triangulation of a manifold
with a smooth volume form admits a subdivision and a new smooth
triangulation whose pullback volume form is constant on each
top-dimensional simplex~\cite[\S~2]{BertelsonDistexhe2021}.
Consequently, after refinement and compatible piecewise smooth
reparametrization, one obtains a triangulation satisfying
Assumption~\ref{ass:volume-compatible}.

An \emph{ordered} simplicial complex is a simplicial complex equipped with
a total order on its vertices, which restricts to an order on each simplex.
Throughout, we assume that $K$ is ordered. Label its vertices
$v_1<\cdots<v_N$.
For a $p$-simplex with vertices $v_{i_0},\ldots,v_{i_p}$, where
$i_0<\cdots<i_p$, write
\[
[i_0,\ldots,i_p]:=[v_{i_0},\ldots,v_{i_p}]
\]
for the simplex oriented by this order. Reordering its vertices
multiplies the oriented simplex by the sign of the permutation.
We write $K^{(p)}$ for the set of $p$-simplices with reference orientations
given by increasing vertex order.

For a tetrahedron $T=\Phi(|\sigma_T|)$, list the vertices $v_0,v_1,v_2,v_3$
of $|\sigma_T|$ in its reference orientation. At a point $\xi$ in the interior
of $|\sigma_T|$, the differential $D(\Phi|_{|\sigma_T|})_\xi$ sends the ordered
edge vectors $(v_1-v_0,v_2-v_0,v_3-v_0)$ to a basis of $T_{\Phi(\xi)}M$.
Set $\varepsilon_T=1$ if this basis is positively oriented and
$\varepsilon_T=-1$ otherwise.
Thus multiplying the reference orientation of $T$ by $\varepsilon_T$
gives the orientation induced by $M$.
Since $M$ is closed and oriented,
each triangular face belongs to exactly two tetrahedra, whose ambient
orientations induce opposite boundary orientations on that face.
Consequently, the integer chain
\begin{equation}
[K]=\sum_{T\in K^{(3)}}\varepsilon_T T\in C_3(K;\mathbb Z),
\qquad \partial[K]=0
\label{eq:fundamental-cycle}
\end{equation}
is defined for every such triangulation and is a cycle. Indeed, under the
identification $H_3(K;\mathbb Z)\cong H_3(M;\mathbb Z) \cong \mathbb{Z}$ , $[K]$ represents the
       generator of this group
\cite{Hatcher2002}. 

Unless a coefficient ring is specified, chains and cochains have real
coefficients. We write $Z_\bullet(K)$, $B_\bullet(K)$, $Z^\bullet(K)$, and $B^\bullet(K)$
for  cycles, boundaries, cocycles, and coboundaries, respectively.
Also, by convention, a \emph{reference tetrahedron} is a tetrahedron in
Euclidean space, and a \emph{physical tetrahedron} is a tetrahedron in $M$.
Keeping this notion in mind will be helpful throughout the rest of the paper.
By a slight abuse of notation, we sometimes omit the hat on a reference
tetrahedron and use the terms \emph{reference} and \emph{physical} to
distinguish the two.
\begin{remark}

We emphasize that the orientation of a tetrahedron induced by the vertex order need
not agree with the orientation induced by $M$
\cite{Hirani2003}.
Unless otherwise specified, we orient each tetrahedron $T$ by increasing vertex order,
but when computing the face fluxes in Lemma~\ref{lem:local-flux-whitney} we orient
$T$ by $M$ and give its faces the induced boundary orientation. Readers must be very careful about this!
 Since $M$ is closed, it is compact, and its total volume
$V$ is finite.
Also, any simplicial complex becomes ordered after barycentric
subdivision~\cite{BerwickEvansEtAl2026}, so the assumption that $K$ is ordered is easily satisfied.

\end{remark}

\subsection{Products}
\label{sec:products}
For $a\in C^p(K)$ and $b\in C^q(K)$, define their cup product on simplices
written in increasing vertex order by
\begin{equation}
(a\smile b)[v_0,\ldots,v_{p+q}]
=a[v_0,\ldots,v_p]\,b[v_p,\ldots,v_{p+q}]
\label{eq:cup-product}
\end{equation}
and extend by orientation signs.
This product is associative and satisfies the Leibniz rule
\begin{equation}
\delta(a\smile b)=\delta a\smile b+(-1)^p a\smile\delta b,
\label{eq:cup-leibniz}
\end{equation}
where $\delta$ is the  coboundary map. However, it is generally not
commutative on cochains~\cite{Hatcher2002}.

The discrete wedge product, denoted by $\wedge $, is defined by antisymmetrizing the cup product~\cite{Hirani2003,SchubelEtAl2024}:
\begin{equation}
(a\wedge b)[v_0,\ldots,v_{p+q}]
=\frac{1}{(p+q+1)!}\sum_{\pi\in S_{p+q+1}}\operatorname{sgn}(\pi)\,
a[v_{\pi(0)},\ldots,v_{\pi(p)}]\,
b[v_{\pi(p)},\ldots,v_{\pi(p+q)}].
\label{eq:discrete-wedge-product}
\end{equation}
Here $S_{p+q+1}$ is the permutation group of $\{0,\ldots,p+q\}$.
The resulting product is independent of the chosen total vertex order.
Indeed, Hirani et al. showed in~\cite{SchubelEtAl2024} that the discrete
wedge product has a formula involving two successive averages. It is graded commutative and satisfies
\eqref{eq:cup-leibniz} with $\wedge $ in place of $\smile$,
but is generally not associative on cochains.

In \S\ref{sec:whitney-de-rham}, we will use the Whitney and de Rham
maps to explain the relation between this discrete wedge and the usual
wedge product of differential forms.

The cap product provides the corresponding pairing with chains.
For $a\in C^p(K)$ and an $m$-simplex with $p\leq m$, write its vertices
in increasing order and set
\begin{equation}
[v_0,\ldots,v_m]\frown a
=a[v_0,\ldots,v_p]\,[v_p,\ldots,v_m],
\label{eq:cap-product}
\end{equation}
and extend linearly. 
With this convention, we have the following cup--cap duality and boundary formula:
\begin{align}
\langle a\smile b,z\rangle
&=\langle b,z\frown a\rangle,
\qquad z\in C_{p+q}(K),
\label{eq:cup-cap-adjunction}\\
\partial(z\frown a)
&=(-1)^p\bigl((\partial z)\frown a-z\frown\delta a\bigr).
\label{eq:cap-boundary}
\end{align}

\begin{remark}
The cup and discrete wedge products induce the same associative,
graded-commutative product on $H^\bullet(K;\R)$: for cocycles $a$ and $b$,
\[
[a\wedge b]=[a\smile b].
\]
In particular, the identity on $H^\bullet(K;\R)$ is an isomorphism
between the two resulting graded rings~\cite{Whitney1957,Wilson2007}.
\end{remark}

\subsection{Whitney and de Rham map}
\label{sec:whitney-de-rham}

In this part we'll introduce two operators connecting  de Rham complex and simplicial cochain complex. 

Whitney interpolation map turns a cochain into a piecewise smooth differential
form.  On each reference tetrahedron, let $\lambda_v$ be the barycentric
coordinate associated with a vertex $v$. For an oriented simplex
$\sigma=[v_0,\ldots,v_p]$, let $\chi_\sigma$ denote the elementary
cochain that takes the value one on $\sigma$ and zero on the other
$p$-simplices with their reference orientations. Define
\begin{equation}
 W_p\chi_\sigma
 =p!\sum_{i=0}^{p}(-1)^i\lambda_{v_i}
 \,d\lambda_{v_0}\wedge\cdots\wedge
 \widehat{d\lambda_{v_i}}\wedge\cdots\wedge d\lambda_{v_p},
 \label{eq:whitney-form}
\end{equation}
and transport this form to each physical tetrahedron $T$ containing $\sigma$
via the pullback $(\Phi_T^{-1})^*$. We set $W_p\chi_\sigma=0$ outside
the closed star
\[
\overline{\operatorname{st}}(\sigma)
:=\bigcup_{\substack{\tau\in K\\ \sigma\subseteq\tau}}\tau.
\]  
Extending linearly gives $W_p a=\sum_\sigma a(\sigma)W_p\chi_\sigma$.
The resulting Whitney forms are piecewise smooth differential forms on $M$,
with matching tangential traces across shared faces, but are generally not
globally smooth. Their exterior derivatives and wedge products are taken
elementwise: 
\[
(d\omega)|_T:=d(\omega|_T),
\qquad
(\omega\wedge\eta)|_T:=(\omega|_T)\wedge(\eta|_T).
\]

The de Rham map converts a form into a cochain:
\begin{equation}
 (R_p\omega)(\sigma)=\int_\sigma\omega.
 \label{eq:de-rham-map}
\end{equation}
We have standard Whitney and de Rham identities:
\begin{equation}
\delta R_p=R_{p+1}d,\qquad
dW_p=W_{p+1}\delta,\qquad R_pW_p=\operatorname{id}.
\label{eq:whitney-identities}
\end{equation}
Indeed, the first one is Stokes' theorem on each simplex.  The last follows
from the integral normalization of the elementary Whitney forms. For a proof of the middle identity, see~\cite[\S~1]{Dodziuk1976}. In particular, reconstruction preserves exactness.
 We use these standard constructions as in
\cite{Whitney1957,Wilson2007}.

We can now give a geometric interpretation of the discrete wedge product.
Indeed, for $a\in C^p(K)$ and $b\in C^q(K)$, the following identity
holds~\cite{SchubelEtAl2024}:
\begin{equation}
a\wedge b=R_{p+q}(W_pa\wedge W_qb).
\label{eq:whitney-product}
\end{equation}
Here $\wedge$ on the left denotes the discrete product, and on the right
the wedge product of differential forms. When $p+q=3$, by the definition
of integration, we have
\begin{equation}
\langle a\wedge b,[K]\rangle
=\int_M W_pa\wedge W_qb.
\label{eq:whitney-integrated-product}
\end{equation}
Thus, using $dW_1a=W_2\delta a$, the Whitney helicity defined in
\eqref{eq:intro-whitney-helicity} can be expressed purely in terms of
simplicial cochains:
\[
\Hw(a)=\langle a\wedge\delta a,[K]\rangle.
\]

\subsection{Dual complex and linking}
\label{sec:dual-complex-linking}

We briefly recall Hirani's construction of dual cells~\cite{Hirani2003}.
Here we use the barycentric dual complex. Write $b_\sigma$ for the barycenter of a reference
simplex. All tetrahedra and their subdivisions below are understood in
the reference sense and are mapped elementwise to $M$ by $\Phi_T$,
yielding a simplicial complex on $M$.
The barycentric subdivision $\operatorname{sd}K$ consists of the simplices
$[b_{\sigma_0},\ldots,b_{\sigma_r}]$ associated with chains of faces
$\sigma_0\subsetneq\cdots\subsetneq\sigma_r$.
For $\sigma\in K^{(p)}$, its closed dual cell is
\[
\sigma^*:=\bigcup_{\substack{\sigma=\sigma_p\subsetneq\sigma_{p+1}\subsetneq\cdots\subsetneq\sigma_3\\
\dim\sigma_j=j}}
[b_{\sigma_p},b_{\sigma_{p+1}},\ldots,b_{\sigma_3}].
\]
The simplices in this union are called elementary dual simplices.
The associated \emph{support volume} is the three-dimensional region obtained
by taking, in each incident reference tetrahedron, the convex hull of
$\sigma$ and its local dual cell, and then forming their union.
This is the barycentric analogue of Hirani's
construction.

\begin{figure}[!htbp]
\centering
\begingroup
\definecolor{primalred}{RGB}{184,49,58}
\definecolor{dualblue}{RGB}{43,112,173}
\definecolor{dualgreen}{RGB}{31,142,112}
\definecolor{supportorange}{RGB}{220,139,36}
\begin{tikzpicture}[x=1.3cm,y=1cm,font=\small,line join=round,line cap=round]
\coordinate (A) at (-1.35,-.45);
\coordinate (B) at (2.4,1);
\coordinate (C) at (-3,2.25);
\coordinate (D) at (-.68,-2.3);
\coordinate (E) at (barycentric cs:A=1,B=1);
\coordinate (P) at (barycentric cs:A=1,B=1,C=1);
\coordinate (Q) at (barycentric cs:A=1,B=1,D=1);
\draw[gray!65,line width=.7pt] (A)--(C)--(B)--(D)--cycle;
\filldraw[fill=supportorange!25,draw=supportorange!70!black,line width=.55pt]
 (A)--(P)--(B)--cycle;
\filldraw[fill=supportorange!38,draw=supportorange!70!black,line width=.55pt]
 (A)--(Q)--(B)--cycle;
\draw[primalred,line width=1.5pt] (A)--(B);
\draw[dualblue,line width=2pt] (P)--(E);
\draw[dualgreen,line width=2pt] (E)--(Q);
\foreach \X in {A,B,C,D}{\fill[gray!75] (\X) circle (1.1pt);}
\foreach \X in {E,P,Q}{\fill (\X) circle (1.5pt);}
\path (A)--(B) node[pos=.18,text=primalred,above=2pt] {$\sigma$};
\node[below right=2pt] at (E) {$b_\sigma$};
\node[right=5pt,yshift=6pt] at (P) {$b_{\tau_+}$};
\node[right=5pt] at (Q) {$b_{\tau_-}$};
\path (E)--(P) node[midway,text=dualblue,left=4pt] {$D_+$};
\path (E)--(Q) node[midway,text=dualgreen,left=4pt] {$D_-$};
\path (P)--(C) node[pos=.52,text=gray!70!black] {$\tau_+$};
\path (Q)--(D) node[pos=.55,text=gray!70!black] {$\tau_-$};
\end{tikzpicture}
\endgroup
\caption{A 2-dimensional example. The red edge is $\sigma$; its dual
$\sigma^*=D_+\cup D_-$ is blue and green, and its support area is orange.}
\label{fig:dual-cell-support-volume}
\end{figure}
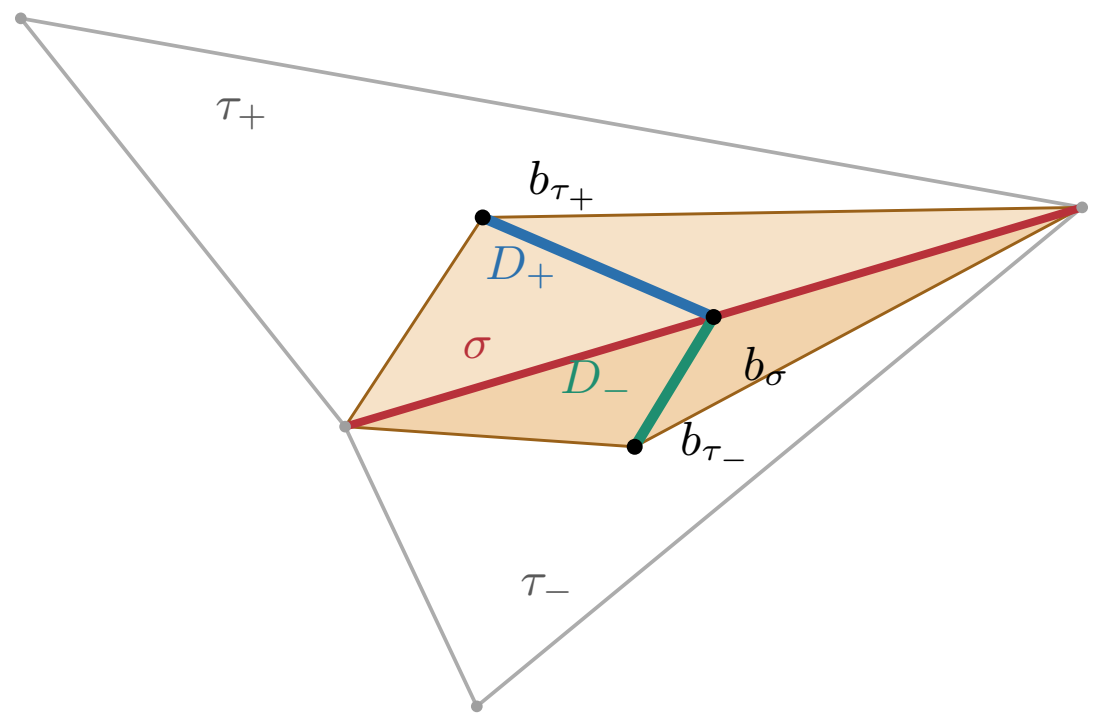

The cells $\sigma^*$ form the barycentric dual complex $K^*$, with
$\dim\sigma^*=3-p$. Their relative interiors partition $M$, and
$\tau^*$ is a face of $\sigma^*$ precisely when $\sigma$ is a face of $\tau$.
In particular, $T^*=\{b_T\}$; if a triangle $f$ separates $T_+$ and $T_-$,
then $f^*=[b_{T_+},b_f]\cup[b_f,b_{T_-}]$.

Throughout the paper, we use the term \emph{mesh} for either the primal
triangulation $K$ or its barycentric dual complex $K^*$, depending on the context.

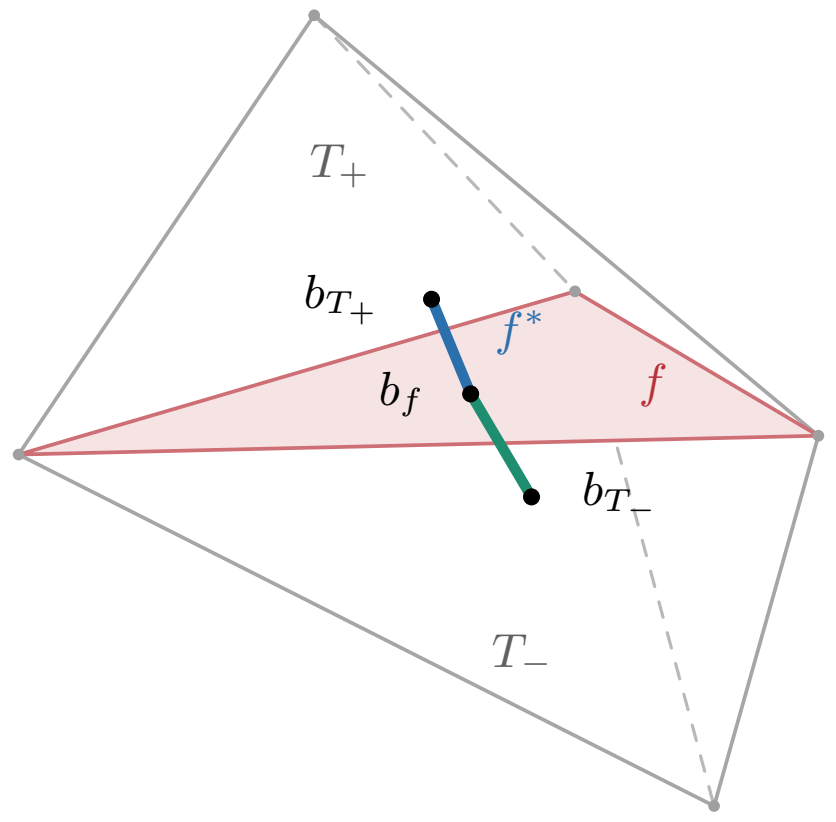
\begin{figure}[!htbp]
\centering
\begingroup
\definecolor{primalred}{RGB}{184,49,58}
\definecolor{dualblue}{RGB}{43,112,173}
\definecolor{dualgreen}{RGB}{31,142,112}
\begin{tikzpicture}[x=1.15cm,y=.83cm,font=\small,line join=round,line cap=round]
\coordinate (A) at (-2.4,0);
\coordinate (B) at (2.2,.15);
\coordinate (C) at (.8,1.3);
\coordinate (U) at (-.7,3.5);
\coordinate (V) at (1.6,-2.8);
\coordinate (F) at (.2,.48333);
\coordinate (P) at (-.025,1.2375);
\coordinate (Q) at (.55,-.3375);
\draw[gray!55,dashed,line width=.6pt] (U)--(C)--(V);
\draw[gray!70,line width=.7pt] (A)--(U)--(B) (A)--(V)--(B);
\filldraw[fill=primalred!13,draw=primalred!70,line width=.7pt]
 (A)--(B)--(C)--cycle;
\draw[dualblue,line width=2pt] (P)--(F) node[pos=.35,right=5pt,text=dualblue] {$f^*$};
\draw[dualgreen,line width=2pt] (F)--(Q);
\foreach \X in {A,B,C,U,V}{\fill[gray!75] (\X) circle (1.1pt);}
\foreach \X in {F,P,Q}{\fill (\X) circle (1.6pt);}
\node[left=6pt] at (P) {$b_{T_+}$};
\node[left=5pt] at (F) {$b_f$};
\node[right=6pt] at (Q) {$b_{T_-}$};
\node[text=primalred] at (1.25,.55) {$f$};
\node[text=gray!75!black] at (-.55,2.3) {$T_+$};
\node[text=gray!75!black] at (.5,-1.6) {$T_-$};
\end{tikzpicture}
\endgroup
\caption{Two tetrahedra sharing $f$ (red). The dual edge
$f^*=[b_{T_+},b_f]\cup[b_f,b_{T_-}]$ joins their barycenters through $b_f$.}
\label{fig:shared-face-dual-edge}
\end{figure}

Consider in the reference tetrahedra $\hat{T}$. Write $\sigma=[v_0,\ldots,v_p]$
with its vertices in increasing order, so this order gives its reference
orientation. An elementary dual simplex has the form
$D=[b_\sigma,b_{\sigma_{p+1}},\ldots,b_{\sigma_3}]$, where $\sigma_3=\hat{T}$.
Set $u_i=v_i-v_0$ for $1\le i\le p$ and
$w_j=b_{\sigma_{p+j}}-b_\sigma$ for $1\le j\le3-p$.
Choose $\varepsilon_D\in\{1,-1\}$ so that
\begin{equation}
 \varepsilon_D\,u_1\wedge\cdots\wedge u_p
 \wedge w_1\wedge\cdots\wedge w_{3-p}
 =c_D\,e_1\wedge e_2\wedge e_3,
 \qquad c_D>0,
 \label{eq:dual-orientation}
\end{equation}
where $e_1=(1,0,0)$, $e_2=(0,1,0)$, and $e_3=(0,0,1)$ are the standard
basis vectors of the reference space $\R^3$. We give $D$ the orientation
$\varepsilon_D[b_\sigma,b_{\sigma_{p+1}},\ldots,b_{\sigma_3}]$.
 Thus the vectors along
$\sigma$, followed by those along the oriented dual simplex, give the
positive orientation of the reference tetrahedron. 
This is Hirani's orientation rule~\cite{Hirani2003}.
Adjacent elementary dual simplices induce opposite boundary orientations
on their common faces, thus their orientations agree.
Their orientations therefore fit together to orient the entire dual cell
$\sigma^*$.
All dual cells and elementary dual simplices below are understood to carry
the orientations defined above. We write a cochain as a dual chain by
\begin{equation}
 \mathfrak D_p:C^p(K)\longrightarrow C_{3-p}(K^*),
 \qquad
 \mathfrak D_p a=\sum_{\sigma\in K^{(p)}}a(\sigma)\sigma^*.
 \label{eq:geometric-duality}
\end{equation}
Let $\partial:C_k(K^*)\to C_{k-1}(K^*)$ denote the cellular boundary operator.
With $[\tau:\sigma]$ denoting the incidence number of $\sigma$ in
$\partial\tau$, the induced boundary orientations give, for $0\leq p\leq2$,
\begin{equation}
 \partial\sigma^*
 =(-1)^{p+1}\sum_{\substack{\tau\supset\sigma\\\dim\tau=p+1}}
 [\tau:\sigma]\tau^*,
 \qquad
 \partial\mathfrak D_p=(-1)^{p+1}\mathfrak D_{p+1}\delta.
 \label{eq:dual-boundary}
\end{equation}
It satisfies $\partial^2=0$.

\begin{definition}\label{def:dual-cellular-coboundary}
For $0\le q\le2$, the dual exterior derivative is the cellular coboundary
$\delta_q^{\mathrm{du}}\colon C^q(K^*)\to C^{q+1}(K^*)$ given by
\[
 (\delta_q^{\mathrm{du}}\eta)(\sigma^*)=\eta(\partial\sigma^*)
\]
on each oriented dual $(q+1)$-cell $\sigma^*$.
\end{definition}

Thus $(\delta^{\mathrm{du}})^2=0$. Writing $\delta_j$ for the primal
coboundary in degree $j$, we obtain from \eqref{eq:dual-boundary} with
$p=2-q$ the matrix identity in corresponding cell bases
\begin{equation}\label{eq:dual-coboundary-sign}
 \delta_q^{\mathrm{du}}=(-1)^{q+1}\delta_{2-q}^{\top},\qquad 0\le q\le2.
\end{equation}

The dual complex $K^*$ is generally a cell complex rather than a
simplicial complex, whereas the elementary dual simplices and all their
faces form $\operatorname{sd}K$, a common simplicial subdivision of $K$
and $K^*$ whose tetrahedra inherit the orientation of $M$.
The 1-skeletons of $K$ and $K^*$ are finite simple graphs, denoted by
$G_{\mathrm{pr}}$ and $G_{\mathrm{du}}$, respectively. They are disjoint, since dual
edges lie in tetrahedron interiors except where they pass through
barycenters in the relative interiors of triangular faces.

We next introduce the basic ingredients needed for linking numbers,
starting with the smooth picture. Let $S_1,S_2\subset M$ be compact
oriented smooth submanifolds of complementary dimensions, transverse in
their interiors, with $\partial S_1\cap S_2=S_1\cap\partial S_2=\varnothing$.
Their intersection number is
\[
 I(S_1,S_2):=\sum_{x\in S_1\cap S_2}\varepsilon_x.
\]
At each intersection point $x$, let $p=\dim S_1$ and choose bases
$(u_1,\ldots,u_p)$ of $T_xS_1$ and
$(v_1,\ldots,v_{3-p})$ of $T_xS_2$ that agree with the given
orientations of $S_1$ and $S_2$, respectively. By transversality, the
ordered list
\[
 (u_1,\ldots,u_p,v_1,\ldots,v_{3-p})
\]
is a basis of $T_xM$. We set $\varepsilon_x=1$ if this ordered basis
agrees with the orientation of $M$, and $\varepsilon_x=-1$ if it has
the opposite orientation. 

Since $M$ is simply connected, every oriented embedded piecewise smooth
closed curve in $M$ is the boundary of some compact oriented piecewise
smooth surface. We define the linking number:
\begin{definition}
\label{def:real-linking}
Let $C_1,C_2\subset M$ be disjoint oriented embedded piecewise smooth
closed curves. Choose a compact oriented piecewise smooth surface $S$
with $\partial S=C_1$ and arrange that all intersections of $S$ with
$C_2$ are transverse and occur at smooth points of both. Define
\begin{equation}
 \Lk(C_1,C_2):=I(S,C_2).
 \label{eq:linking-definition}
\end{equation}
\end{definition}

To see that this is well defined, let $S'$ be another such surface.
Gluing $S$ to $-S'$ along their common boundary gives a closed oriented
surface  in $M$. Since $C_2$ is null-homotopic, its algebraic
intersection with this closed surface vanishes. Thus
\[
 I(S,C_2)-I(S',C_2)=0.
\]
The linking number is therefore independent of the choice of $S$ and
is an integer. Its symmetry follows from the description of linking
numbers in terms of Thom forms and Stokes' theorem~\cite[\S~17]{BottTu1982}:
\begin{equation}
 \Lk(C_1,C_2)=\Lk(C_2,C_1).
 \label{eq:linking-symmetry}
\end{equation}
The following figure illustrates the definition.

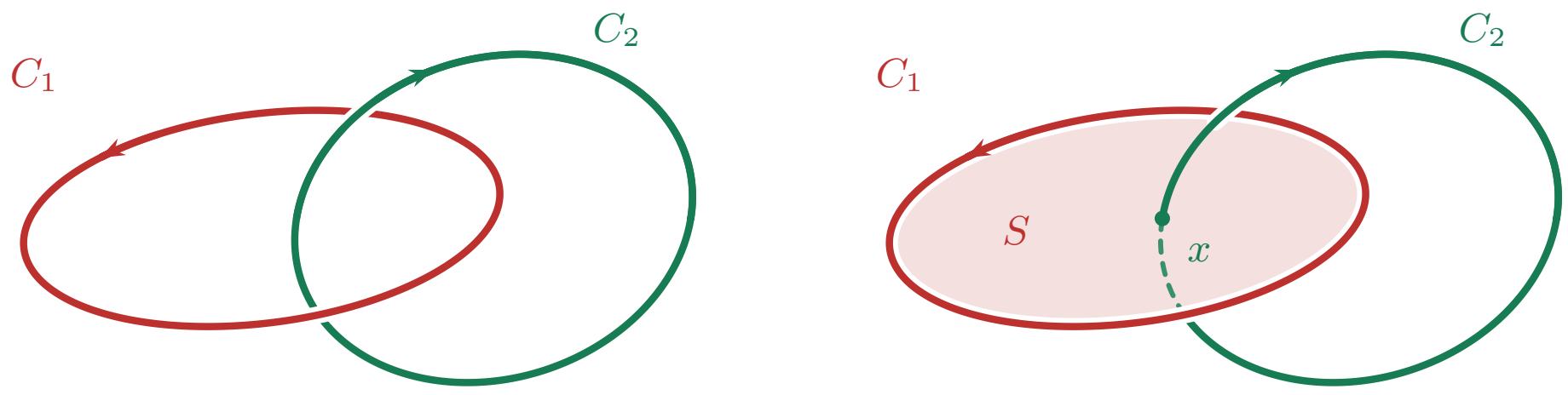
\begin{figure}[htbp]
\centering
\begingroup
\definecolor{linkred}{RGB}{188,48,46}
\definecolor{linkgreen}{RGB}{23,124,83}
\begin{tikzpicture}[x=1.2cm,y=1.2cm,font=\small,
 line cap=round,line join=round]
\foreach \panel in {0,1}{
\begin{scope}[shift={({5.6*\panel},0)}]
% The same pair of linked circles is used in both panels.
\draw[linkgreen,line width=1.6pt]
 plot[domain=0:360,samples=181,smooth,variable=\t]
 ({1.5*(1+.85*cos(\t))+.17*sin(\t)},{1.0625*sin(\t)});
\ifnum\panel=1\relax
 \fill[linkred!15]
  plot[domain=0:360,samples=181,smooth,variable=\t]
  ({1.5*cos(\t)+.35*sin(\t)},{.7*sin(\t)}) -- cycle;
 % A dashed arc indicates the part behind the spanning disk.
 \begin{scope}
  \clip plot[domain=0:360,samples=181,smooth,variable=\t]
   ({1.5*cos(\t)+.35*sin(\t)},{.7*sin(\t)}) -- cycle;
  \draw[linkgreen!85,line width=1.1pt,dash pattern=on 2.5pt off 2.5pt]
   plot[domain=180:360,samples=91,smooth,variable=\t]
   ({1.5*(1+.85*cos(\t))+.17*sin(\t)},{1.0625*sin(\t)});
 \end{scope}
\fi
\draw[linkred,line width=1.6pt,preaction={draw=white,line width=3.8pt}]
 plot[domain=0:360,samples=181,smooth,variable=\t]
 ({1.5*cos(\t)+.35*sin(\t)},{.7*sin(\t)}) -- cycle;
% The upper green branch passes in front of the red circle.
\draw[white,line width=3.8pt]
 plot[domain=136:144,samples=13,smooth,variable=\t]
 ({1.5*(1+.85*cos(\t))+.17*sin(\t)},{1.0625*sin(\t)});
\draw[linkgreen,line width=1.6pt]
 plot[domain=0:180,samples=91,smooth,variable=\t]
 ({1.5*(1+.85*cos(\t))+.17*sin(\t)},{1.0625*sin(\t)});
\draw[linkred,line width=1.6pt,-{Stealth[length=1.9mm,width=1.4mm]}]
 plot[domain=129:146,samples=13,smooth,variable=\t]
 ({1.5*cos(\t)+.35*sin(\t)},{.7*sin(\t)});
\draw[linkgreen,line width=1.6pt,-{Stealth[length=1.9mm,width=1.4mm]}]
 plot[domain=132:115,samples=13,smooth,variable=\t]
 ({1.5*(1+.85*cos(\t))+.17*sin(\t)},{1.0625*sin(\t)});
\node[text=linkred] at (-1.47,.93) {$C_1$};
\node[text=linkgreen] at (2.3,1.22) {$C_2$};
\ifnum\panel=1\relax
 \node[text=linkred] at (-.72,-.08) {$S$};
 \fill[linkgreen] (.225,0) circle (1.7pt);
 \node[text=linkgreen,below right=2pt] at (.225,0) {$x$};
 \node at (.65,-1.48) {(b) A spanning surface for $C_1$};
\else
 \node at (.65,-1.48) {(a) Linked curves};
\fi
\end{scope}
}
\end{tikzpicture}
\endgroup
\caption{The linked curves $C_1$ (red) and $C_2$ (green).
The indicated orientations give $\Lk(C_1,C_2)=1$.}
\label{fig:linking-definition}
\end{figure}

We then define linking directly for real 1-cycles $c$ and $d$ on the
primal and dual 1-skeletons, respectively. Since $M$ is simply connected,
choose a real 2-chain $S\in C_2(K)$ with $\partial S=c$.
By definition, $S$ and $d$ intersect transversely. We set
\[
 \Lk(c,d):=I(S,d),
\]
The same argument as above shows that the linking number is well defined.

Finally, fix a Gauss--de Rham linking form
$G\in\Omega^2((M\times M)\setminus\Delta)$, where
$\Delta:=\{(x,x):x\in M\}$. For any two disjoint oriented closed curves $C_1,C_2\subset M$, we have
\[
  \int_{C_1\times C_2} G=\Lk(C_1,C_2).
\]
For its existence,
see \cite{ArnoldKhesin2021} and
\cite{Vogel2003}. It defines a bilinear form
\begin{equation}
 \Lk_G:C_1(K)\times C_1(K^*)\longrightarrow\R,
 \qquad
 \Lk_G(c,d):=\int_{c\times d}G.
 \label{eq:linking-form-pairing}
\end{equation}
Since the primal and dual 1-skeletons are compact and disjoint,
$c\times d$ stays away from $\Delta$, so this integral is well defined.
On cycles, $\Lk_G$ agrees with the linking number defined above, thus it can be viewed as an extension of the above linking pairing. However, 
on general chains, it depends on the fixed choice of $G$.

We fix $G$ for the rest of the paper and suppress the subscript,
writing $\Lk$ for the resulting bilinear form on
$C_1(K)\times C_1(K^*)$.
On the finite-dimensional spaces of primal and dual cycles this pairing is
bounded in the coefficient norm $\|\sum_e z_e e\|_{\ell^1}=\sum_e|z_e|$.
Thus for some $C_{\mathrm{lk}}<\infty$ depending only on the triangulation,
\begin{equation}\label{eq:markov-linking-bound}
 |\Lk(z,w)|\le C_{\mathrm{lk}}\,\|z\|_{\ell^1}\,\|w\|_{\ell^1}
\end{equation}
for every primal cycle $z$ and dual cycle $w$.

\begin{proposition}
\label{prop:intersection-evaluation}
For any $S\in C_2(K)$ and $\eta\in C^2(K)$ with $\delta\eta=0$,
\[
 \Lk(\partial S,\mathfrak D_2\eta)=\langle\eta,S\rangle.
\]
\end{proposition}

\begin{proof}
The dual boundary formula gives
$\partial\mathfrak D_2\eta=-\mathfrak D_3\delta\eta=0$.
Thus $\partial S$ and $\mathfrak D_2\eta$ are 1-cycles supported on
the disjoint primal and dual 1-skeletons. Since $M$ is simply connected,
both are boundaries, so their linking number is defined.

Let $f,g\in K^{(2)}$. The dual edge $g^*$ consists of the two segments
joining $b_g$ to the barycenters of the two tetrahedra incident to $g$.
Apart from $b_g$, these segments lie in the interiors of those
tetrahedra. Hence $f\cap g^*$ is empty for $f\neq g$, while
$f\cap f^*=\{b_f\}$. The two branches of $f^*$ lie on opposite sides
of $f$. A small local smoothing therefore gives a single transverse
intersection, whose sign is positive by
\eqref{eq:dual-orientation}. Consequently,
$I(f,g^*)=\delta_{fg}$.

Writing $S=\sum_f s_f f$ and using bilinearity, we obtain
\begin{equation}
 \begin{aligned}
 I(S,\mathfrak D_2\eta)
 &=\sum_{f,g}s_f\eta(g)I(f,g^*)\\
 &=\sum_f s_f\eta(f)=\langle\eta,S\rangle.
 \end{aligned}
 \label{eq:intersection-evaluation}
\end{equation}
Since $S$ bounds $\partial S$, the definition of linking gives
$\Lk(\partial S,\mathfrak D_2\eta)=I(S,\mathfrak D_2\eta)$,
which proves the claim.
\end{proof}

\subsection{Helicity in hydrodynamics}

Helicity is a classical invariant in fluid mechanics and
magnetohydrodynamics that relates a divergence-free field to the
topology of its field lines. In 1969, Moffatt showed that the helicity
of linked vortex filaments is determined by their linking numbers
and circulations~\cite{Moffatt1969}. This interpretation connects an
integral defined from local field data with the global topology of
the field lines.

For a general divergence-free field, trajectories need not return to
their starting points, and the familiar picture of linked rings no
longer applies directly. Arnold's insight was to recover linking in the long-time limit by
closing trajectory segments with short return paths and dividing
their linking number by the product of the observation times. His theorem identifies the volume
average of the resulting limit with helicity~\cite{ArnoldKhesin2021}.
This passage from local differential forms to the long-time geometry
of trajectories is the model for our discrete construction. We recall
the analytic definition before making this dynamical picture precise.

For the flow $\varphi_t$ of $X$, Cartan's formula and the definition of the
Lie derivative give the equivalences
\[
 \mathcal L_X\mu=0
 \quad\Longleftrightarrow\quad d(\iota_X\mu)=0
 \quad\Longleftrightarrow\quad \varphi_t^*\mu=\mu\quad(t\in\R).
\]
Since $M$ is simply connected, Poincar\'e duality and de Rham's theorem
imply that the closed 2-form $\iota_X\mu$ is exact~\cite{Hatcher2002}. We may therefore choose a smooth 1-form
$\alpha$ satisfying $d\alpha=\iota_X\mu$ and call $\alpha$ a
potential for $X$.

\begin{definition}
\label{def:smooth-helicity}
The helicity of $X$ is
\begin{equation}
 \Hsm(X)=\Hsm(\alpha):=\int_M\alpha\wedge d\alpha.
 \label{eq:smooth-helicity}
\end{equation}
\end{definition}

Helicity is independent of the choice of potential, a property known
as \emph{gauge invariance}. Indeed, any other potential has the form
$\alpha+\eta$ with $d\eta=0$, and
\[
 \Hsm(\alpha+\eta)-\Hsm(\alpha)
 =\int_M\eta\wedge d\alpha
 =\int_M d(\alpha\wedge\eta)=0
\]
by Stokes' theorem. Helicity is also independent
of the metric and invariant under volume-preserving
diffeomorphisms~\cite{ArnoldKhesin2021}. We retain the notation
$\Hsm(\alpha)$ when comparing smooth and discrete potentials.

We now turn to Arnold's helicity theorem, which identifies helicity
with the average asymptotic linking of trajectories. Since $M$ is compact, $\varphi_t$ is defined for all $t\in\R$. Choose a system of oriented minimizing
geodesics $\sigma(p,q)$ from $p$ to $q$ as in Vogel's
construction~\cite{Vogel2003}. The paths depend measurably on their
endpoints and have length at most $\operatorname{diam}(M)$;
$\sigma(p,p)$ is the constant path. For $T>0$, write
$\gamma_T(x):[0,T]\to M$, $\gamma_T(x)(t)=\varphi_t(x)$, for the
oriented trajectory segment, and close it by setting
\begin{equation}
 \Gamma_T(x):=\gamma_T(x)*\sigma(\varphi_T(x),x),
 \qquad
 \Gamma_S(y):=\gamma_S(y)*\sigma(\varphi_S(y),y),
 \label{eq:smooth-closed-trajectories}
\end{equation}
where $*$ denotes concatenation of paths.  For each fixed $T,S>0$, the chosen system makes them
disjoint for $(\mu\otimes\mu)$-almost every $(x,y)$, so their linking
number is defined a.s.

\begin{definition}
\label{def:smooth-asymptotic-linking}
Define $F_{T,S}\in L^1(M\times M,\mu\otimes\mu)$ by
\[
 F_{T,S}(x,y):=
 \frac{\Lk(\Gamma_T(x),\Gamma_S(y))}{TS}.
\] When the limit exists,
the asymptotic linking function of $X$ is
\begin{equation}
 \lambda_X:=\lim_{T,S\to\infty}F_{T,S}
 \quad\text{in }L^1(M\times M,\mu\otimes\mu).
 \label{eq:smooth-asymptotic-linking}
\end{equation}
Its average asymptotic linking number is
\[
 \Lambda(X):=\int_{M\times M}\lambda_X(x,y)\,\mu(x)\mu(y).
\]
\end{definition}

The analytic and dynamical descriptions meet in Arnold's helicity
theorem. We use the $L^1$ formulation established by
Vogel~\cite{ArnoldKhesin2021,Vogel2003}.

\begin{theorem}
\label{thm:smooth-helicity}
The limit
\eqref{eq:smooth-asymptotic-linking} exists as $T$ and $S$ tend to
infinity independently and is independent of the chosen geodesic
system above. Moreover,
\begin{equation}
 \Lambda(X)
 =\Hsm(X)
 \label{eq:smooth-helicity-theorem}
\end{equation}
\end{theorem}

Thus helicity records the average long-term entanglement of a vector
field's trajectories throughout the manifold, connecting its description in terms of differential forms with its
global topological properties.

While much of the existing work on discrete helicity has focused on
computational aspects~\cite{Kotiuga1990,ChammasKotiuga1996}, this
intrinsic topological characterization has received less attention,
and our aim is to establish its analogue in the discrete setting.

\subsection{Continuous-time Markov chains on a finite set}
\label{sec:finite-markov-chains}

We assume familiarity with the basics of stochastic processes and refer
the reader to Durrett~\cite{Durrett2016} for an introduction. We briefly
recall the facts about continuous-time Markov chains used below.

Let $\mathcal V$ be a finite set, called the \emph{state space}, its elements are called \emph{states}.  In \S\ref{sec:markov}, these states
will be the vertices of the primal or dual 1-skeleton.
A time-homogeneous continuous-time Markov chain on $\mathcal V$ is a jump process $(X_t)_{t\geq 0}$ with natural filtration $\mathcal F_t=\sigma(X_r:0\leq r\leq t)$.
Thus $\mathcal F_t$ records all the information obtained by observing the chain up to and including time $t$.
The Markov property requires that, for every $s,t\geq0$ and $y\in\mathcal V$,
\[
 \mathbb P(X_{t+s}=y\mid\mathcal F_t)
 =\mathbb P(X_{t+s}=y\mid X_t)\quad\text{a.s.}
\]
 For $s\geq 0$, the
\emph{transition matrix} at time $s$ is the matrix indexed by $\mathcal V\times\mathcal V$ whose entries are
\[
 P_s=\bigl(P_s(x,y)\bigr)_{x,y\in\mathcal V},\qquad
 P_s(x,y)=\mathbb P_x(X_s=y).
\]
The row indexed by $x$ gives the distribution of the chain after time $s$
when it starts at $x$. Its entries are nonnegative and sum to one.  On a function $f\colon\mathcal V\to\R$, the matrix
$P_s$ acts by taking the expected value of $f$ after time $s$,
\[
 (P_s f)(x)=\sum_{y\in\mathcal V}P_s(x,y)f(y)=\mathbb E_x[f(X_s)].
\]
More generally, by time homogeneity, $\mathbb E[f(X_{t+s})\mid X_t]=(P_s f)(X_t)$.

A state $x$ is \emph{absorbing} if the chain stays at $x$ a.s.\ once it
reaches $x$, equivalently if $P_t(x,x)=1$ for every $t\geq0$.
At a nonabsorbing state $x$, the chain waits an exponential time with
parameter $\ell_x>0$. It then jumps to $y\ne x$ with probability $p_{xy}$,
independently of the waiting time~\cite{Norris1997}. We set $\ell_x=0$
at an absorbing state.

The transition matrices satisfy $P_{t+s}=P_tP_s$ and $P_0=I$. Their
\emph{ generator} is the right derivative at time zero,
\[
 Q=\left.\frac{d}{dt}P_t\right|_{t=0}.
\]
Write $Q=(q_{xy})_{x,y\in\mathcal V}$. For $x\ne y$, we have
$P_0(x,y)=0$, so $q_{xy}=\lim_{h\downarrow0}P_h(x,y)/h$.
Starting at a nonabsorbing state $x$, the probability that the first
jump occurs by time $h$ and lands at $y$ is $(1-e^{-\ell_xh})p_{xy}$.
With $\Lambda=\max_{z\in\mathcal V}\ell_z$, the probability of at least
two jumps by time $h$ is at most $(\Lambda h)^2$, since each jump has
conditional probability at most $\Lambda h$ within the time remaining.
Thus
\[
 P_h(x,y)=(1-e^{-\ell_xh})p_{xy}+O(h^2),\qquad x\ne y.
\]
It follows that $q_{xy}=\ell_xp_{xy}$, the jump rate from $x$ to $y$.
At an absorbing state $x$, we have $q_{xy}=0$ for $y\ne x$.

Since each row of $P_h$ sums to one, each row of $Q$ sums to zero. Thus
$q_{xx}=-\sum_{y\ne x}q_{xy}=-\ell_x$ and
\[
 (Qf)(x)=\sum_{y\ne x}q_{xy}\bigl(f(y)-f(x)\bigr).
\]
On a finite set the generator determines the transition matrices by $P_t=e^{tQ}$.
A \emph{path of jumps with positive rates} from $x$ to $y$ is a finite sequence
$x=x_0,x_1,\ldots,x_m=y$ with $q_{x_i x_{i+1}}>0$ for $0\leq i<m$.
The chain is \emph{irreducible} if such a path exists for every ordered
pair of distinct states $x,y\in\mathcal V$.

A nonempty set $C\subseteq\mathcal V$ is \emph{closed} if $q_{xy}=0$
for every $x\in C$ and $y\notin C$. It is a \emph{closed irreducible class}
if, in addition, every ordered pair of distinct states in $C$ can be
joined by a path of jumps with positive rates lying entirely in $C$.

A probability distribution $\pi$ is \emph{stationary} if a chain started
with distribution $\pi$ has that same distribution at every later time.
Viewing $\pi$ as a row vector, this means $\pi P_t=\pi$ for every $t\geq0$,
or equivalently $\pi Q=0$. In coordinates this is the balance of the
incoming and outgoing probability fluxes at each state,
\[
 \sum_{y\ne x}\pi(y)q_{yx}=\pi(x)\ell_x.
\]
In stationarity, the expected number of jumps from $x$ to $y$ in time $t$
is $t\pi(x)q_{xy}$ for $x\ne y$. The chain is \emph{reversible} with
respect to $\pi$ if the two fluxes agree on each pair of states, that is,
$\pi(x)q_{xy}=\pi(y)q_{yx}$. For these conventions on finite sets, see
Norris~\cite{Norris1997}.

We use the following facts to estimate the number of jumps along an edge in
\S\ref{sec:markov}.

\begin{lemma}
\label{lem:markov-class-poisson}
Let $Q$ be the generator of a continuous-time Markov chain on the finite
set $\mathcal V$, with finite jump rates.
\begin{enumerate}[label=(\roman*),ref=\roman*,beginpenalty=10000,midpenalty=10000]
\item\label{item:markov-nonexplosion}
There are a.s.\ only finitely many jumps on any bounded time
interval.
\item\label{item:markov-poisson}
Let $C$ be a closed irreducible class with stationary distribution
$\pi_C$. For $f\colon C\to\R$, the Poisson equation $Qg=f$ on $C$
has a solution if and only if
$\sum_{x\in C}\pi_C(x)f(x)=0$.
\item\label{item:markov-closed-classes}
Suppose $\pi(x)>0$ at every state and $\pi Q=0$. Join $x$ and $y$
whenever $q_{xy}>0$ or $q_{yx}>0$. Each connected component $C$ is a
closed irreducible class, with stationary distribution
$\pi_C=\pi|_C/\pi(C)$, where $\pi(C):=\sum_{x\in C}\pi(x)$ is the
total mass of $C$.
\end{enumerate}
\end{lemma}

\begin{proof}
Part~(\ref{item:markov-nonexplosion}) follows from
\cite[\S~2.7]{Norris1997}. For part~(\ref{item:markov-poisson}),
apply the solvability criterion and semigroup construction in
\cite[\S~III]{KhodabandehlouEtAl2024} to the chain restricted to $C$.

For part~(\ref{item:markov-closed-classes}), use the decomposition of
stationary distributions in \cite[\S~2.4]{KuntzEtAl2021}.
Strict positivity of $\pi$ implies that every state belongs to a closed
irreducible class. Since no jumps connect different such classes,
they are precisely the connected components in the statement.
Restricting $\pi Q=0$ to $C$ and normalizing gives
$\pi_C=\pi|_C/\pi(C)$.
\end{proof}

For a directed jump $e:e^-\to e^+$, write $e^-$ and $e^+$ for its initial
and final states and $q_e$ for its rate. Let $N_e(t)$ denote the number of jumps along the directed edge $e$
during the time interval $[0,t]$.

\begin{lemma}
\label{lem:markov-jump-martingales}
For every initial state and every $g\colon\mathcal V\to\R$, the processes
\[
 \begin{aligned}
 A_e(t)&=N_e(t)-q_e\int_0^t\mathbf 1_{\{X_{s-}=e^-\}}\,ds,\\
 D_g(t)&=g(X_t)-g(X_0)-\int_0^t (Qg)(X_s)\,ds
 \end{aligned}
\]
are square-integrable martingales . Their
second moments satisfy
\[
 \mathbb E|A_e(t)|^2\le q_e t,\qquad
 \mathbb E|D_g(t)|^2\le t\max_x\sum_{y\ne x}q_{xy}\bigl(g(y)-g(x)\bigr)^2.
\]
\end{lemma}

\begin{proof}
For the compensated jump count $A_e$, see \cite[\S~8.3]{DarlingNorris2008}
and \cite[\S~3.3]{PangTalrejaWhitt2007}. For Dynkin's martingale $D_g$ and
its second moment, see \cite[\S~8.4]{DarlingNorris2008}.
\end{proof}

\begin{definition}
For square-integrable random variables $Y_t$ and $Y$ on a common probability
space, with values in a finite-dimensional normed space, we say that $Y_t$
\emph{converges strongly} to $Y$ as $t\to\infty$ if
\[
 Y_t\to Y\quad\text{a.s.},\qquad
 \lim_{t\to\infty}\mathbb E\|Y_t-Y\|^2=0.
\]
In the stochastic statements below, \emph{strong convergence} always refers
to this joint mode of convergence.
\end{definition}

\begin{lemma}
\label{lem:markov-martingale-law}
If $D$ is a right continuous square-integrable martingale with $D(0)=0$
and $\mathbb E|D(t)|^2\le Kt$, then $D(t)/t\to0$ strongly
as $t\to\infty$ .
\end{lemma}

\begin{proof}
The bound on the second moment makes
$Z_t=\int_0^t(1+s)^{-1}\,dD(s)$ an $L^2$-bounded martingale.
Applying the martingale convergence argument in \cite[\S~1]{Lepingle1978}
with normalizing process $1+t$ gives $D(t)/t\to0$ a.s. The $L^2$ estimate is
$\mathbb E|D(t)/t|^2\le K/t$.
\end{proof}

We next use a spanning tree to close walks on a discrete 1-skeleton,
just as we used geodesic return paths to close smooth trajectory segments.
To close a walk in a finite connected graph, fix a rooted spanning tree
and let $\tau(v)$ be the oriented tree path from the root to $v$.
On each edge set $\Pi e=e+\tau(e^-)-\tau(e^+)$ and extend linearly.
Collecting the coefficients at each vertex gives
\begin{equation}\label{eq:markov-cycle-projection}
 \Pi z=z-\sum_v(\partial z)_v\,\tau(v).
\end{equation}
\begin{figure}[H]
\centering
\begin{tikzpicture}[x=1.08cm,y=1.1cm,line cap=round,line join=round,
  >={Stealth[length=2.6mm,width=1.8mm]},
  tree/.style={draw=black!35,line width=.8pt},
  route/.style={draw=orange!80!black,line width=1.55pt},
  shared/.style={route,dash pattern=on 3pt off 2pt},
  chord/.style={draw=blue!68!black,line width=1.6pt}]
\def\treeArrow#1#2#3{%
  \draw[#1] (#2) -- (#3);
  \path (#2) -- coordinate[pos=.58] (arrowtip) (#3);
  \draw[#1,->] (#2) -- (arrowtip);}
\foreach \panel/\offset in {1/0cm,2/7cm}{
  \begin{scope}[xshift=\offset]
    \coordinate (o) at (2.8,3.36);
    \coordinate (w) at (1.65,2.24); \coordinate (c) at (4.35,2.24);
    \coordinate (a) at (.62,1.12);  \coordinate (b) at (2.80,1.12);
    \coordinate (d) at (5.0,1.12);
    \coordinate (u) at (.15,0);     \coordinate (v) at (3.25,0);
    \draw[tree] (u)--(a)--(w)--(o)--(c)--(d) (w)--(b)--(v);
    \ifnum\panel=1
      \treeArrow{shared}{o}{w}
      \treeArrow{route}{w}{a} \treeArrow{route}{a}{u}
      \treeArrow{route}{w}{b} \treeArrow{route}{b}{v}
      \path (a)--(u) node[midway,left=5pt] {$\tau(u)$};
      \path (b)--(v) node[midway,right=5pt] {$\tau(v)$};
      \node at (2.45,-.84) {(a) Root paths to $u$ and $v$};
    \else
      \treeArrow{route}{v}{b} \treeArrow{route}{b}{w}
      \treeArrow{route}{w}{a} \treeArrow{route}{a}{u}
      \node at (2.45,-.84) {(b) The closed cycle $\Pi e$};
    \fi
    \treeArrow{chord}{u}{v}
    \path (u)--(v) node[pos=.46,above=3pt] {$e$};
    \foreach \pt in {o,w,a,b,c,d,u,v}{\fill (\pt) circle (2.1pt);}
    \node[above=3pt] at (o) {$o$}; \node[left=3pt] at (w) {$w$};
    \node[below=3pt] at (u) {$u$}; \node[below=3pt] at (v) {$v$};
  \end{scope}}
\end{tikzpicture}
\caption{Closing the blue edge $e=[u,v]$ in a rooted spanning tree.
The root paths share the dashed edge $[o,w]$, which cancels in
$\tau(u)-\tau(v)$. The orange arrows in (b) trace the return path
from $v$ through $w$ to $u$, completing $\Pi e=e+\tau(u)-\tau(v)$.}
\label{fig:tree-closure}
\end{figure}
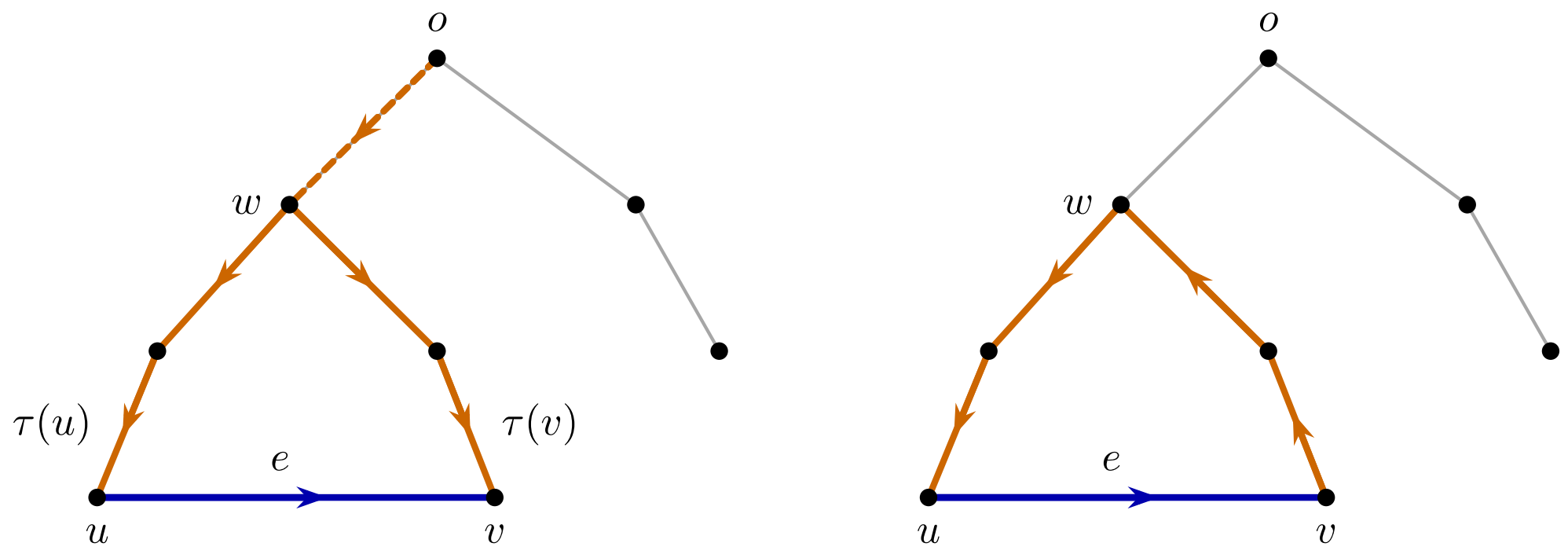
Since $\partial\tau(v)=[v]-[o]$ for the root $o$, $\Pi z$ is a cycle and
$\Pi$ fixes every cycle. The correction for a walk from $u$ to $v$ is the
unique tree path from $v$ to $u$. Its coefficient norm is at most $|C|-1$
on a component with vertex set $C$.

\section{Discrete helicity and a first look at linking numbers}
\label{sec:discrete-helicity}

We introduce cup helicity alongside Whitney helicity and use cup--cap
duality to express the former as the linking number of a primal 1-cycle
and a dual 1-cycle. Comparing this formula with the helicity of linked
flux tubes suggests interpreting the cycle weights as fluxes through
cross sections.

\subsection{Cup and Whitney helicity}

To construct a discrete helicity, we need to discretize the 1-form,
the product, and integration over $M$. We first seek a discrete
counterpart of a 1-form. The discussion of the de Rham and Whitney
maps above suggests that a 1-cochain is an appropriate object.
Throughout this section, we therefore fix $a\in C^1(K)$ and set
$b=\delta a$, so that $\delta b=0$.

For the product, the discrete wedge is a natural choice, but it does
not directly give the required duality on chains. Motivated by
its relation to the cup product and by cup--cap duality, we also
introduce an auxiliary cup helicity and study its properties.

\begin{definition}
\label{def:discrete-helicities}
The  cup helicity and Whitney helicity of $a$ are
\begin{align}
 \Hc(a)&:=\langle a\smile\delta a,[K]\rangle,
 \label{eq:cup-helicity}\\
 \Hw(a)&:=\int_MW_1a\wedge W_2\delta a
 =\langle a\wedge \delta a,[K]\rangle.
 \label{eq:whitney-helicity}
\end{align}
\end{definition}

The following proposition
shows that both discrete helicities also possess the gauge invariance property.

\begin{proposition}
\label{prop:gauge-invariance}
For every $z\in Z^1(K)$,
\begin{equation}
 H_\star(a+z)=H_\star(a),
 \qquad \star\in\{\mathrm c,\mathrm w\}.
 \label{eq:discrete-gauge-invariance}
\end{equation}
In particular, each helicity depends only on the exact 2-cochain
$\delta a$.
\end{proposition}
\begin{proof}
For either product $\diamond\in\{\smile,\wedge\}$, the Leibniz identity
and $\delta z=0$ give
\[
z\diamond\delta a=-\delta(z\diamond a).
\]
Its pairing with $[K]$ is zero by \eqref{eq:fundamental-cycle}.
Since $\delta(a+z)=\delta a$, this proves gauge invariance of both
helicities.
\end{proof}

\begin{remark}
\label{rem:vertex-order}
The cup product in $\Hc(a)$ uses the vertex order. Reordering the vertices
leaves Whitney helicity $\Hw(a)$ unchanged~\cite{SchubelEtAl2024}.
\end{remark}
We then provide a discrete analogue of the Arnold's energy inequality.

The Whitney mass matrix $M_2$ represents the Galerkin Hodge star on 2-cochains and is assembled from local mass matrices as in Bell and Hirani~\cite[\S~9]{BellHirani2012}.
Its quadratic form evaluates the energy integral of the Whitney reconstruction $W_2b$ exactly and thus gives a discrete approximation to kinetic energy.

\begin{proposition}
\label{cor:whitney-energy}
For every $a\in C^1(K)$, let $\mathbf b=(b(f))_{f\in K^{(2)}}$ be the face coefficients of $b=\delta a$.
Then the discrete kinetic energy satisfies
\begin{equation}
 E_K(b):=\tfrac12\mathbf b^{\mathsf T}M_2\mathbf b
 =\tfrac12\int_M W_2b\wedge *W_2b
 =\tfrac12\|W_2b\|_{L^2(M)}^2
 \ge C_M\,|\Hw(a)|,
 \label{eq:whitney-energy}
\end{equation}
where $C_M>0$ depends only on $M$.
\end{proposition}
\begin{proof}
The first three equalities follow from the definition of the Whitney mass matrix~\cite[\S~9]{BellHirani2012}.
The inequality follows from Arnold's continuous energy bound~\cite{ArnoldKhesin2021} applied to the potential $W_1a$, with $d(W_1a)=W_2b$ by \eqref{eq:whitney-identities}.
\end{proof}

\subsection{Cup helicity as primal--dual linking}

The cap product and the map to dual cells assign to $b=\delta a$
the primal and dual 1-chains
\begin{equation}
 C_{\mathrm{pr}}^{\mathrm c}(a):=[K]\frown b \in C_1(K),\qquad
 C_{\mathrm{du}}^{\mathrm c}(a):=\mathfrak D_2 b \in C_1(K^*).
 \label{eq:cup-chains}
\end{equation}

\begin{theorem}
\label{thm:cup-linking}
For every $a\in C^1(K)$, the 1-chains in \eqref{eq:cup-chains} are
boundaries with explicit fillings
\begin{equation}
 C_{\mathrm{pr}}^{\mathrm c}(a)=\partial([K]\frown a),\qquad
 C_{\mathrm{du}}^{\mathrm c}(a)=\partial(\mathfrak D_1a).
 \label{eq:cup-chains-exact}
\end{equation}
Their linking number is the cup helicity of $a$,
\begin{equation}
 \Hc(a)=\Lk\bigl(C_{\mathrm{pr}}^{\mathrm c}(a),
                 C_{\mathrm{du}}^{\mathrm c}(a)\bigr).
 \label{eq:cup-linking}
\end{equation}
\end{theorem}
\begin{proof}
Take the primal 2-chain $S=[K]\frown a$. Since $a$ has degree $1$
and $\partial[K]=0$, the cap boundary formula \eqref{eq:cap-boundary}
gives
\[
 \partial S=-(\partial[K])\frown a+[K]\frown\delta a
           =[K]\frown b.
\]
In the same degree, the sign in the dual boundary formula
\eqref{eq:dual-boundary} is $(-1)^{1+1}=1$, so
$\partial(\mathfrak D_1a)=\mathfrak D_2\delta a=\mathfrak D_2b$.
This proves \eqref{eq:cup-chains-exact}. The two boundaries lie on the
disjoint primal and dual 1-skeletons, so Definition~\ref{def:real-linking}
computes their linking number using the filling $S$.

Since $\delta b=0$, Proposition~\ref{prop:intersection-evaluation} and
cup--cap duality \eqref{eq:cup-cap-adjunction} give
\[
 \Lk\bigl(C_{\mathrm{pr}}^{\mathrm c}(a),C_{\mathrm{du}}^{\mathrm c}(a)\bigr)
 =\langle b,[K]\frown a\rangle
 =\langle a\smile b,[K]\rangle
 =\Hc(a).
\]

\end{proof}

The theorem above gives cup helicity a linking interpretation,
but leaves open why these particular 1-chains should represent the field
and how a deterministic linking identity can lead to an average over long
trajectories. We explore the first question through the flux tube model
in the next subsection and the geometric construction in
\S\ref{sec:whitney-flux}. \S\ref{sec:markov} answers the
second for Whitney helicity through a discrete analogue of Arnold's
helicity theorem.

\subsection{Cycle weights and the flux tube model}

A directed cycle in either 1-skeleton is a closed edge path that follows the edge directions and repeats no vertex except its initial vertex at the end. We identify it with the 1-chain obtained by summing its oriented edges. To decompose either 1-cycle in \eqref{eq:cup-chains}, discard edges of coefficient zero and orient each remaining edge so that its coefficient is positive. Since the boundary vanishes, the total weight entering each vertex equals the total weight leaving it.  By the  decomposition theorem~\cite[\S~3.5]{AhujaMagnantiOrlin1993}, every 1-cycle is a finite sum of positively weighted directed cycles. Applying this to the two 1-skeletons gives
\[
 C_{\mathrm{pr}}^{\mathrm c}(a)=\sum_i r_i\gamma_i,\qquad
 C_{\mathrm{du}}^{\mathrm c}(a)=\sum_j s_j\gamma_j^*.
\]
Here $\gamma_i$ is a directed cycle on the primal 1-skeleton with weight $r_i>0$, and $\gamma_j^*$ is a directed cycle on the dual 1-skeleton with weight $s_j>0$.  Theorem~\ref{thm:cup-linking} and bilinearity give
\begin{equation}
 \Hc(a)=\Lk\Bigl(\sum_i r_i\gamma_i,\sum_j s_j\gamma_j^*\Bigr)
       =\sum_{i,j}r_i s_j\Lk(\gamma_i,\gamma_j^*).
 \label{eq:weighted-cup-linking}
\end{equation}

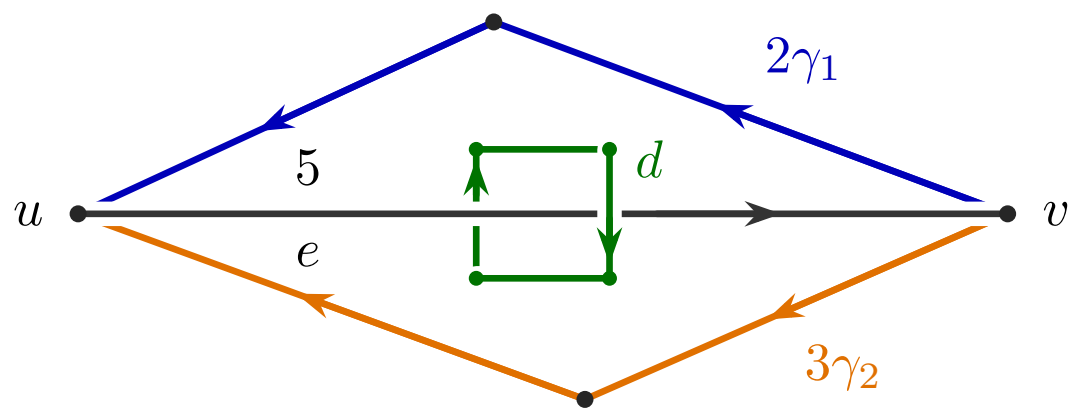
\begin{figure}[H]
\centering
\begin{tikzpicture}[x=1.16cm,y=1.04cm,line cap=round,line join=round,
  >={Stealth[length=2.4mm,width=1.65mm]},
  upper/.style={draw=blue!72!black,line width=1.25pt},
  lower/.style={draw=orange!88!black,line width=1.25pt},
  meridian/.style={draw=green!45!black,line width=1.25pt},
  common/.style={draw=black!80,line width=1.25pt}]
\coordinate (u) at (-2.65,0); \coordinate (v) at (2.65,0);
\coordinate (p) at (-.28,1.22); \coordinate (q) at (.24,-1.18);
\def\cycleedge#1#2#3{\draw[#1] (#2)--(#3);
  \path (#2)--coordinate[pos=.57] (tip) (#3);
  \draw[#1,->] (#2)--(tip);}
\cycleedge{upper}{v}{p} \cycleedge{upper}{p}{u}
\cycleedge{lower}{v}{q} \cycleedge{lower}{q}{u}
\draw[meridian] (-.38,-.41)--(-.38,.41)--(.38,.41);
\draw[meridian,->] (-.38,.12)--(-.38,.35);
\draw[white,line width=4.6pt] (u)--(v);
\draw[common] (u)--(v); \draw[common,->] (.65,0)--(1.36,0);
\draw[meridian,preaction={draw=white,line width=4.6pt}] (.38,.41)--(.38,-.41)--(-.38,-.41);
\draw[meridian,->] (.38,-.03)--(.38,-.32);
\foreach \pt in {u,v,p,q}{\fill[black!85] (\pt) circle (1.6pt);}
\foreach \xx in {-.38,.38}{\foreach \yy in {-.41,.41}{\fill[green!45!black] (\xx,\yy) circle (1.4pt);}}
\node[left=3pt] at (u) {$u$}; \node[right=3pt] at (v) {$v$};
\path (v)--(p) node[pos=.53,above right=2pt,text=blue!72!black] {$2\gamma_1$};
\path (v)--(q) node[pos=.55,below right=2pt,text=orange!88!black] {$3\gamma_2$};
\node[above=2pt] at (-1.34,0) {$5$};
\node[below=2pt] at (-1.34,0) {$e$};
\node[anchor=west,text=green!45!black,inner sep=2pt] at (.46,.35) {$d$};
\end{tikzpicture}
\caption{The blue and orange cycles traverse the black edge in the same direction. The green cycle $d=\partial e^*$ links each of them once.}
\label{fig:directed-cycle-decomposition}
\end{figure}
This weighted linking formula has a continuous counterpart for smooth flux tubes.

\begin{example}
In this example the ambient space is $\R^3$. Suppose that $X$ is supported
in the union of pairwise disjoint solid tori $U_1,\ldots,U_N$. Choose product coordinates
$U_i\cong S^1\times D^2$ in which $X$ is tangent to the circles
$S^1\times\{z\}$, and assume that any two distinct such circles in the
same tube have linking number zero. Orient the core $\gamma_i$ and a
meridional disk $D_i$ so that their intersection is positive, and write
$\Phi_i=\int_{D_i}\iota_X\mu$. Then helicity is the quadratic form
\begin{equation}
 \Hsm(X)=\sum_{i\ne j}\Phi_i\Phi_j\Lk(\gamma_i,\gamma_j)
 =2\sum_{1\le i<j\le N}\Phi_i\Phi_j\Lk(\gamma_i,\gamma_j).
 \label{eq:two-flux-tubes}
\end{equation}

\end{example}
\begin{proof}
Fix a tube $U_i$ and write $\beta=\iota_X\mu=d\alpha$. Use product
coordinates $(\theta,u,v)$ on $U_i=S^1\times D_i$, with $\theta$
increasing in the direction of the oriented core. Holding $z=(u,v)$
fixed and going once around $S^1$ traces the longitude
$\ell_z=S^1\times\{z\}$. The center $z=0$ gives the core
$\gamma_i=\ell_0$. Since $X$ is tangent to the longitudes, we can write
\[
 X=f(\theta,u,v)\,\partial_\theta,\qquad
 \mu=\rho(\theta,u,v)\,d\theta\wedge du\wedge dv.
\]
Contracting with $X$ gives $\beta=\rho f\,du\wedge dv$. The field is
divergence-free, so
\[
 0=d\beta=\partial_\theta(\rho f)\,d\theta\wedge du\wedge dv.
\]
Thus $\rho f=g_i(u,v)$ is independent of $\theta$. A small patch of the
cross section carries the same flux as it moves around the tube, although
the flux density $g_i$ may vary across the disk.

Writing $\alpha=A\,d\theta+B\,du+C\,dv$ and substituting $\beta=g_i\,du\wedge dv$ gives
\[
\begin{aligned}
 \alpha\wedge\beta
 &= (A\,d\theta+B\,du+C\,dv)\wedge(g_i\,du\wedge dv)\\
 &= A g_i\,d\theta\wedge du\wedge dv
    + B g_i\,\underbrace{du\wedge du}_{0}\wedge dv
    + C g_i\,\underbrace{dv\wedge du\wedge dv}_{0}\\
 &= A g_i\,d\theta\wedge du\wedge dv.
\end{aligned}
\]
Fubini's theorem now gives
\[
 \int_{U_i}\alpha\wedge\beta
 =\int_{D_i}\left(\int_{S^1}A(\theta,u,v)\,d\theta\right)g_i(u,v)\,du\,dv
 =\int_{D_i}\left(\oint_{\ell_z}\alpha\right)g_i(u,v)\,du\,dv.
\]
The inner integral is the circulation around $\ell_z$, since $u$ and
$v$ are fixed along that curve.

Draw the radial segment from the center of the  disk $D_i$
to any other point $z$ in the disk. Carrying this segment once
around the tube sweeps out an annulus $A_z$ bounded by $\gamma_i$ and
$\ell_z$. The field $X$ is tangent to the annulus, so its flux through
$A_z$ is zero. Orienting $A_z$ so that $\partial A_z=\ell_z-\gamma_i$,
Stokes' theorem gives
\[
 \oint_{\ell_z}\alpha-\oint_{\gamma_i}\alpha
 =\int_{A_z}d\alpha=\int_{A_z}\iota_X\mu=0.
\]
The coefficient $A$ may vary from one longitude to another, but its
integral around a full turn has the same value on every longitude.
We can therefore take that common value outside the integral and obtain
\[
 \int_{U_i}\alpha\wedge\beta
 =\left(\oint_{\gamma_i}\alpha\right)
   \left(\int_{D_i}g_i(u,v)\,du\,dv\right)
 =\Phi_i\oint_{\gamma_i}\alpha.
\]

It remains to compute the core circulation. Choose a compact oriented
smooth surface $S_i$ with $\partial S_i=\gamma_i$. By Stokes' theorem
and the support condition $\operatorname{supp}(\beta)\subset\bigcup_{j=1}^N U_j$,
\[
 \oint_{\gamma_i}\alpha=\int_{S_i}\beta
 =\sum_{j=1}^N\int_{S_i\cap U_j}\beta.
\]
Let $p_j\colon U_j\to D_j$ be the projection $p_j(\theta,u,v)=(u,v)$.
By Sard's theorem applied to $p_j$ on the interior of $S_i\cap U_j$,
$S_i$ meets the longitude over almost every $z\in D_j\setminus\{0\}$
transversely. At each such intersection $x\in S_i\cap\ell_z$, the inverse
function theorem gives a neighborhood $A$ of $x$ in $S_i\cap U_j$ such
that $p_j|_A\colon A\to B=p_j(A)$ is a diffeomorphism. Since
$\beta|_{U_j}=p_j^*(g_j(u,v)\,du\wedge dv)$, change of variables gives
\[
 \int_A\beta=\varepsilon_x\int_B g_j(u,v)\,du\,dv,\qquad
 \varepsilon_x=
 \begin{cases}
  +1 & \text{if $p_j|_A$ preserves orientation},\\
  -1 & \text{if $p_j|_A$ reverses orientation}.
 \end{cases}
\]
Since $\ell_z$ is oriented by $\partial_\theta$, $\varepsilon_x$ is the
intersection sign of $S_i$ and $\ell_z$ at $x$~\cite{ArnoldKhesin2021}.
For $j\ne i$, the longitude $\ell_z$ can be moved to $\gamma_j$ inside
$U_j$ without meeting $\gamma_i$. Thus its signed crossing count is
$I(S_i,\ell_z)=\Lk(\gamma_i,\gamma_j)$ for almost every $z=(u,v)\in D_j$.
Integrating the flux density with this signed count gives
\[
 \int_{S_i\cap U_j}\beta
 =\int_{D_j}g_j(u,v)\,I(S_i,\ell_z)\,du\,dv
 =\Lk(\gamma_i,\gamma_j)\underbrace{\int_{D_j}g_j(u,v)\,du\,dv}_{\Phi_j}.
\]

The first equality also applies to $U_i$. For almost every $z\ne0$,
$\ell_z$ is transverse to $S_i$ and disjoint from $\partial S_i=\gamma_i=\ell_0$.
Hence $I(S_i,\ell_z)=\Lk(\gamma_i,\ell_z)=0$ by the hypothesis on distinct
longitudes in $U_i$. Since $z=0$ has measure zero in $D_i$, the integral
gives $\int_{S_i\cap U_i}\beta=0$. Summing over $i$ and using symmetry of linking gives
\[
 \Hsm(X)=\sum_i\Phi_i\oint_{\gamma_i}\alpha
 =\sum_{i\ne j}\Phi_i\Phi_j\Lk(\gamma_i,\gamma_j)
 =2\sum_{i<j}\Phi_i\Phi_j\Lk(\gamma_i,\gamma_j),
\]
as claimed in \eqref{eq:two-flux-tubes}.
\end{proof}

\begin{figure}[!htbp]
\centering
\begin{tikzpicture}[x=1cm,y=1cm,>=stealth,font=\small]
  \definecolor{rimblue}{RGB}{35,85,123}
  \definecolor{leftochre}{RGB}{163,77,30}
  \definecolor{rightgreen}{RGB}{28,111,92}
  \newcommand{\drawfluxtube}[2]{%
    \draw[#1!85!black,line width=6.3pt] #2;
    \draw[#1!32!white,line width=4.8pt] #2;
    \draw[#1!90!black,line width=.55pt] #2;}
  \def\spanningrim{(0,0) ellipse[x radius=2.17,y radius=.66]}
  % Shear and bend the right-hand ring smoothly, keeping its two disk-plane crossings.
  \def\rightbranch#1#2{plot[domain=#1:#2,samples=41,smooth,variable=\t]
    ({2.02-.80*cos(\t)+.44*sin(\t)+.13*sin(\t)^2},
     {.89*sin(\t)+.13*sin(\t)^2})}
  \drawfluxtube{leftochre}{(-2.02,0) ellipse[x radius=.80,y radius=.89]}
  \drawfluxtube{rightgreen}{\rightbranch{0}{360} -- cycle}
  \fill[rimblue!8!white] \spanningrim;
  \begin{scope}
    \clip \spanningrim;
    \draw[leftochre!70,densely dashed,line width=.65pt] (-1.22,0)
      arc[start angle=0,end angle=-180,x radius=.80,y radius=.89];
    \draw[rightgreen!70,densely dashed,line width=.65pt] \rightbranch{180}{360};
  \end{scope}
  \drawfluxtube{rimblue}{\spanningrim}
  \drawfluxtube{leftochre}{(-1.22,0)
    arc[start angle=0,end angle=180,x radius=.80,y radius=.89]}
  \drawfluxtube{rightgreen}{\rightbranch{0}{180}}
  % The small ellipses are the two tube cross sections in the spanning disk.
  \foreach \xi/\col in {-1.22/leftochre,1.22/rightgreen}{
    \filldraw[fill=\col!35!white,draw=\col!85!black,line width=.65pt]
      (\xi,0) ellipse[x radius=.115,y radius=.047];
    \fill[\col!85!black] (\xi,0) circle[radius=.8pt];}
  \draw[->,leftochre!85!black,line width=1pt] (-1.22,.02)
    arc[start angle=0,end angle=27,x radius=.80,y radius=.89];
  \draw[->,rightgreen!85!black,line width=1pt] \rightbranch{3}{27};
  \draw[->,rimblue!85!black,line width=1pt]
    ({2.17*cos(250)},{.66*sin(250)})
    arc[start angle=250,end angle=268,x radius=2.17,y radius=.66];
  \node[text=leftochre!85!black] at (-2.02,1.14) {$\gamma_2$};
  \node[text=rightgreen!85!black] at (2.58,1.27) {$\gamma_3$};
  \node[text=rimblue!85!black] at (0,-.95) {$\gamma_1,\ \Phi_1$};
  \node[text=leftochre!85!black] at (-.84,-.10) {$+\Phi_2$};
  \node[text=rightgreen!85!black] at (.84,-.10) {$+\Phi_3$};
  \node[text=rimblue!85!black] at (-.18,.10) {$S$};
  \draw[->,rimblue!85!black,semithick] (.12,-.04) -- (.12,.37);
  \begin{scope}[xshift=6.8cm,
    core strand/.style={draw=#1!85!black,line width=1pt,line join=round,
      preaction={draw=white,line width=3.6pt}}]
    % Use the same projections as the smooth cores, sampled at 12 vertices.
    \foreach \j in {0,...,11}{
      \pgfmathsetmacro{\ang}{30*\j}
      \coordinate (rim\j) at ({2.17*cos(\ang)},{.66*sin(\ang)});
      \coordinate (och\j) at ({-2.02+.80*cos(\ang)},{.89*sin(\ang)});
      \coordinate (grn\j) at
        ({2.02-.80*cos(\ang)+.44*sin(\ang)+.13*sin(\ang)^2},
         {.89*sin(\ang)+.13*sin(\ang)^2});}
    \foreach \prefix/\col in {och/leftochre,grn/rightgreen,rim/rimblue}{
      \draw[core strand=\col] (\prefix0)
        \foreach \j in {1,...,11}{ -- (\prefix\j)} -- cycle;}
    % The upper arcs pass over the rim, and the lower arcs pass under it.
    \foreach \prefix/\col in {och/leftochre,grn/rightgreen}{
      \draw[core strand=\col] (\prefix0)
        \foreach \j in {1,...,6}{ -- (\prefix\j)};}
    \foreach \prefix/\col in {och/leftochre,grn/rightgreen,rim/rimblue}{
      \foreach \j in {0,...,11}{\fill[\col!85!black] (\prefix\j) circle (1pt);}}
    \foreach \prefix/\col/\j/\k in {och/leftochre/0/1,grn/rightgreen/0/1,rim/rimblue/8/9}{
      \draw[->,\col!85!black,line width=1pt] (\prefix\j) -- (\prefix\k);}
    \node[text=leftochre!85!black] at (-2.02,1.14) {$\Phi_2 c_2$};
    \node[text=rightgreen!85!black] at (2.58,1.27) {$\Phi_3 c_3$};
    \node[text=rimblue!85!black] at (0,-.95) {$\Phi_1 c_1$};
  \end{scope}
\end{tikzpicture}
\caption{For $N=3$, each flux tube (left) retracts onto its core, represented on the right by an oriented polygonal $1$-cycle $c_i$ of weight $\Phi_i$. The pairs $(c_1,c_2)$ and $(c_1,c_3)$ each link once positively, while $c_2$ and $c_3$ are unlinked. The discrete weighted linking sum is therefore $2\Phi_1(\Phi_2+\Phi_3)=\Hsm(X)$.}
\label{fig:flux-tube-piercing}
\end{figure}
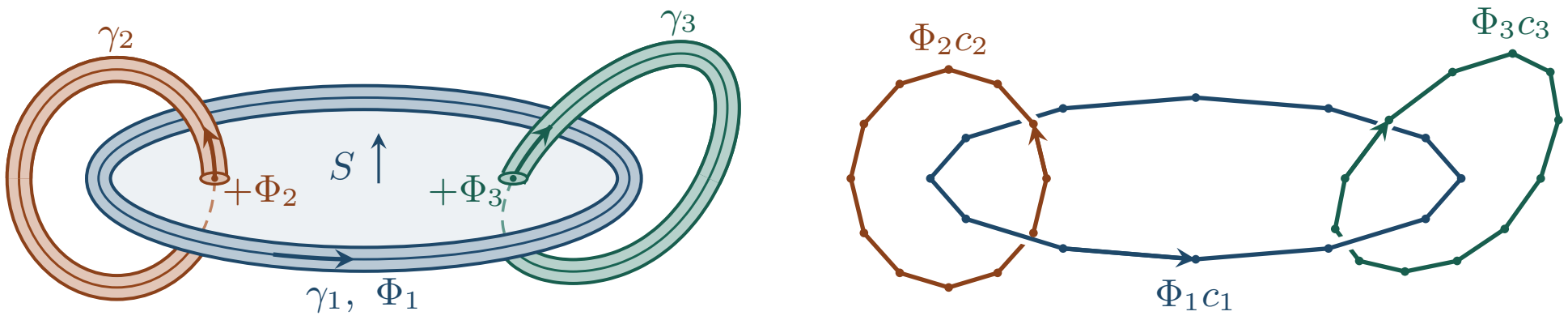

This connection provides us an understanding on how to correspond the discrete objects to the continuous ones. In \eqref{eq:weighted-cup-linking}, $r_i$ and $s_j$ play the role of
the fluxes carried around the chosen cycles. This interpretation can
be read off on a dual edge $f^*$. Its coefficient is
$b(f)=\int_f W_2b$ by \eqref{eq:whitney-identities}, so the oriented
sum of the weights $s_j$ of the cycles traversing $f^*$ is exactly
the flux through $f$.  We call the following explanation the flux tube model.

This model suggests viewing each cycle as the core curve
onto which a solid torus deformation retracts, and its coefficient
as the signed flux through a meridional disk of that torus.
In \S\ref{sec:whitney-flux} we make this interpretation precise for
the primal and dual edge coefficients.
There we will encounter the same idea of retraction again, this time with a cell retracting onto the vertex that represents it. 

\section{A linking number realization for Whitney helicity }
\label{sec:whitney-flux}

The flux through a meridional disk gives the weight of a tube core. On
the mesh we use the flux through the face separating two cells to weight
the edge joining their representative vertices. Applied to $W_2(\delta a)$, this rule gives one cycle on each
1-skeleton. We then use \eqref{eq:discrete-wedge-product}, which expresses
$a\wedge\delta a$ on each tetrahedron as a signed average of cup products
over its vertex orderings, to show that the linking number of these two
cycles equals $\Hw(a)$.

\subsection{Recording a vector field on the 1-skeleton}
\label{sec:recording-vector-field}

We begin with a heuristic explanation of the results that follow. One may picture each tetrahedron as filled with fluid whose velocity is
given by $X$. Suppose that we follow
a single fluid particle over a long time and observe only the cells it visits. Represent each
cell by a vertex and consider a trajectory segment with finitely many
transverse crossings, all through face interiors. Each time the particle crosses from one cell into the next, we record the
oriented edge joining their representative vertices. Its coefficient in the
resulting 1-chain is the number of crossings in its chosen direction minus
the number in the opposite direction.

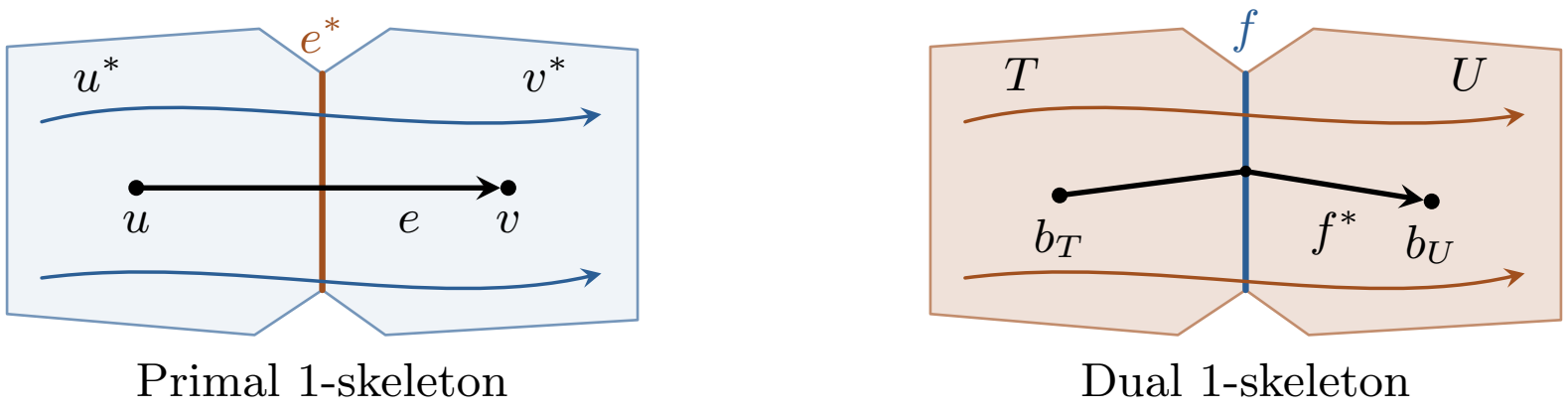
\begin{figure}[!htbp]
\centering
\begin{tikzpicture}[x=1.05cm,y=1cm,>=stealth,font=\small,
  line cap=round,line join=round]
  \definecolor{recordblue}{RGB}{43,94,149}
  \definecolor{recordochre}{RGB}{161,80,30}
  \def\recordcells#1#2#3{%
    \filldraw[fill=#1,draw=#2!65,line width=.5pt]
      (-2,.9)--(-.40,1.02)--(0,.72)--(0,-.72)--(-.43,-1.02)--(-2,-.88)--cycle;
    \filldraw[fill=#1,draw=#2!65,line width=.5pt]
      (0,.72)--(.43,1.02)--(2,.88)--(2,-.92)--(.40,-1.02)--(0,-.72)--cycle;
    \draw[#3,line width=1.2pt] (0,-.72)--(0,.72);
    \draw[->,#2,line width=.65pt] (-1.78,.40) .. controls (-.82,.67) and (.55,.26) .. (1.77,.45);
    \draw[->,#2,line width=.65pt] (-1.78,-.64) .. controls (-.62,-.48) and (.70,-.87) .. (1.77,-.61);}
  \begin{scope}
    \recordcells{recordblue!7}{recordblue}{recordochre}
    \node at (-1.42,.72) {$u^*$}; \node at (1.42,.72) {$v^*$};
    \node[text=recordochre,above] at (0,.74) {$e^*$};
    \draw[->,line width=1pt,shorten >=1.5pt] (-1.18,-.04)--(1.18,-.04);
    \fill (-1.18,-.04) circle (1.5pt) (1.18,-.04) circle (1.5pt);
    \node[below=1pt] at (-1.18,-.04) {$u$}; \node[below=1pt] at (1.18,-.04) {$v$};
    \node[below=1pt] at (.55,-.04) {$e$};
    \node at (0,-1.32) {Primal 1-skeleton};
  \end{scope}
  \begin{scope}[xshift=6.15cm]
    \recordcells{recordochre!17}{recordochre}{recordblue}
    \node at (-1.42,.72) {$T$}; \node at (1.42,.72) {$U$};
    \node[text=recordblue,above] at (0,.74) {$f$};
    \draw[->,line width=1pt,shorten >=1.5pt] (-1.18,-.09)--(0,.07)--(1.18,-.13);
    \fill (-1.18,-.09) circle (1.5pt) (0,.07) circle (1.1pt) (1.18,-.13) circle (1.5pt);
    \node[below=1pt] at (-1.18,-.09) {$b_T$}; \node[below=1pt] at (1.18,-.13) {$b_U$};
    \node[below=1pt] at (.57,-.03) {$f^*$};
    \node at (0,-1.32) {Dual 1-skeleton};
  \end{scope}
\end{tikzpicture}
\caption{Recording interface crossings, shown in a planar schematic.
The curved arrows follow the flow. The black arrow records a passage
across $e^*$ by the primal edge $e$ (left), or across $f$ by the dual
edge $f^*$ (right).}
\label{fig:record-interface}
\end{figure}

The primal and dual 1-skeletons give two parallel ways to record these
crossings (Figure~\ref{fig:record-interface}).
On the primal 1-skeleton, think of the dual 3-cell $u^*$ as contracted to
its vertex $u$. If a particle leaves $u^*$ for $v^*$ across their common
dual face $e^*$, we record $e=[u,v]$. Dually, when the particle crosses a face $f$ from $T$ to $U$, we record
the dual edge $f^*$ oriented from $b_T$ to $b_U$.

We therefore take the coefficient of each oriented edge to be the flux of
$X$ through the interface it crosses. For a primal edge $e$ this flux is
$\int_{e^*}\iota_X\mu$, and for a dual edge $f^*$ it is $\int_f\iota_X\mu$.
In each case the flux measures the net volume transported across the
interface per unit time, with positive sign in the direction of the edge.
For a closed $2$-cochain $b$, we compute these fluxes using $W_2b$ in
place of $\iota_X\mu$, summing over the oriented barycentric triangles
when the interface is a dual face.

\begin{definition}
\label{def:flux-chains}
For $b\in Z^2(K)$ and $\beta=W_2b$, define the primal and dual 1-chains by
\begin{equation}
 \cpr:=\sum_{e\in K^{(1)}}\left(\int_{e^*}\beta\right)e,
 \qquad
 \cdu:=\sum_{f\in K^{(2)}}\left(\int_f\beta\right)f^*.
 \label{eq:flux-chain-definitions}
\end{equation}
\end{definition}

Every edge in these two sums is oriented. Reversing its orientation changes
the signs of both the coefficient and the edge it multiplies. On a primal
face $f$ the coefficient is simply $b(f)$, by $R_2W_2=\mathrm{id}$.
On a dual face $e^*$ it receives a contribution from every tetrahedron
containing $e$. In \S\ref{sec:markov} we use them to construct trajectories on the
two skeletons.

\subsection{The linking formula}

Our goal is to prove that Whitney helicity indeed equals the linking number of
$\cpr$ and $\cdu$, the two 1-chains suggested by the flux interpretation
above. Before turning to the proof, we carry out a few elementary
calculations on the tetrahedron.

\begin{lemma}
\label{lem:local-flux-whitney}
Let $b\in Z^2(K)$ and set $\beta=W_2b$. Orient each tetrahedron $T$ by the
orientation of $M$. Let $f_i$ be the face opposite vertex $i$, with the
outward boundary orientation induced by $T$, and set $B_{T,i}=b(f_i)$.
For each oriented edge $e=[u,v]$ of $T$,
\begin{equation}
 \int_{e^*\cap T}\beta
 =\int_T W_1\chi_e\wedge\beta
 =\frac{B_{T,u}-B_{T,v}}{12}.
 \label{eq:local-flux-whitney}
\end{equation}
\end{lemma}

\begin{figure}[!htbp]
\centering
\begin{tikzpicture}[font=\small,line join=round,line cap=round]
  \definecolor{patchblue}{RGB}{45,105,149}
  \definecolor{patchedge}{RGB}{165,55,48}
  \tikzset{patch rim/.style={draw=patchblue!85!black,fill=patchblue!19,line width=.9pt}}
  \begin{scope}[x=1.07cm,y=1.04cm]
    \coordinate (A) at (-2.72,-.60);
    \coordinate (B) at (1.12,-1.37);
    \coordinate (C) at (2.32,-.18);
    \coordinate (D) at (.32,1.89);
    \coordinate (M) at (barycentric cs:A=1,B=1);
    \coordinate (P) at (barycentric cs:A=1,B=1,C=1);
    \coordinate (R) at (barycentric cs:A=1,B=1,D=1);
    \coordinate (G) at (barycentric cs:A=1,B=1,C=1,D=1);
    \draw[gray!65,densely dashed,line width=.5pt] (A)--(C);
    \draw[gray!75,line width=.65pt] (A)--(D)--(C)--(B) (B)--(D);
    \filldraw[patch rim] (M)--(P)--(G)--(R)--cycle;
    \draw[patchedge,line width=1.3pt] (A)--(B);
    \path (A)--(B) node[pos=.24,below=2pt,text=patchedge] {$e$};
    \foreach \v/\n/\pos in {A/0/left,B/1/below,C/2/right,D/3/above}{
      \fill[black!75] (\v) circle (1.2pt) node[\pos=2pt] {$\n$};}
    \foreach \v/\n/\pos in {M/m/below left,P/p/right,R/r/above left,G/t/above right}{
      \fill[patchblue!85!black] (\v) circle (1.4pt) node[\pos=1pt,text=black] {$\n$};}
  \end{scope}
  \begin{scope}[shift={(4.20cm,.16cm)},x=.80cm,y=.80cm]
    \coordinate (Mflat) at (0,0);
    \coordinate (Tflat) at (3,0);
    \coordinate (Pflat) at (2,{-sqrt(8/3)});
    \coordinate (Rflat) at (2,{sqrt(8/3)});
    \filldraw[patch rim] (Mflat)--(Pflat)--(Tflat)--(Rflat)--cycle;
    \foreach \v in {Mflat,Tflat,Pflat,Rflat}{\fill[patchblue!85!black] (\v) circle (1.4pt);}
    \node[left=3pt] at (Mflat) {$m=b_{01}$};
    \node[below=3pt] at (Pflat) {$p=b_{012}$};
    \node[right=3pt] at (Tflat) {$t=b_T$};
    \node[above=3pt] at (Rflat) {$r=b_{013}$};
  \end{scope}
\end{tikzpicture}
\caption{The edge $e=[0,1]$ in a reference tetrahedron $T$ (left).
The shaded region in both views is $e^*\cap T$, the part of the dual cell $e^*$ inside $T$.}
\label{fig:local-dual-patch}
\end{figure}
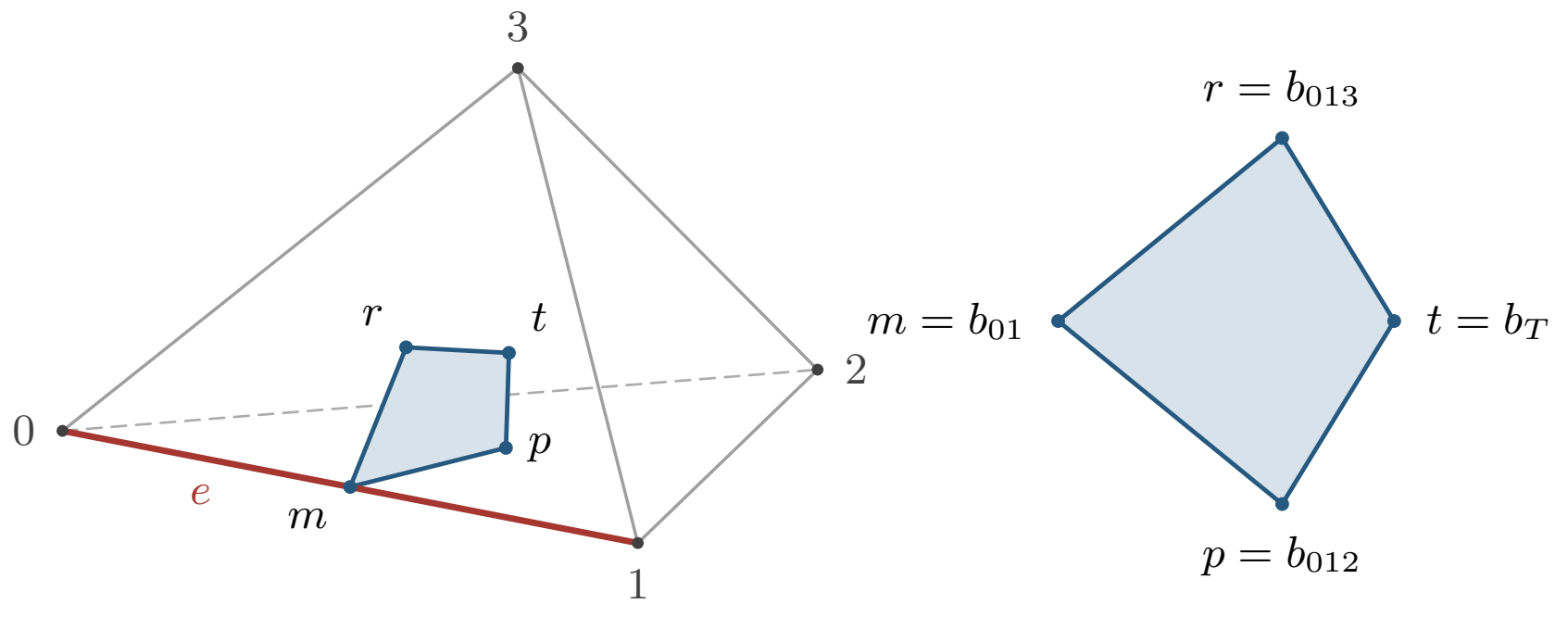

\begin{proof}
Write $\widehat e=\Phi_T^{-1}(e)$ and $\widehat\beta_T=\Phi_T^*(\beta|_T)$.
By construction, $\Phi_T$ maps the two triangles of the dual cell
$\widehat e^*$ inside $\widehat T$ onto those of $e^*\cap T$, preserving
their orientations. Moreover, the Whitney form on $T$ was defined in
\S\ref{sec:whitney-de-rham} by transporting the reference form, so
$\Phi_T^*(W_1\chi_e|_T)=W_1\chi_{\widehat e}$. Change of variables gives
\[
\begin{aligned}
 \int_{e^*\cap T}\beta
 &=\int_{\widehat e^*\cap\widehat T}\widehat\beta_T,\\
 \int_T W_1\chi_e\wedge\beta
 &=\int_{\widehat T}W_1\chi_{\widehat e}\wedge\widehat\beta_T.
\end{aligned}
\]

It therefore suffices to prove the identity on the reference tetrahedron.
Let $p_0,\ldots,p_3$ be the vertices of $\widehat T=[0,e_1,e_2,e_3]$,
whose volume is $\widehat V=1/6$. Fix $i$ and order the other three vertices so that
$[p_i,p_j,p_k,p_\ell]$ is positively oriented. Its outward face opposite
$p_i$ is then $f_i=[p_j,p_k,p_\ell]$. Let $\chi_{f_i}$ be the elementary
2-cochain that takes the value $1$ on $f_i$ and $0$ on the other faces.
On the Euclidean reference tetrahedron, the canonical identification
$T_{\mathbf{x}}\R^3\simeq\R^3$ lets us regard the displacement
$\mathbf{x}-p_i$ as a tangent vector at $\mathbf{x}$. Thus
$\mathbf{x}\mapsto\mathbf{x}-p_i$ defines a smooth vector field on
$\widehat T$. We claim that the Whitney form $W_2\chi_{f_i}$ satisfies
\[
 (W_2\chi_{f_i})_{\mathbf{x}}
 =\frac{1}{3\widehat V}\,\iota_{\mathbf{x}-p_i}dV_{\widehat T}.
\]
By \eqref{eq:whitney-form},
\[
 W_2\chi_{f_i}=2\bigl(\lambda_j\,d\lambda_k\wedge d\lambda_\ell
 -\lambda_k\,d\lambda_j\wedge d\lambda_\ell
 +\lambda_\ell\,d\lambda_j\wedge d\lambda_k\bigr).
\]
The barycentric coordinates are affine functions and their differentials
are constant. Thus the coefficients of this Whitney form in Cartesian
coordinates are affine in $\mathbf{x}$. For the right-hand side, write $\mathbf{x}=(x,y,z)$ and $p_i=(x_i,y_i,z_i)$.
Contracting the Euclidean volume form with the vector field $\mathbf{x}\mapsto\mathbf{x}-p_i$ gives
\[
 \iota_{\mathbf{x}-p_i}dV_{\widehat T}
 =(x-x_i)\,dy\wedge dz+(y-y_i)\,dz\wedge dx+(z-z_i)\,dx\wedge dy.
\]
The three coefficients are affine functions of $\mathbf{x}$, and the factor $(3\widehat V)^{-1}$ is constant. Now any affine function $g$ on $\widehat T$ satisfies
\[
 g(\mathbf{x})=\sum_{s=0}^3\lambda_s(\mathbf{x})\,g(p_s),
\]
because $\mathbf{x}=\sum_s\lambda_s(\mathbf{x})p_s$ and
$\sum_s\lambda_s(\mathbf{x})=1$. Applying this to each coefficient of
the difference of the two forms shows that equality at the four vertices
implies equality throughout $\widehat T$.

At $p_i$ all three barycentric coordinates in this expression vanish,
and the vector $\mathbf{x}-p_i$ is zero. At the other three vertices,
write $u_s=p_s-p_i$ for $s=j,k,\ell$. The covectors
$d\lambda_j,d\lambda_k,d\lambda_\ell$ are dual to $u_j,u_k,u_\ell$.
Since $\det(u_j,u_k,u_\ell)=6\widehat V$, we have
$dV_{\widehat T}=6\widehat V\,d\lambda_j\wedge d\lambda_k\wedge d\lambda_\ell$.
Evaluating the Whitney form and contracting this volume form at each
of the three vertices gives
\[
\begin{aligned}
 (W_2\chi_{f_i})_{p_j}&= 2\,d\lambda_k\wedge d\lambda_\ell
       =\frac{1}{3\widehat V}\,\iota_{u_j}dV_{\widehat T},\\
 (W_2\chi_{f_i})_{p_k}&=-2\,d\lambda_j\wedge d\lambda_\ell
       =\frac{1}{3\widehat V}\,\iota_{u_k}dV_{\widehat T},\\
 (W_2\chi_{f_i})_{p_\ell}&= 2\,d\lambda_j\wedge d\lambda_k
       =\frac{1}{3\widehat V}\,\iota_{u_\ell}dV_{\widehat T}.
\end{aligned}
\]
This proves the claimed formula throughout $\widehat T$. Now
$\widehat\beta_T=\sum_i B_{T,i}W_2\chi_{f_i}$. The terms containing
$\mathbf{x}$ cancel because $\sum_i B_{T,i}=(\delta b)(T)=0$, giving
\begin{equation}
 \widehat\beta_T=\iota_{v_T}dV_{\widehat T},
 \qquad v_T=-\frac{1}{3\widehat V}\sum_{i=0}^3 B_{T,i}p_i.
 \label{eq:flux-reference-field}
\end{equation}

We now compute the two integrals for the oriented edge $e=[u,v]$. After an affine change of coordinates on $\widehat T$ if necessary, we may assume that
$\Phi_T(0)=u$ and $\Phi_T(e_1)=v$. We order the other two vertices so that
$\Phi_T$ preserves orientation. In these coordinates the reference
edge $\widehat e=[0,1]$ points in the positive $x$-direction. Write
$v_T=(q_1,q_2,q_3)$ for the constant vector in
\eqref{eq:flux-reference-field}. Its flux form is
\begin{equation}
 \widehat\beta_T=\iota_{v_T}(dx\wedge dy\wedge dz)
 =q_1\,dy\wedge dz+q_2\,dz\wedge dx+q_3\,dx\wedge dy.
 \label{eq:flux-standard-form}
\end{equation}

We first compute the left-hand side of \eqref{eq:local-flux-whitney}. Inside $\widehat T$ the edge $[0,1]$
belongs to exactly two faces  $[0,1,2]$ and $[0,1,3]$. Each contributes
one triangle to $\widehat e^*\cap\widehat T$, the part of the dual cell
inside $\widehat T$. Its vertices are the barycenters of the
edge, the face, and the tetrahedron. Averaging the corresponding vertex
coordinates gives
\[
\begin{aligned}
 m=b_{01}&=(\tfrac12,0,0), & p=b_{012}&=(\tfrac13,\tfrac13,0),\\
 r=b_{013}&=(\tfrac13,0,\tfrac13), & t=b_{\widehat T}&=(\tfrac14,\tfrac14,\tfrac14).
\end{aligned}
\]
For an oriented triangle $[A,B,C]$, its area vector is
$\frac12(B-A)\times(C-A)$. It points in the oriented normal direction
and its length is the area of the triangle. The flux of a constant
field through the triangle is its normal component times this area,
hence the dot product of $v_T$ with the area vector. For our two triangles,
\[
 \tfrac12(p-m)\times(t-m)=\tfrac1{48}(2,1,1),\qquad
 \tfrac12(t-m)\times(r-m)=\tfrac1{48}(2,1,1).
\]
Both vectors have positive dot product with the edge tangent $e_1=(1,0,0)$,
so $[m,p,t]$ and $[m,t,r]$ have the dual orientation prescribed by
\eqref{eq:dual-orientation}. Adding their fluxes gives
\begin{equation}
 \int_{\widehat e^*\cap\widehat T}\widehat\beta_T
 =v_T\cdot\tfrac1{24}(2,1,1)=\frac{2q_1+q_2+q_3}{24}.
 \label{eq:flux-standard-patch}
\end{equation}

We next compute the volume integral for the same edge. Its endpoint
barycentric coordinates are $\lambda_0=1-x-y-z$ and $\lambda_1=x$.
Substitution into the Whitney formula gives
\[
 W_1\chi_{01}=\lambda_0\,d\lambda_1-\lambda_1\,d\lambda_0
 =(\lambda_0+\lambda_1)\,dx+\lambda_1\,dy+\lambda_1\,dz.
\]
When we wedge this 1-form with \eqref{eq:flux-standard-form}, the $dx$
term pairs only with $q_1\,dy\wedge dz$, since all other products with
$dx$ contain $dx$ twice. The $dy$ and $dz$ terms likewise contribute
$q_2\lambda_1$ and $q_3\lambda_1$ to the coefficient of $dV_{\widehat T}$.
Thus
\[
 W_1\chi_{01}\wedge\widehat\beta_T
 =\bigl(q_1(\lambda_0+\lambda_1)+(q_2+q_3)\lambda_1\bigr)dV_{\widehat T}.
\]
The four barycentric coordinates have equal integrals, since affine
permutations of the vertices preserve Euclidean volume. Their sum is $1$,
so each has integral $\widehat V/4=1/24$. It follows that
\begin{equation}
 \int_{\widehat T}W_1\chi_{01}\wedge\widehat\beta_T
 =q_1\Bigl(\frac1{24}+\frac1{24}\Bigr)+(q_2+q_3)\frac1{24}
 =\frac{2q_1+q_2+q_3}{24}.
 \label{eq:flux-standard-pairing}
\end{equation}
This agrees with the flux through $\widehat e^*\cap\widehat T$ computed in
\eqref{eq:flux-standard-patch}.

It remains to express this common value in terms of the original 2-cochain
$b$. Whitney interpolation reproduces the integral on each oriented face,
so $B_{T,i}=b(f_i)$ is the flux of $v_T$ through the corresponding outward
reference face. The face opposite $0$ is $[1,2,3]$, with outward area
vector $\tfrac12(e_2-e_1)\times(e_3-e_1)=\tfrac12(1,1,1)$. The face
opposite $1$ lies in the plane $x=0$ and is oriented as $[0,3,2]$,
with outward area vector $\tfrac12 e_3\times e_2=\tfrac12(-1,0,0)$.
Taking their dot products with $v_T$ gives
\[
 B_{T,0}=\frac{q_1+q_2+q_3}{2},\qquad B_{T,1}=-\frac{q_1}{2},\qquad
 \frac{B_{T,0}-B_{T,1}}{12}=\frac{2q_1+q_2+q_3}{24}.
\]
The reference vertices $0,1$ correspond to $u,v$, respectively. Applying
the change of variables from the beginning of the proof now gives
\[
 \int_{e^*\cap T}\beta=\int_T W_1\chi_e\wedge\beta
 =\frac{B_{T,u}-B_{T,v}}{12}.
\]
\end{proof}

\begin{corollary}
\label{cor:global-flux-pairing}
For $b\in Z^2(K)$ and every $\eta\in C^1(K)$,
\begin{equation}
 \langle\eta,\cpr\rangle=\int_M W_1\eta\wedge\beta,
 \qquad \cdu=\mathfrak D_2 b.
 \label{eq:global-flux-pairing}
\end{equation}
In particular, when $b=\delta a$,
\begin{equation}
 \langle a,\cpr\rangle=\Hw(a).
 \label{eq:flux-helicity-pairing}
\end{equation}
With the outward face orientations induced by the orientation of each $T$
in $M$, the coefficient of the primal edge $[u,v]$ is
\begin{equation}
 \int_{[u,v]^*}\beta
 =\sum_{T\supset[u,v]}\frac{B_{T,u}-B_{T,v}}{12}.
 \label{eq:primal-interface-flux}
\end{equation}
\end{corollary}

\begin{proof}
By definition of $\cpr$, each edge $e$ has coefficient $\int_{e^*}\beta$.
 Applying \eqref{eq:local-flux-whitney} gives
\[
\begin{aligned}
 \langle\eta,\cpr\rangle
 &=\sum_{e\in K^{(1)}}\eta(e)\int_{e^*}\beta\\
 &=\sum_{e\in K^{(1)}}\sum_{T\supset e}\eta(e)\int_{e^*\cap T}\beta\\
 &=\sum_{T\in K^{(3)}}\sum_{e\subset T}\eta(e)\int_{e^*\cap T}\beta\\
 &=\sum_T\int_T\left(\sum_{e\subset T}\eta(e)W_1\chi_e\right)\wedge\beta
 =\int_M W_1\eta\wedge\beta.
\end{aligned}
\]
The dual identity follows from $\int_f W_2b=b(f)$. If $b=\delta a$, then
$\beta=dW_1a$ by \eqref{eq:whitney-identities}, so taking $\eta=a$ yields
\eqref{eq:flux-helicity-pairing}. Summing the local contributions at a
fixed edge gives \eqref{eq:primal-interface-flux}.
\end{proof}

\begin{theorem}
\label{thm:whitney-linking}
For $a\in C^1(K)$, the chains associated to $b=\delta a$ by
Definition~\ref{def:flux-chains} are 1-cycles and satisfy
\begin{equation}
 \Hw(a)=\Lk(\cpr,\cdu).
 \label{eq:whitney-flux-linking}
\end{equation}
\end{theorem}

\begin{proof}
By \eqref{eq:cup-chains} and \eqref{eq:global-flux-pairing},
$C_{\mathrm{du}}^{\mathrm c}(a)=\mathfrak D_2 b=\cdu$.
We now compare the edge coefficients of the primal cup chain
$C_{\mathrm{pr}}^{\mathrm c}(a)=[K]\frown b$ with those of $\cpr$. The cup linking formula \eqref{eq:cup-linking} holds
for each total order on the vertices of $K$. We average over all such orders,
so each of the $24$ orders on a tetrahedron occurs equally often. On a
tetrahedron labeled $T=[0,1,2,3]$ in its orientation induced by $M$, the
order $0<1<2<3$ gives $T\frown b=b[0,1,2]\,[2,3]$ by
\eqref{eq:cap-product}. Antisymmetrizing this expression gives
\[
 \frac1{24}\sum_{\pi\in S_4}\operatorname{sgn}(\pi)\,
 b[\pi(0),\pi(1),\pi(2)]\,[\pi(2),\pi(3)].
\]
To read off the coefficient of $[0,1]$, keep only the four permutations
ending in $(0,1)$ or $(1,0)$. Using $[1,0]=-[0,1]$, this coefficient is
\[
 \frac{b[2,3,0]-b[3,2,0]+b[2,3,1]-b[3,2,1]}{24}
 =\frac{b[0,2,3]+b[1,2,3]}{12}
 =\frac{B_{T,0}-B_{T,1}}{12}.
\]

Relabeling gives the coefficient $(B_{T,u}-B_{T,v})/12$ for any oriented
edge $e=[u,v]$ of $T$. By \eqref{eq:local-flux-whitney}, this is exactly
$\int_{e^*\cap T}\beta$. Summing over the tetrahedra containing $e$ gives
$\int_{e^*}\beta$, its coefficient in $\cpr$ by \eqref{eq:flux-chain-definitions}.

Write $T=[v_0,v_1,v_2,v_3]$ in its reference orientation for the fixed
vertex order. Averaging \eqref{eq:cup-linking} and using bilinearity of
linking gives
\[
\begin{aligned}
 \Lk(\cpr,\cdu)
 &=\frac{1}{24}\sum_{T\in K^{(3)}}\varepsilon_T
   \sum_{\pi\in S_4}\operatorname{sgn}(\pi)\,
   a[v_{\pi(0)},v_{\pi(1)}]\,b[v_{\pi(1)},v_{\pi(2)},v_{\pi(3)}]\\
 &=\sum_{T\in K^{(3)}}\varepsilon_T\,(a\wedge b)(T)\\
 &=\langle a\wedge b,[K]\rangle=\Hw(a).
\end{aligned}
\]
\end{proof}

The static linking formula is the appetizer; our main course is the discrete
helicity theorem in \S\ref{sec:markov}. There we use the fluxes to
introduce Markov dynamics on the primal and dual 1-skeletons. Recording the
jumps gives a discrete trajectory on each skeleton, and the volume-weighted
average of their asymptotic linking recovers $\Hw(a)$.

\section{Markov trajectories and the discrete helicity theorem}
\label{sec:markov}

This section gives a trajectory interpretation of Whitney helicity. Given
$a\in C^1(K)$ and $b=\delta a$, we use the primal and dual 1-cycles $\cpr$ and $\cdu$ of
Definition~\ref{def:flux-chains}. On each 1-skeleton, the edge fluxes
determine how a particle jumps from one vertex to the next.

\subsection{A stationary Markov process for a 1-cycle}
\label{subsec:stationary-marked-process}

We first work on the primal graph $G_{\mathrm{pr}}$. Write $\mathcal V$
for its vertex set and $c=\cpr$ for the 1-cycle defined in
\S\ref{sec:whitney-flux}. Each oriented edge $[v,w]$ of $K$ crosses the
dual face $[v,w]^*$ of $K^*$. Its coefficient in $c$ is the signed flux
\begin{equation}\label{eq:markov-signed-flux}
 J_{vw}=\int_{[v,w]^*}W_2b.
\end{equation}
We have $J_{wv}=-J_{vw}$ and set $J_{vw}=0$ for nonadjacent vertices.
Let $E$ consist of the oriented edges $e=(v,w)$ with $J_{vw}>0$, and
write $J_e=J_{vw}$ on $E$.

\begin{lemma}
\label{lem:markov-flux-balance}
If $\partial c=0$, then at every vertex $v\in\mathcal V$ the total incoming
flux equals the total outgoing flux,
\begin{equation}\label{eq:markov-flux-balance}
 \sum_{\substack{e\in E\\e^+=v}}J_e=\sum_{\substack{e\in E\\e^-=v}}J_e.
\end{equation}
\end{lemma}

\begin{proof}
Since $\partial e=[e^+]-[e^-]$, each edge entering $v$ contributes $J_e$
to the coefficient of $[v]$ in $\partial c$, and each edge leaving $v$
contributes $-J_e$. This coefficient vanishes because $\partial c=0$.
\end{proof}

Recall from \S\ref{sec:recording-vector-field} that each primal vertex
$v$ represents the dual 3-cell $v^*$ contracted to it. We therefore take
its geometric mass to be
\begin{equation}\label{eq:markov-cell-masses}
 m_v=\int_{v^*}\mu.
\end{equation}
The resulting de Rham 3-cochain is the discrete volume form on the dual
complex $K^*$, with value $m_v$ on each dual cell $v^*$~\cite{Hirani2003}.
Since the dual 3-cells partition $M$ up to their boundaries,
\begin{equation}\label{eq:total-mass-pi}
 V=\sum_{v\in\mathcal V}m_v,\qquad \pi(v)=\frac{m_v}{V}.
\end{equation}
Thus $\pi(v)$ is the probability that a uniformly distributed point lies in $v^*$.

We next determine the transition rates for our Markov process. Recall that the flow $\varphi_t$ of $X$ satisfies $\varphi_t^*\mu=\mu$
for every $t$~\cite[\S~9]{Arnold1966Hydrodynamics}. A particle initially
distributed according to $\mu/V$ therefore has the same law at every later
time. We seek the same invariance for the distribution $\pi$ determined by the cell volumes
on the primal graph. Solving the flux relation for the jump rate gives
\begin{equation}\label{eq:transition-rates}
 q_e=\frac{J_e}{m_{e^-}},\qquad \ell_v=\sum_{e^-=v}q_e.
\end{equation}
Only edges in $E$ carry jumps. In the notation of
\S\ref{sec:finite-markov-chains}, the generator therefore has off-diagonal entries
\begin{equation}\label{eq:markov-generator}
 q_{xy}=\frac{(J_{xy})_+}{m_x}\qquad (x\ne y),
\end{equation}
where $(r)_+=\max\{r,0\}$.

\begin{proposition}
\label{prop:stationarity}
The distribution $\pi$ determined by the cell volumes is stationary for the resulting Markov
chain.
\end{proposition}

\begin{proof}
By Lemma~\ref{lem:markov-flux-balance},
\[
\sum_{e^+ = v} \pi(e^-) q_e
= \frac{1}{V} \sum_{e^+ = v} J_e
= \frac{1}{V} \sum_{e^- = v} J_e
= \pi(v) \ell_v
\]
for every $v\in\mathcal V$. This is the balance condition $\pi Q=0$ of
\S\ref{sec:finite-markov-chains}, so $\pi$ is stationary.
\end{proof}

Starting at $x=X_0$, we record each jump along the edge traversed.
By Lemma~\ref{lem:markov-class-poisson}(\ref{item:markov-nonexplosion}),
the chain makes only finitely many jumps on $[0,t]$ a.s.\ for every finite $t$.
With the jump counts $N_e(t)$ of \S\ref{sec:finite-markov-chains}, the
recorded trajectory up to time $t$ is therefore the 1-chain
\begin{equation}\label{eq:markov-trajectory-chain}
J_t(x) = \sum_{e \in E} N_e(t) e.
\end{equation}
Successive intermediate endpoints cancel, yielding $\partial J_t(x) = [X_t] - [x]$.

Different parts of the graph may carry separate circulations. Let $\mathcal C$ be the set of weakly connected components of the directed graph $(\mathcal V,E)$.
In the proof of Proposition~\ref{thm:trajectory-limit}, we use Lemma~\ref{lem:markov-flux-balance} to show
that each of these components is strongly connected. An isolated vertex forms a singleton
component. For each $C\in\mathcal C$, let $E_C$ be its set of oriented edges and set
\begin{equation}\label{eq:markov-components}
m_C = \sum_{v \in C} m_v, \qquad
c_C = \sum_{e\in E_C} J_e e, \qquad
\pi_C(v) = \mathbf{1}_{\{v \in C\}} \frac{m_v}{m_C}.
\end{equation}
Since each jump follows an edge in $E$, a chain started at $v\in C$ satisfies
$X_t\in C$ for all $t\ge0$ a.s.
For a fixed vertex $v\in\mathcal V$, write $\mathbb P_v$ for the law of the
chain started at $v$ and $\mathbb E_v$ for expectation under this law.

\begin{proposition}
\label{thm:trajectory-limit}
Fix $v\in\mathcal V$ and let $C$ be its weakly connected component. Then
\begin{equation}\label{eq:markov-chain-limit}
\frac{J_t(v)}{t}\longrightarrow\frac{c_C}{m_C}
\qquad\text{strongly as }t\to\infty
\end{equation}
under $\mathbb P_v$, with respect to the $\ell^1$ norm on 1-chains.
\end{proposition}

\begin{proof}
Since $E_C$ is finite , it suffices to prove
that $N_e(t)/t\to J_e/m_C$ strongly for each $e\in E_C$.
Indeed, the coefficient limits then hold simultaneously a.s., and the bound
\[
 \Bigl\|\sum_{e\in E_C}z_e e\Bigr\|_{\ell^1}^{\,2}
 \le |E_C|\sum_{e\in E_C}|z_e|^2
\]
gives $L^2(\mathbb P_v)$ convergence of the 1-chains from that of their coefficients.

If $E_C$ is empty, then $J_t(v)=c_C=0$ for all $t$ and the conclusion follows. Hence we assume $E_C$ is nonempty. 

We first show that $C$ is strongly connected. Fix $u\in C$ and let $R$ be
the set of vertices reachable from $u$. No edge leaves $R$. Summing
the identity in Lemma~\ref{lem:markov-flux-balance} over $R$ cancels the internal fluxes, so
the total incoming flux is also zero. Since $J_e>0$ on $E$, no edge enters
$R$ either. Weak connectivity gives $R=C$. As $u$ was arbitrary, $C$ is
strongly connected.

On $C$ the generator acts by
\[
 (Qg)(w)=\sum_{\substack{e\in E_C\\ e^-=w}}q_e\bigl(g(e^+)-g(w)\bigr).
\]
Since $\pi_C(e^-)q_e=J_e/m_C$, Lemma~\ref{lem:markov-flux-balance} gives $\pi_C Q=0$.
By Lemmas~\ref{lem:markov-jump-martingales} and~\ref{lem:markov-martingale-law},
$A_e(t)/t$ and $D_g(t)/t$ both converge strongly to zero as $t\to\infty$.

Fix $u\in C$ and set $f_u(w)=\mathbf 1_{\{w=u\}}-\pi_C(u)$. Since $\pi_C$ is a probability distribution on $C$,
\[
 \sum_{w\in C}\pi_C(w)f_u(w)
 =\sum_{w\in C}\pi_C(w)\mathbf 1_{\{w=u\}}
  -\pi_C(u)\sum_{w\in C}\pi_C(w)
 =\pi_C(u)-\pi_C(u)=0.
\]
By Lemma~\ref{lem:markov-class-poisson}(\ref{item:markov-poisson}),
there exists $g_u\colon C\to\R$ such that $Qg_u=f_u$. Then
\[
 \frac1t\int_0^t\mathbf 1_{\{X_s=u\}}\,ds-\pi_C(u)
 =\frac1t\int_0^t(Qg_u)(X_s)\,ds
 =\frac{g_u(X_t)-g_u(v)-D_{g_u}(t)}{t}.
\]
Since $g_u$ is bounded on $C$, the proportion of time spent at $u$
converges strongly to $\pi_C(u)$. Taking $u=e^-$ and using $X_{s-}=X_s$ for Lebesgue-a.e.\ $s$ gives,
for each $e\in E_C$,
\[
 \frac{N_e(t)}{t}
 =q_e\frac1t\int_0^t\mathbf 1_{\{X_s=e^-\}}\,ds+\frac{A_e(t)}t
 \xrightarrow[\text{strongly}]{t\to\infty} q_e\pi_C(e^-)
 =\frac{J_e}{m_{e^-}}\frac{m_{e^-}}{m_C}=\frac{J_e}{m_C}.
\]
\end{proof}

The same construction applies to $G_{\mathrm{du}}$ with $c=\cdu$.
The dual vertex $b_T$ represents the primal tetrahedron $T$ and has
geometric mass $m_{b_T}=V_T=\int_T\mu$. Across the common face $f_{TU}$,
oriented so that its dual edge runs from $b_T$ to $b_U$, the signed flux
is $J_{b_T b_U}=b(f_{TU})$. Since $\partial\cdu=0$, the stationarity and
trajectory limit in Propositions~\ref{prop:stationarity} and
\ref{thm:trajectory-limit} apply verbatim.
We distinguish trajectories on $G_{\mathrm{pr}}$ and $G_{\mathrm{du}}$ by superscripts,
writing $J_t^{\mathrm{pr}}(x)$ and $J_s^{\mathrm{du}}(y)$ for their recorded 1-chains, respectively.

\subsection{Closed trajectories and asymptotic linking}
\label{subsec:closed-trajectories}

We still work on the primal 1-skeleton $G_{\mathrm{pr}}$. Fix a rooted
spanning tree in each component. In the notation of
\eqref{eq:markov-cycle-projection}, let $\tau(w)$ denote the tree path from
the root to $w$ in its component and $\Pi$ the corresponding projection onto 1-cycles. To close
$J_t(x)$, we return from its endpoint $X_t$ to its starting vertex $x$
along the tree.

\begin{lemma}
\label{lem:tree-closure}
Let $C$ be the component containing the initial vertex $x$. The closed
trajectory is the 1-cycle
\begin{equation}\label{eq:markov-tree-closure}
 \Gamma_t(x):=\Pi J_t(x)=J_t(x)+\tau(x)-\tau(X_t).
\end{equation}
For all $t\ge0$,
\begin{equation}\label{eq:markov-closure-bound}
 \|\Gamma_t(x)-J_t(x)\|_{\ell^1}\le |C|-1.
\end{equation}
Consequently $\Gamma_t(x)/t$ has the same strong limit as $J_t(x)/t$.
\end{lemma}

\begin{proof}
The boundary of $J_t(x)$ is $[X_t]-[x]$, so \eqref{eq:markov-cycle-projection}
gives \eqref{eq:markov-tree-closure}. The common part of the two root paths
cancels in $\tau(x)-\tau(X_t)$, leaving the simple tree path from $X_t$
to $x$. It has at most $|C|-1$ edges. Dividing \eqref{eq:markov-closure-bound}
by $t$ and applying Proposition~\ref{thm:trajectory-limit} gives the claim.

\end{proof}

The same construction and proof apply verbatim to the dual 1-skeleton $G_{\mathrm{du}}$.

Now start independent walks at a primal vertex $x$ and a dual vertex $y$.
Write $\mathbb P_{x,y}=\mathbb P_x^{\mathrm{pr}}\otimes\mathbb P_y^{\mathrm{du}}$
be their joint probability distribution and $\mathbb E_{x,y}$ for its expectation.  We next form
their normalized linking number
\[
 L_{t,s}(x,y):=\frac{\Lk\bigl(\Gamma_t^{\mathrm{pr}}(x),
                                \Gamma_s^{\mathrm{du}}(y)\bigr)}{ts},
 \qquad t,s>0.
\]
 and show that it has a strong limit.
\begin{theorem}
\label{thm:asymptotic-linking}
Let $C$ and $D$ be the components of $x$ and $y$,
respectively. Set
\begin{equation}\label{eq:markov-linking-kernel}
 \lambda(x,y)=\frac{\Lk(c_C^{\mathrm{pr}},c_D^{\mathrm{du}})}
 {m_C^{\mathrm{pr}}m_D^{\mathrm{du}}}.
\end{equation}
Under $\mathbb P_{x,y}$,
\begin{equation}\label{eq:markov-asymptotic-linking}
 L_{t,s}(x,y)\xrightarrow[\text{strongly}]{t,s\to\infty}\lambda(x,y).
\end{equation}
Moreover,
the limit is constant on $C\times D$ and independent of the spanning trees.
\end{theorem}

\begin{proof}
Under $\mathbb P_{x,y}$, Proposition~\ref{thm:trajectory-limit} and
Lemma~\ref{lem:tree-closure} give
\begin{equation}\label{eq:markov-closed-chain-limit}
 \begin{aligned}
 A_t:=\frac{\Gamma_t^{\mathrm{pr}}(x)}{t}
 &\xrightarrow[\text{strongly}]{t\to\infty} A:=\frac{c_C^{\mathrm{pr}}}{m_C^{\mathrm{pr}}},\\
 B_s:=\frac{\Gamma_s^{\mathrm{du}}(y)}{s}
 &\xrightarrow[\text{strongly}]{s\to\infty} B:=\frac{c_D^{\mathrm{du}}}{m_D^{\mathrm{du}}}.
 \end{aligned}
\end{equation}
We have $L_{t,s}(x,y)=\Lk(A_t,B_s)$ and $\lambda(x,y)=\Lk(A,B)$.
Bilinearity and \eqref{eq:markov-linking-bound} give
\[
 |L_{t,s}(x,y)-\lambda(x,y)|
 \le C_{\mathrm{lk}}\bigl(\|A_t-A\|_{\ell^1}\,\|B_s\|_{\ell^1}
                      +\|A\|_{\ell^1}\,\|B_s-B\|_{\ell^1}\bigr).
\]
On the event of probability one where both chains converge, the right-hand side tends to zero
as $t,s\to\infty$ independently. This proves a.s.\ convergence.
To prove convergence in mean square, square the same bound and use $(a+b)^2\le 2a^2+2b^2$.
Since the two walks are independent and the trees are fixed,
\[
 \mathbb E_{x,y}\!\left[\|A_t-A\|_{\ell^1}^2\,\|B_s\|_{\ell^1}^2\right]
 =\mathbb E_x^{\mathrm{pr}}\|A_t-A\|_{\ell^1}^2\,
  \mathbb E_y^{\mathrm{du}}\|B_s\|_{\ell^1}^2.
\]
The first factor tends to zero as $t\to\infty$, and the second is bounded
for all sufficiently large $s$ because $B_s\to B$ in $L^2$. The other
term in the squared bound tends to zero as well, since $A$ is deterministic
and $\mathbb E_y^{\mathrm{du}}\|B_s-B\|_{\ell^1}^2\to0$. Together with a.s.\ convergence, this proves
\eqref{eq:markov-asymptotic-linking}.
\end{proof}

Averaging $\lambda$ over the two starting vertices gives the discrete
helicity theorem.

\begin{theorem}
\label{thm:discrete-helicity}
Whitney helicity is the mass-weighted sum of the asymptotic linking numbers,
\begin{equation}\label{eq:discrete-helicity-sum}
 \Hw(a)=\sum_{v\in\mathcal V_{\mathrm{pr}}}\sum_{T\in K^{(3)}}
 m_v^{\mathrm{pr}}m_{b_T}^{\mathrm{du}}\,\lambda(v,b_T).
\end{equation}

\end{theorem}

\begin{proof}
Group the sum by component pairs $C\times D$. The weights on each pair
sum to $m_C^{\mathrm{pr}}m_D^{\mathrm{du}}$, cancelling the denominator in
\eqref{eq:markov-linking-kernel}. By bilinearity and
Theorem~\ref{thm:whitney-linking}, the resulting sum is
\[
 \sum_{C,D}\Lk(c_C^{\mathrm{pr}},c_D^{\mathrm{du}})
 =\Lk\!\left(\sum_C c_C^{\mathrm{pr}},\sum_D c_D^{\mathrm{du}}\right)
 =\Lk(\cpr,\cdu)=\Hw(a).
\]
\end{proof}

We extend the definition of $\lambda$ to $M\times M$ by setting
\[
 \lambda(x,y)=\lambda(v,b_T),\qquad
 x\in\operatorname{int}(v^*),\quad y\in\operatorname{int}(T).
\]
On the remaining set of measure zero we set $\lambda=0$.
As in \S\ref{sec:recording-vector-field}, we regard the dual cell $v^*$ as
contracted to its primal vertex $v$. Thus every starting point in the
interior of $v^*$ inherits the value at $v$ in the first argument of
$\lambda$. Likewise, every starting point in the interior of $T$ inherits
the value at $b_T$ in the second argument.

\begin{remark}
The original $\lambda$ takes a primal vertex as its first argument and a
dual vertex as its second. Exchanging $x$ and $y$ therefore selects the
dual cell containing $y$ and the primal tetrahedron containing $x$.
Since the primal and dual cell decompositions are staggered, the extension
need not satisfy $\lambda(x,y)=\lambda(y,x)$.
\end{remark}

\begin{corollary}
\label{cor:discrete-helicity-integral}
With the above extended $\lambda$ , we have
\begin{equation}\label{eq:discrete-helicity-integral}
 \Hw(a)=\int_{M\times M}\lambda(x,y)\,d\mu(x)\,d\mu(y).
\end{equation}
\end{corollary}

\begin{proof}
Each product cell $v^*\times T$ has volume $m_v^{\mathrm{pr}}m_{b_T}^{\mathrm{du}}$,
so integrating $\lambda$ over the cells gives the sum in
\eqref{eq:discrete-helicity-sum}.
\end{proof}

\subsection{Comparison with Arnold's helicity theorem}
\label{subsec:arnold-motivation}

Arnold's helicity theorem (Theorem~\ref{thm:smooth-helicity}) recovers
helicity from the average asymptotic linking of two trajectories of the
same vector field~\cite{ArnoldKhesin2021}. To see the correspondence on
the triangulation, let $\alpha$ be a potential for $X$, so that
$d\alpha=\iota_X\mu$. The cochains $a=R_1\alpha$ and
$b=\delta a=R_2(d\alpha)$ record the integrals of $\alpha$ over edges and
of $d\alpha$ over faces. Both 1-cycles in Definition~\ref{def:flux-chains}
are determined by this same $b$, one on the primal 1-skeleton and the
other on the dual 1-skeleton.

In the smooth setting, each flow segment is closed by a short geodesic
from its endpoint back to its starting point. On the mesh, a path in the
chosen spanning tree closes each recorded walk in the same way. In both
settings we divide the linking number of the two closed trajectories by
the product of their observation times. Letting both observation times tend to infinity gives $\lambda_X$ in
\eqref{eq:smooth-asymptotic-linking} and $\lambda$ in
\eqref{eq:markov-asymptotic-linking}. For fixed starting
vertices the two walks are random, but their asymptotic linking depends only
on the components containing those vertices.

In Arnold's theorem the integral runs over the two starting points in $M$.
For the discrete walks, a pair of points in
$\operatorname{int}(v^*)\times\operatorname{int}(T)$ selects the starting
vertices $v$ and $b_T$. The product cell has volume
$m_v^{\mathrm{pr}}m_{b_T}^{\mathrm{du}}$, so integrating the piecewise
constant extension of $\lambda$ gives exactly the weighted sum in
Theorem~\ref{thm:discrete-helicity}. The two helicity theorems thus take
the same integral form,
\begin{equation}\label{eq:normalized-comparison}
 \begin{aligned}
  \Hsm(X) &= \int_M \alpha\wedge d\alpha
           = \int_{M\times M}\lambda_X(x,y)\,d\mu(x)\,d\mu(y),\\
  \Hw(a)  &= \int_M W_1a\wedge d(W_1a)
           = \int_{M\times M}\lambda(x,y)\,d\mu(x)\,d\mu(y).
 \end{aligned}
\end{equation}
The second identity is the discrete helicity theorem on a fixed
triangulation. It holds for every $a\in C^1(K)$. When $a=R_1\alpha$,
the convergence of $\Hw(R_1\alpha)$ to $\Hsm(X)$ under refinement is
addressed in \S\ref{sec:convergence}.

Representing each cell by a vertex removes the distinction between starting points within that cell. Taking the long-time limit also identifies different starting cells in the same component, since the limiting edge frequencies depend only on the component of the starting vertex in the graph of jumps with positive rates. By \eqref{eq:markov-linking-kernel}, $\lambda$ is constant on each product of a primal and a dual component. If the directed graphs $(\mathcal V_{\mathrm{pr}},E_{\mathrm{pr}})$ and
$(\mathcal V_{\mathrm{du}},E_{\mathrm{du}})$ formed by the edges of positive flux
are both weakly connected, then $\lambda(x,y)=\Hw(a)/V^2$ for $(\mu\otimes\mu)$-almost every $(x,y)\in M\times M$.

\section{Extrusion and discrete transport}\label{sec:transport}

Since our construction of the Markov process was motivated by a heuristic
picture of transport, we now seek a geometric explanation for the choice of
jump rates.
In the smooth setting, a particle starting at $x$ is carried to $\varphi_t(x)$
by the flow of $X$. In the Markov model of \S\ref{sec:markov}, the position
at time $t$ of a particle starting at a vertex $v$ is the random vertex
$X_t$. This gives two ways to transport scalar functions back to their
starting points. For a smooth function $f$ and a vertex function $g$, the
corresponding operators and their generators are
\[
 \begin{aligned}
 (U_t f)(x)&=f(\varphi_t(x)), &
 \left.\frac{d}{dt}\right|_{t=0}(U_t f)(x)&=(\mathcal L_X f)(x),\\
 (P_t g)(v)&=\mathbb E_v[g(X_t)], &
 \left.\frac{d}{dt}\right|_{t=0^+}(P_t g)(v)&=(Qg)(v).
 \end{aligned}
\]
These identities suggest that $Q$ plays the same role in the discrete setting
as the Lie derivative $\mathcal L_X$ does in the smooth setting. Following
the approach of discrete exterior calculus, can we interpret $Q$ as a
suitable discrete Lie derivative on functions?

Throughout this section we impose
\hyperref[ass:volume-compatible]{Assumption~\ref*{ass:volume-compatible}}
to better capture the properties of the smooth Lie derivative and to define
the discrete inner products.
By the \emph{scalar Lie derivative}, we mean the Lie derivative acting on functions.

\subsection{Weighted extrusion and Cartan's formula}

Extrusion gives contraction a geometric meaning in every degree.
Let $\varphi_s$ be the flow of $X$ and let $S$ be an oriented $(q-1)$-dimensional
surface. Its sweep $E_X(S,t)$ is the parameterized $q$-surface given by
\[
 [0,t]\times S\longrightarrow M,\qquad (s,p)\longmapsto\varphi_s(p).
\]
We orient the parameter domain by placing the time coordinate first. For
any smooth $q$-form $\omega$,
\[
 \int_{E_X(S,t)}\omega
 =\int_0^t\int_S\varphi_s^*(\iota_X\omega)\,ds,
 \qquad
 \left.\frac{d}{dt}\right|_{t=0}
       \int_{E_X(S,t)}\omega=\int_S\iota_X\omega.
\]
This is the extrusion interpretation of contraction
in \cite[\S~8.3]{Hirani2003}.

The boundary of a smooth sweep also explains the two terms in Cartan's
formula. With the time coordinate first, we have
$\partial E_X(S,t)=\varphi_t(S)-S-E_X(\partial S,t)$.
Applying Stokes to a form on this identity and differentiating at $t=0$
separates the contribution from the swept interior from that of the swept
boundary, producing $\iota_Xd+d\iota_X$.  We retain this algebraic structure
after defining a discrete contraction.

As before, $a\in C^1(K)$, $b=\delta a$, and $\beta=W_2b$.
Orient each physical tetrahedron $T$ by the orientation of $M$. Write
$p_0,\ldots,p_3$ for the vertices of $\widehat T$ and $\lambda_{T,i}$ for
their barycentric coordinates. Let $f_i$ be the face of $T$ opposite
$\Phi_T(p_i)$, with the outward boundary orientation induced by $T$, and
set $B_{T,i}=b(f_i)$. The local calculation in Lemma~\ref{lem:local-flux-whitney} gives
\begin{equation}\label{eq:transport-reference-field}
 \widehat\beta_T:=\Phi_T^*(\beta|_T)
      =\iota_{\nu_T}\widehat\mu_T,\qquad
 \nu_T=-\frac{1}{3V_T}\sum_{i=0}^3 B_{T,i}p_i,
 \qquad
 d\lambda_{T,i}(\nu_T)=-\frac{B_{T,i}}{3V_T}.
\end{equation}
In the notation of \eqref{eq:flux-reference-field}, the reference velocity is
$\nu_T=v_T/(6V_T)$. For a shared simplex, our idea is to compute the sweep separately
in each incident tetrahedron using its own reference velocity and Whitney
form, then average the first variations.

At a vertex $v$, the weight of an incident tetrahedron $T$ is the fraction
of the mass of $v^*$ contained in $T$.

\begin{lemma}
\label{lem:vertex-cell-mass}
Let $v$ be a vertex of $K$ and let $m_v$ denote the volume of its dual cell as defined in
\eqref{eq:markov-cell-masses}. For every tetrahedron $T$ containing $v$,
\[
 \int_{v^*\cap T}\mu=\frac{V_T}{4}.
\]
Consequently,
\begin{equation}\label{eq:transport-vertex-mass}
 m_v=\int_{v^*}\mu=\frac14\sum_{T\ni v}V_T.
\end{equation}
\end{lemma}

\begin{proof}
By Assumption~\ref{ass:volume-compatible}, it suffices to work on the reference
tetrahedron $\widehat T$. The rest is Euclidean trivial things.
\end{proof}

We use these vertex weights to reproduce the identity for the smooth scalar Lie derivative
obtained by integration by parts.
Write $\beta_X=\iota_X\mu$.
Since $\mathcal L_X f=df(X)$, Stokes' theorem gives
\begin{equation}\label{eq:continuous-weak-transport}
 \int_M h\,\mathcal L_Xf\,\mu
 =\int_M h\,df\wedge\beta_X
 =-\int_M f\,dh\wedge\beta_X
 =-\int_M f\,\mathcal L_Xh\,\mu.
\end{equation}
Thus $\mathcal L_X$ is skew-symmetric in the volume pairing.

\begin{definition}
\label{def:weighted-extrusion}
Let $S$ be an oriented $(q-1)$-simplex, where $q=1,2,3$, and let
$\widehat S_T$ be its oriented reference copy in an incident tetrahedron $T$.
For signed $t$, let $E_{\nu_T}(\widehat S_T,t)$ be the oriented surface
parameterized by
\begin{equation}\label{eq:element-sweep}
 E_{T,S,t}\colon[0,1]\times\widehat S_T\longrightarrow\R^3,\qquad
 E_{T,S,t}(r,z)=z+rt\nu_T,
\end{equation}
with the time coordinate first in the product orientation.  For
$c\in C^q(K)$, extend
$\widehat\omega_{c,T}:=\Phi_T^*(W_qc|_T)$ to $\R^3$ by
its polynomial Whitney expression.  Set
\begin{equation}\label{eq:extrusion-weights}
 \theta_{S,T}=
 \begin{cases}
 \displaystyle\frac{V_T}{\sum_{U\supset S}V_U},&\dim S=0,1,\\[6pt]
 \displaystyle\frac12,&\dim S=2,
 \end{cases}
\end{equation}
and define $\iota^q\colon C^q(K)\to C^{q-1}(K)$ by
\begin{equation}\label{eq:weighted-contraction}
 (\iota^q c)(S)=\sum_{T\supset S}\theta_{S,T}
     \left.\frac{d}{dt}\right|_{t=0}
     \int_{E_{\nu_T}(\widehat S_T,t)}\widehat\omega_{c,T}.
\end{equation}
Put $\iota^0=\iota^4=0$ . Following Cartan's formula~\cite{Cartan1922}, we define
\begin{equation}\label{eq:primal-cartan}
 L^q=\delta\iota^q+\iota^{q+1}\delta,
 \qquad 0\le q\le3.
\end{equation}
\end{definition}

No vector field is prescribed on the lower-dimensional skeleton. To sweep a
simplex $S$ in this skeleton, we use the local vector field in each incident
tetrahedron and average the resulting first variations with weights
$\theta_{S,T}$. We explain the choice of these weights below.

On the top degree, we assign weight $1/2$ to each of the two tetrahedra sharing
a face so that the scalar Lie derivative on the dual complex has the same
form as its primal counterpart.

\begin{proposition}\label{prop:cartan-properties}
Let $1\le q\le 3$ and $c\in C^q(K)$. On every oriented
$(q-1)$-simplex $S$ , the  contraction operator satisfies
\[
 (\iota^q c)(S)
 =\frac{q}{3}\sum_{T\supset S}\frac{\theta_{S,T}}{V_T}
   \sum_{i:\,v_i\notin S} B_{T,i}\,c([v_i,S]).
\]
Here $v_i=\Phi_T(p_i)$, and for $S=[s_0,\ldots,s_{q-1}]$ we write
$[v_i,S]=[v_i,s_0,\ldots,s_{q-1}]$.
Thus, once the tetrahedral volumes $V_T$ are known, we can compute
$\iota^q c$ directly from the cochains $b$ and $c$.
\end{proposition}

\begin{proof}
Fix a tetrahedron $T\supset S$ and write $\tau_s(z)=z+s\nu_T$ for translation
along its reference velocity. The sweep is $(r,z)\mapsto\tau_{rt}(z)$.
With the time direction first, the extrusion formula  gives
\[
 \int_{E_{\nu_T}(\widehat S_T,t)}\widehat\omega_{c,T}
 =\int_0^t\left(\int_{\widehat S_T}
      \tau_s^*(\iota_{\nu_T}\widehat\omega_{c,T})\right)ds.
\]
The inner integral is continuous in $s$, since $\widehat\omega_{c,T}$ is
polynomial. Differentiating at $t=0$, where $\tau_0$ is the identity,
and taking the weighted sum gives
\begin{equation}\label{eq:contraction-trace}
 (\iota^q c)(S)=\sum_{T\supset S}\theta_{S,T}
    \int_{\widehat S_T}\iota_{\nu_T}\widehat\omega_{c,T}.
\end{equation}

To evaluate the integral on $\widehat S_T$, expand $\widehat\omega_{c,T}$
in Whitney basis forms $w_F$, one for each oriented $q$-face $F$ of $T$.
Only faces containing $S$ contribute. Indeed, if $F$ omits a vertex of
$S$, it has at least two vertices outside $S$. In each term of
\eqref{eq:whitney-form}, each such vertex contributes either a barycentric
coordinate or its differential. The coordinate vanishes on $\widehat S_T$.
If both occur as differentials, at most one can be evaluated on $\nu_T$,
and the other vanishes on every tangent vector to $\widehat S_T$.

Now write $S=[s_0,\ldots,s_{q-1}]$. Orient each contributing face as
$F=[v_i,S]$, so its coefficient is $c([v_i,S])$. Write $\lambda_i$ for
$\lambda_{T,i}$ and $\lambda_{s_j}$ for the coordinate of $s_j$.
Separating the term containing $\lambda_i$ in \eqref{eq:whitney-form} gives
\[
 w_{[v_i,S]}=q!\,\lambda_i\,d\lambda_{s_0}\wedge\cdots\wedge d\lambda_{s_{q-1}}
             -q\,d\lambda_i\wedge w_S.
\]
 On $\widehat S_T$,
$\lambda_i=0$ and $d\lambda_i$ annihilates tangent vectors. Thus only
the evaluation of $d\lambda_i$ on $\nu_T$ survives. By
\eqref{eq:whitney-identities} and \eqref{eq:transport-reference-field},
\[
 \int_{\widehat S_T}\iota_{\nu_T}w_{[v_i,S]}
 =-q\,d\lambda_i(\nu_T)\int_{\widehat S_T}w_S
 =-q\,d\lambda_i(\nu_T)
 =\frac{q}{3V_T}B_{T,i}.
\]
Multiplying by $c([v_i,S])$ and summing over the vertices of $T$ outside
$S$ gives $\int_{\widehat S_T}\iota_{\nu_T}\widehat\omega_{c,T}$.
Equation~\eqref{eq:contraction-trace} then sums these integrals with
weights $\theta_{S,T}$, proving the claim.
\end{proof}

\begin{proposition}\label{prop:cartan-conservation}
For the primal volume cochain $R_3\mu$ of \eqref{eq:de-rham-map}, we have
\begin{equation}\label{eq:cartan-conservation}
 \iota^3(R_3\mu)=b,\qquad \iota^2b=0,\qquad
 \delta L^q=L^{q+1}\delta\ (0\le q\le2),\qquad
 L^2b=0,\qquad L^3(R_3\mu)=0.
\end{equation}
Moreover, for any $c\in C^3(K)$ and any oriented face $f$ shared by positively
oriented tetrahedra $T$ and $U$,
\begin{equation}\label{eq:top-contraction-face}
 (\iota^3c)(f)
 =\frac{b(f)}2\left(\frac{c(T)}{V_T}+\frac{c(U)}{V_U}\right).
\end{equation}
\end{proposition}

\begin{proof}
Fix the orientation of $f$ and let $A$ be either $T$ or $U$. In the
formula of Proposition~\ref{prop:cartan-properties} with $q=3$ and $S=f$,
the only vertex of $A$ outside $f$ is its opposite vertex $v_i$.
Let $\varepsilon_{A,f} \in \{\pm1\}$ be the coefficient of $f$ in $\partial A$.
The outward face opposite $v_i$ is $\varepsilon_{A,f}f$, so
$B_{A,i}=\varepsilon_{A,f}b(f)$. The face opposite the first vertex in
$[v_i,f]$ occurs in its boundary with sign $+1$, so
$[v_i,f]=\varepsilon_{A,f}A$ as oriented tetrahedra. Hence
\[
 B_{A,i}\,c([v_i,f])=\varepsilon_{A,f}^{\,2}\,b(f)c(A)=b(f)c(A).
\]
The contribution of $A$ is therefore $b(f)c(A)/(2V_A)$, since
$\theta_{f,A}=1/2$. Adding the two contributions proves
\eqref{eq:top-contraction-face}. 

Since $(R_3\mu)(A)=V_A$, taking $c=R_3\mu$ gives $(\iota^3(R_3\mu))(f)=b(f)$
on every oriented face.

For $c=b$, the reference Whitney form in each tetrahedron is
$\widehat\beta_T=\iota_{\nu_T}\widehat\mu_T$ by
\eqref{eq:transport-reference-field}. Its contraction with $\nu_T$ vanishes
by alternation. Thus every edge integral in \eqref{eq:contraction-trace}
is zero, giving $\iota^2b=0$.

Finally, Cartan's formula \eqref{eq:primal-cartan} and $\delta^2=0$ give
\[
 \delta L^q=\delta\iota^{q+1}\delta=L^{q+1}\delta,
 \qquad 0\le q\le2.
\]
Since $b=\delta a$, we have $\delta b=0$. Substitution into Cartan's
formula now yields
\[
 L^2b=\delta(\iota^2b)+\iota^3(\delta b)=0,
 \qquad L^3(R_3\mu)=\delta(\iota^3(R_3\mu))=\delta b=0.
\]
\end{proof}

\begin{remark}
The identities in \eqref{eq:cartan-conservation} reproduce the following
properties of contraction and Lie differentiation on smooth forms.
\[
\begin{array}{@{}l@{\qquad\qquad}l@{}}
\text{Smooth forms} & \text{Simplicial cochains} \\[3pt]
\iota_X\mu=\beta_X & \iota^3(R_3\mu)=b \\
\iota_X\beta_X=\iota_X^2\mu=0 & \iota^2b=0 \\
 d\mathcal L_X=\mathcal L_Xd & \delta L^q=L^{q+1}\delta \\
\mathcal L_X\beta_X=0,\quad \mathcal L_X\mu=0
 & L^2 b=0,\quad L^3(R_3\mu)=0
\end{array}
\]

\end{remark}

\subsection{The primal scalar operator}

We use superscripts $\mathrm{pr}$ and $\mathrm{du}$ to distinguish primal and dual quantities.
Let $D^{\mathrm{pr}}=\diag(m_v)$. Motivated by the contraction interpretation in \S\ref{sec:recording-vector-field},
we weight $v$ by the volume $m_v$ of $v^*$ and define
$\langle h,f\rangle_{D^{\mathrm{pr}}}=\sum_v m_v h(v)f(v)$ to be the discrete analogue of the $L^2$ inner product. Write
$L^{\mathrm{pr}}=L^0=\iota^1\delta$.  The following calculation
connects its weak form to the signed edge flux from \eqref{eq:markov-signed-flux},
$J_{vw}=\int_{[v,w]^*}\beta$.

\begin{proposition}\label{prop:primal-scalar-stencil}
For $0$-cochains $h,f$,
\begin{align}
 \langle h,L^{\mathrm{pr}}f\rangle_{D^{\mathrm{pr}}}
   &=\int_M(W_0h)\,d(W_0f)\wedge\beta,
       \label{eq:primal-weak-transport}\\
 (L^{\mathrm{pr}}f)(v)
   &=\frac1{2m_v}\sum_{w\sim v}J_{vw}\bigl(f(w)-f(v)\bigr).
       \label{eq:primal-centered-stencil}
\end{align}
In particular, 
$L^{\mathrm{pr}}$ is skew-adjoint for $\langle\cdot,\cdot\rangle_{D^{\mathrm{pr}}}$.
\end{proposition}

\begin{proof}
Write $\lambda_v=W_0\chi_v$ for the Whitney basis function at $v$ and
$f_T=(W_0f)\circ\Phi_T$ for the affine interpolant on $\widehat T$.
Since $W_1\delta f=d(W_0f)$, the contraction formula
\eqref{eq:contraction-trace} with $m_v\theta_{v,T}=V_T/4$ gives
\[
 m_v(L^{\mathrm{pr}}f)(v)
 =\sum_{T\ni v}\frac{V_T}{4}\,df_T(\nu_T)
 =\int_M\lambda_v\,d(W_0f)\wedge\beta.
\]
For the last equality, pull the integral on each $T\ni v$ back to
$\widehat T$. By \eqref{eq:transport-reference-field}, its integrand is
$\lambda_{T,v}\,df_T(\nu_T)\,\widehat\mu_T$. The factor $df_T(\nu_T)$ is
constant, and $\int_{\widehat T}\lambda_{T,v}\,\widehat\mu_T=V_T/4$.
Multiplying by $h(v)$ and summing over $v$ proves
\eqref{eq:primal-weak-transport}, since $W_0h=\sum_v h(v)\lambda_v$.

Since $d\beta=0$, interchanging $h$ and $f$ in the weak formula and
adding gives $\int_M d\bigl((W_0h)(W_0f)\beta\bigr)=0$ by Stokes'
theorem on each tetrahedron. The face integrals cancel because the traces
agree and $M$ is closed. Thus $L^{\mathrm{pr}}$ is skew-adjoint.

For an oriented edge $[v,w]$, apply the weak formula to $\chi_v,\chi_w$
and use skew-adjointness to obtain
\[
 2m_v(L^{\mathrm{pr}}\chi_w)(v)
 =\int_M(\lambda_v\,d\lambda_w-\lambda_w\,d\lambda_v)\wedge\beta
 =\int_M W_1\chi_{[v,w]}\wedge\beta=J_{vw},
\]
where the last equality is \eqref{eq:global-flux-pairing}. If $v$ and $w$ are distinct and nonadjacent, 
no tetrahedron contains both, so the weak formula gives $(L^{\mathrm{pr}}\chi_w)(v)=0$.
The diagonal term vanishes by skew-adjointness, since
$m_v(L^{\mathrm{pr}}\chi_v)(v)=\langle\chi_v,L^{\mathrm{pr}}\chi_v\rangle_{D^{\mathrm{pr}}}=0$.
Expanding $f=\sum_w f(w)\chi_w$ therefore gives
\[
 2m_v(L^{\mathrm{pr}}f)(v)
 =\sum_{w\sim v}J_{vw}f(w)
 =\sum_{w\sim v}J_{vw}\bigl(f(w)-f(v)\bigr),
\]
where the last equality uses
Lemma~\ref{lem:markov-flux-balance}, which gives $\sum_{w\sim v}J_{vw}=0$.
\end{proof}

\begin{remark}
We can now explain the choice of weights in \eqref{eq:extrusion-weights}.
Recall the weak identity \eqref{eq:continuous-weak-transport},
\[
 \int_M h\,\mathcal L_X f\,\mu=\int_M h\,df\wedge\beta_X.
\]
We want the discrete Lie derivative to reproduce this weak formula on the
finite element space of Whitney $0$-forms, as in \eqref{eq:primal-weak-transport}.
Taking $h=\chi_v$ and matching the coefficient of $df_T(\nu_T)$ on each
tetrahedron $T\ni v$ gives the vertex weights
\[
 m_v\theta_{v,T}=\int_{\widehat T}\lambda_{T,v}\,\widehat\mu_T=\frac{V_T}{4},
 \qquad
 \theta_{v,T}=\frac{V_T}{4m_v}=\frac{V_T}{\sum_{U\ni v}V_U}.
\]
For consistency, we also weight by volume on edges. On faces,
the two weights are $1/2$ so that the dual scalar Lie derivative has the
same form as its primal counterpart, as shown in \S\ref{subsec:dual-transport}.
\end{remark}

\subsection{Dual transport from the primal--dual pairing}
\label{subsec:dual-transport}

Following Adams, we calculate the wedge of primal and dual cochains of complementary degrees
by reconstructing them on the common barycentric subdivision $\operatorname{sd}K$ and apply Whitney interpolation. Given $c\in C^{3-q}(K)$
and $\eta\in C^q(K^*)$, distribute the value on each cell equally among
the simplices of the same dimension that subdivide it, oriented by the
cell. Set the values on the remaining simplices of that degree to zero.
Denote the resulting cochains by $B^{\mathrm{pr}}c$ and $B^{\mathrm{du}}\eta$,
and let $W_{\mathrm{sd}K}$ be Whitney interpolation on $\operatorname{sd}K$.
Define the primal--dual pairing by integrating their wedge product over $M$.
Adams's formula~\cite{Adams1996} gives
\begin{equation}\label{eq:transport-primal-dual-pairing}
 P_q(c,\eta):=
 \int_M W_{\mathrm{sd}K}(B^{\mathrm{pr}}c)\wedge
        W_{\mathrm{sd}K}(B^{\mathrm{du}}\eta)
 =\frac{(3-q)!\,q!}{4!}
  \sum_{\sigma\in K^{(3-q)}}c(\sigma)\eta(\sigma^*).
\end{equation}
The coefficient $\kappa_q:=(3-q)!\,q!/4!$ is the absolute value of the
Whitney wedge integral for two complementary faces of a tetrahedron
meeting at a single vertex~\cite{Adams1996}.
Each dual cell carries the orientation fixed in \S\ref{sec:dual-complex-linking}.
In particular, a positively oriented tetrahedron has a positively oriented
dual vertex.

The continuous identity motivating dual transport is
\[
 \int_M(\mathcal L_X\alpha)\wedge\gamma
 +\int_M\alpha\wedge\mathcal L_X\gamma=0
 \qquad(\deg\alpha+\deg\gamma=3).
\]
Indeed, the smooth Lie derivative obeys the product rule, and Cartan's
formula applied to the top form $\alpha\wedge\gamma$ gives
\[
 (\mathcal L_X\alpha)\wedge\gamma
       +\alpha\wedge\mathcal L_X\gamma
 =\mathcal L_X(\alpha\wedge\gamma)
 =d\bigl(\iota_X(\alpha\wedge\gamma)\bigr).
\]
Its integral is zero by Stokes' theorem. 

We require dual transport to satisfy the corresponding identity
\begin{equation}\label{eq:dual-transport-pairing}
 P_q(L^{3-q}c,\eta)+P_q(c,L^{q,\mathrm{du}}\eta)=0
\end{equation}
for every $c\in C^{3-q}(K)$ and $\eta\in C^q(K^*)$. Since the matrix of
$P_q$ in corresponding cell bases is $\kappa_q I$, this identity
determines the dual Lie derivative uniquely.

To summarize, we define the dual operators as follows.
\begin{definition}\label{def:dual-cartan}
For $0\le q\le3$, the dual Lie derivative
$L^{q,\mathrm{du}}\colon C^q(K^*)\to C^q(K^*)$ is the operator determined
by \eqref{eq:dual-transport-pairing}. In corresponding cell bases it is
\begin{equation}\label{eq:dual-negative-transpose}
 L^{q,\mathrm{du}}:=-(L^{3-q})^\top.
\end{equation}
To realize this operator by Cartan's formula, define the dual contraction
$\iota^{q,\mathrm{du}}\colon C^q(K^*)\to C^{q-1}(K^*)$ for $1\le q\le3$ by
\begin{equation}\label{eq:dual-contraction}
 \iota^{q,\mathrm{du}}:=(-1)^{q+1}(\iota^{4-q})^\top,
 \qquad \iota^{0,\mathrm{du}}=\iota^{4,\mathrm{du}}=0.
\end{equation}
With this choice, \eqref{eq:dual-coboundary-sign} and the primal Cartan
formula \eqref{eq:primal-cartan} give
\begin{equation}\label{eq:dual-cartan}
 L^{q,\mathrm{du}}=\delta^{\mathrm{du}}\iota^{q,\mathrm{du}}
       +\iota^{q+1,\mathrm{du}}\delta^{\mathrm{du}},\qquad 0\le q\le3.
\end{equation}
\end{definition}

For a dual 0-cochain $g$, Cartan's formula reduces to
$L^{0,\mathrm{du}}g=(\iota^3)^\top\delta^{\mathrm{du}}g$.
As in Proposition~\ref{prop:primal-scalar-stencil}, we compute the coefficient
at $b_T$ using the elementary $3$-cochain $\chi_T$ of $T$. Here $U$ ranges over all tetrahedra that share a face with $T$.
Orient the common face $f_{TU}=T\cap U$ outward from $T$
and recall $J_{TU}=b(f_{TU})$. Then \eqref{eq:top-contraction-face} gives
$(\iota^3\chi_T)(f_{TU})=J_{TU}/(2V_T)$, and $f_{TU}^*$ runs from $b_T$ to $b_U$. Hence
\begin{equation}\label{eq:dual-scalar-stencil}
\begin{aligned}
 (L^{0,\mathrm{du}}g)(b_T)
 &=\sum_U(\iota^3\chi_T)(f_{TU})\,(\delta^{\mathrm{du}}g)(f_{TU}^*)\\
 &=\frac{1}{2V_T}\sum_UJ_{TU}\bigl(g(b_U)-g(b_T)\bigr).
\end{aligned}
\end{equation}
Since the dual vertex $b_T$ represents
the tetrahedron $T$, the corresponding volume inner product on $C^0(K^*)$ is
$\langle h,f\rangle_{D^{\mathrm{du}}}:=\sum_{T\in K^{(3)}}V_T\,h(b_T)f(b_T)$.
The dual mass matrix is $D^{\mathrm{du}}=\diag(V_T)$. We write $L^{\mathrm{du}}:=L^{0,\mathrm{du}}$
for brevity. We can now see why the top degree weights were chosen to be $1/2$.
With this choice, $L^{\mathrm{du}}$ and $L^{\mathrm{pr}}$ have exactly the same
form, as \eqref{eq:dual-scalar-stencil} and \eqref{eq:primal-centered-stencil} show.
In fact, the edge weights in the contraction $\iota^2$ on 2-cochains can
be varied without affecting the properties proved above.

\subsection{Transport and the Markov generators}

The jump rates in \eqref{eq:markov-generator} use the positive part of the
signed flux. The scalar Lie derivatives in \eqref{eq:primal-centered-stencil}
and \eqref{eq:dual-scalar-stencil} use the signed flux itself. To compare
these operators, write $Q$ and $L$ for the generator and scalar Lie
derivative on the same 1-skeleton.
The adjoint with respect to the inner product weighted by cell volumes is $Q^\dagger=D^{-1}Q^\top D$, with
$D=D^{\mathrm{pr}}=\diag(m_v)$ for $Q^{\mathrm{pr}}$ and
$D=D^{\mathrm{du}}=\diag(V_T)$ for $Q^{\mathrm{du}}$.

Suppose a continuous-time Markov chain on a finite set has a stationary
distribution $\pi$ with $\pi(x)>0$ for every state $x$. Its generator admits
the transport--diffusion decomposition
\[
 Q=S+\Delta_{\mathrm{rev}},\qquad S:=\frac{Q-Q^\dagger}{2},\qquad \Delta_{\mathrm{rev}}:=\frac{Q+Q^\dagger}{2},
\]
where $Q^\dagger$ is the adjoint in $L^2(\pi)$ ~\cite[\S~3.7]{Norris1997}. The skew-adjoint
transport part $S$ changes sign under time reversal. The diffusion part $\Delta_{\mathrm{rev}}$
averages the forward and reversed jump rates and is reversible with respect
to $\pi$. For our two chains, $\pi(x)=m_x/V>0$ because each vertex represents a cell
of positive volume $m_x$. Dividing the inner product weighted by cell volumes by $V$ gives the
$L^2(\pi)$ inner product, so their adjoints agree. Averaging the rates across an edge gives
$\bigl((J_{xy})_++(-J_{xy})_+\bigr)/(2m_x)=|J_{xy}|/(2m_x)$, hence
\begin{equation}\label{eq:reversible-part}
 (\Delta_{\mathrm{rev}} f)(x)=\frac{1}{2m_x}\sum_y |J_{xy}|\bigl(f(y)-f(x)\bigr),
\end{equation}
where $y$ ranges over the vertices joined to $x$ by an edge.

\begin{theorem}\label{thm:scalar-transport}
The scalar Lie derivative on each 1-skeleton is the transport part $S$ of its
Markov generator,
\begin{equation}\label{eq:both-scalar-transport}
 L^{\mathrm{pr}}=\frac{Q^{\mathrm{pr}}-(Q^{\mathrm{pr}})^\dagger}{2},
 \qquad
 L^{\mathrm{du}}=\frac{Q^{\mathrm{du}}-(Q^{\mathrm{du}})^\dagger}{2}.
\end{equation}
\end{theorem}

\begin{proof}
For adjacent vertices $x,y$, the generator and its adjoint have entries
\[
 q_{xy}=\frac{(J_{xy})_+}{m_x},\qquad
 (Q^\dagger)_{xy}=\frac{m_y}{m_x}q_{yx}=\frac{(-J_{xy})_+}{m_x}.
\]
Since $\pi(x)=m_x/V$ is stationary by Proposition~\ref{prop:stationarity},
the rows of $Q^\dagger$ sum to zero, as do those of $Q$. Taking half the
difference of the two entries therefore gives
\[
 (Sf)(x)
 =\frac{1}{2m_x}\sum_y J_{xy}\bigl(f(y)-f(x)\bigr).
\]
At a primal vertex $v$, this is \eqref{eq:primal-centered-stencil}.
At the dual vertex $b_T$, it is \eqref{eq:dual-scalar-stencil}
with $m_{b_T}=V_T$ and $J_{b_T b_U}=J_{TU}$.
\end{proof}

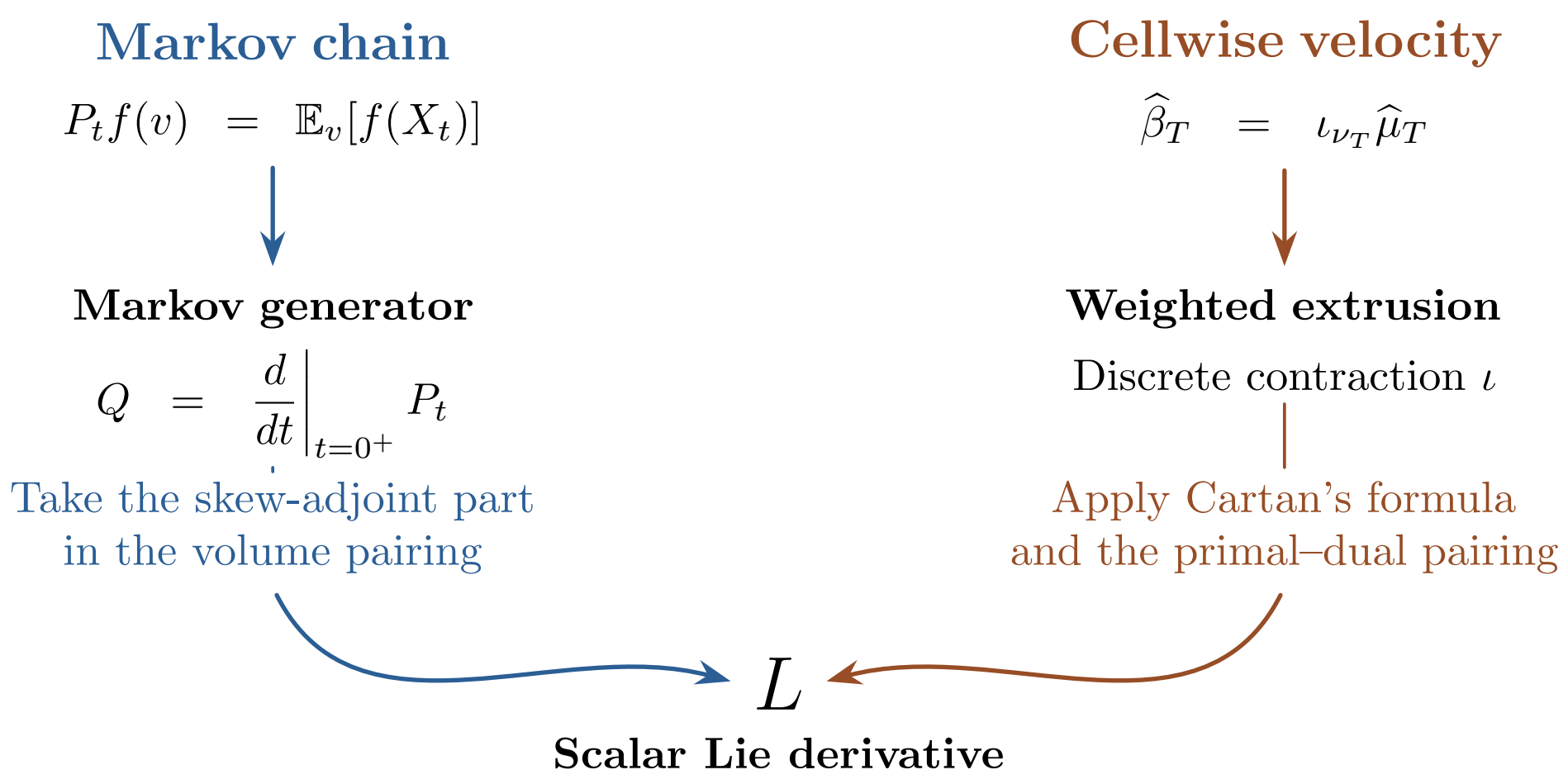
\begin{figure}[H]
\centering
\begin{tikzpicture}[
  x=1cm,y=1cm, line cap=round, line join=round,
  every node/.style={align=center,inner sep=3pt},
  stage/.style={text width=6.0cm,font=\normalsize},
  action/.style={font=\normalsize,text width=5.5cm},
  flow/.style={line width=1pt,-{Stealth[length=2.8mm,width=2.0mm]},
               shorten <=2pt,shorten >=2pt}]
  \definecolor{markovink}{RGB}{43,94,149}
  \definecolor{extrusionink}{RGB}{151,77,37}

  \node[stage] (markov) at (-4.05,0) {
    {\large\bfseries\color{markovink} Markov chain}\\[5pt]
    $P_t f(v)=\mathbb E_v[f(X_t)]$};
  \node[stage] (extrusion) at (4.05,0) {
    {\large\bfseries\color{extrusionink} Cellwise velocity}\\[5pt]
    $\widehat\beta_T=\iota_{\nu_T}\widehat\mu_T$};

  \node[stage,anchor=north] (generator) at (-4.05,-1.53) {
    {\bfseries Markov generator}\\[4pt]
    $\displaystyle Q=\left.\frac{d}{dt}\right|_{t=0^+}P_t$};
  \node[stage,anchor=north] (contraction) at (4.05,-1.53) {
    {\bfseries Weighted extrusion}\\[4pt]
    Discrete contraction $\iota$};
  \draw[flow,markovink] (markov.south) -- (generator.north);
  \draw[flow,extrusionink] (extrusion.south) -- (contraction.north);

  \node[action,text=markovink] (transport) at (-4.05,-3.55) {
    Take the skew-adjoint part\\in the volume pairing};
  \node[action,text=extrusionink] (cartan) at (4.05,-3.55) {
    Apply Cartan's formula\\and the primal--dual pairing};
  \draw[markovink,line width=.7pt] (generator.south) -- (transport.north);
  \draw[extrusionink,line width=.7pt] (contraction.south) -- (cartan.north);

  \node[font=\LARGE,inner sep=3pt] (lie) at (0,-4.80) {$L$};
  \draw[flow,markovink] (transport.south) to[out=-63,in=165] (lie.west);
  \draw[flow,extrusionink] (cartan.south) to[out=-117,in=15] (lie.east);
  \node[anchor=north,font=\normalsize\bfseries] at (lie.south) {Scalar Lie derivative};
\end{tikzpicture}
\caption{Two paths, one destination.}
\label{fig:markov-extrusion-lie}
\end{figure}

The scalar Lie derivative $\mathcal L_X$
is formally skew-adjoint in the volume pairing by \eqref{eq:continuous-weak-transport}.
This motivates taking the skew-adjoint part $(Q-Q^\dagger)/2$ of the Markov
generator. Figure~\ref{fig:markov-extrusion-lie} shows the correspondence at a
glance. Taking this part on the left gives the same scalar Lie derivative $L$
as weighted extrusion and Cartan's formula on the right.

\begin{remark}
Hirani uses the contraction formula 
$\iota_X\alpha=\star\bigl((\star\alpha)\wedge X^\flat\bigr)$ along with Cartan's formula to define the discrete Lie derivative~\cite{Hirani2003}.
He also discusses Bossavit's geometric interpretation of contraction by extrusion
and works out the case of a field aligned with an edge~\cite[\S~8.3]{Hirani2003}.
Mullen et al. used finite-volume approximations of extrusions to discretize
Lie advection~\cite{MullenEtAl2011}. Here we realize the extrusion idea for
cellwise velocities by averaging the contributions from the tetrahedra incident
to each simplex. Even on scalar cochains, the resulting Lie derivative can
differ from the primal analogue of Hirani's algebraic construction with
circumcentric Hodge stars and the Whitney wedge on a simplicial dual.
For the dual operator we draw on the adjoint formulation of Heumann and
Hiptmair~\cite{HeumannHiptmair2011}. Moving the Lie derivative
between complementary forms in the integrated wedge pairing changes its sign.
We impose the same relation on primal and dual cochains in \eqref{eq:dual-transport-pairing}.
\end{remark}

\section{Convergence}
\label{sec:convergence}

We now refine the triangulation $K$. The two questions are
whether the helicities of \S\ref{sec:discrete-helicity} recover smooth helicity
and whether the primal and dual Lie derivatives of \S\ref{sec:transport}
recover $\mathcal L_X$ on differential forms. Hirani takes Stokes' theorem and naturality under pullback as guides for
constructing discrete operators on a fixed mesh. He treats convergence of
these operators under refinement as a separate question that requires
further analysis~\cite{Hirani2003}.

We compare cochains with smooth forms using Whitney interpolation, whose
approximation properties under refinement were established by
Dodziuk~\cite{Dodziuk1976}. Whitney forms constitute the lowest-order
de Rham complex in finite element exterior calculus (FEEC)~\cite{ArnoldFalkWinther2006}.
Wilson's first-order estimate for the cochain product~\cite{Wilson2007}
gives the error bound for Whitney helicity. For cup helicity we also
compare the cup product with the discrete wedge.
For transport, Heumann and Hiptmair~\cite{HeumannHiptmair2011} used finite
element differential forms to approximate the Lie derivative in advection
problems. Here the velocity enters through its face fluxes, and we study
how accurately the resulting operator approximates $\mathcal L_X\omega$.

The tetrahedra of $K$ are the
individual mesh elements, and their vertices and edges form the 1-skeleton.
These elements cover $M$, and any two distinct
tetrahedra that intersect share a triangular face, an edge, or a vertex.
We index each triangulation $K_h$ by its mesh diameter
$h:=\max_{T\in K_h^{(3)}} h_T$, where $h_T:=\operatorname{diam}(T)$.
The subscript $h$ on $a_h$ and $b_h$ indicates that these cochains are
defined on the mesh of diameter $h$. On $K_h$ we keep the notation $R_p,W_p,\delta$ of \S\ref{sec:preliminaries}
and $L^{\mathrm{pr}}$ of \S\ref{sec:transport}. For a potential $\alpha$, the
cochains on this mesh are
\[
 a_h=R_1\alpha,\qquad b_h=\delta a_h=R_2(d\alpha).
\]
 
Write $d(x,y)$ for the geodesic distance between $x$ and $y$.
For a triangular face $F$ and a tetrahedron $T$, let
$\operatorname{area}(F)$ and $\operatorname{vol}(T)$ denote their area
and volume, respectively.

\subsection{Convergence of two discrete helicities}
\label{sec:helicity-convergence}

In this subsection, $(K_h)$ is a family of triangulations that refines a fixed finite smooth triangulation $K_0$
by compatible affine subdivisions of its reference tetrahedra, with $h\to0$.
A tetrahedron $R$ in the subdivision of the reference tetrahedron for $T_0$
is called a \emph{reference subcell}. Its image is $T=\Phi_{T_0}(R)$.

\begin{definition}\label{def:controlled-refinement}
We call the refinement \emph{controlled} if its fullness is uniformly
bounded below~\cite[\S~4]{Wilson2007}, that is, if some
$\Theta_0>0$ independent of $h$ satisfies
\begin{equation}\label{eq:mesh-fullness}
 \inf_{T\in K_h^{(3)}}\frac{\operatorname{vol}(T)}{\eta_h^3}\ge\Theta_0,
 \qquad \eta_h:=\max_{[v,w]\in K_h^{(1)}}d(v,w)\le h.
\end{equation}
\end{definition}

The estimate for Whitney helicity follows from convergence of the discrete
wedge. The cup product requires a separate comparison using algebraic topology because the
potential cochain $a_h$ need not be closed.

\begin{proposition}\label{prop:whitney-product-convergence}
Let $(K_h)$ be a controlled refinement and fix smooth forms
$\omega\in\Omega^p(M)$ and $\gamma\in\Omega^q(M)$ with $p+q\le3$.
There is a constant $C_{\omega,\gamma}$ independent of $h$ such that
\begin{equation}\label{eq:wilson-product-estimate}
 \|W_{p+q}(R_p\omega\wedge R_q\gamma)-\omega\wedge\gamma\|_{L^2(M)}
 \le C_{\omega,\gamma}\,h.
\end{equation}
\end{proposition}
\begin{proof}
This follows from Wilson's estimate~\cite[\S~5]{Wilson2007}.
\end{proof}

\begin{lemma}\label{lem:local-product-homotopy}
On every finite simplicial complex $K$ of dimension at most $3$, there are
bilinear maps $\mathcal K:C^p(K)\times C^q(K)\to C^{p+q-1}(K)$ such that
\begin{equation}\label{eq:product-homotopy}
 a\smile b-a\wedge b=\delta\mathcal K(a,b)+\mathcal K(\delta a,b)
                       +(-1)^p\mathcal K(a,\delta b).
\end{equation}
For $p=0$ or $q=0$,  $\mathcal K(a,b)=0$, and for $p,q\ge1$ and every
$(p+q-1)$-simplex $\sigma$, we have
\begin{equation}\label{eq:product-homotopy-bound}
 |\mathcal K(a,b)(\sigma)|\le 496
 \max_{\substack{F\subseteq\sigma\\\dim F=p}}|a(F)|\,
 \max_{\substack{G\subseteq\sigma\\\dim G=q}}|b(G)|.
\end{equation}
\end{lemma}
\begin{proof}
For each product $\diamond\in\{\smile,\wedge\}$ and each degree $r$, define
its chain diagonal
\[
 \Delta_{\diamond}\colon C_r(K)\longrightarrow
 \bigoplus_{i+j=r} C_i(K)\otimes C_j(K)
\]
by $\langle a\otimes b,\Delta_{\diamond}c\rangle=\langle a\diamond b,c\rangle$
for $c\in C_r(K)$, $a\in C^i(K)$ and $b\in C^j(K)$ with $i+j=r$.
Since $K$ is finite, transposing the  product formulas gives the
two diagonals, and transposing the Leibniz rules shows that they are chain
maps~\cite{SchubelEtAl2024}.
Both send each vertex $v$ to $v\otimes v$.
Set $\Delta=\Delta_{\smile}-\Delta_{\wedge}$ in each degree.

Write a simplex as $T=[v_0,\ldots,v_m]$ in increasing vertex order.
Let $\pi_T$ send every vertex to $v_0$ and vanish in positive degrees. Let $s_T$ be the cone
operator that prepends $v_0$ to a face and is zero if the face already
contains $v_0$. Then $\partial s_T+s_T\partial=\mathrm{id}-\pi_T$ and the $\ell^1$--operator norm
$\|s_T\|\le1$. On homogeneous tensors put
\[
 S_T(x\otimes y)=s_Tx\otimes y+(-1)^{\deg x}\pi_Tx\otimes s_Ty.
\]
With $\partial_\otimes(x\otimes y)=\partial x\otimes y+
(-1)^{\deg x}x\otimes\partial y$, this gives
\[
 \partial_\otimes S_T+S_T\partial_\otimes=\mathrm{id}-\pi_T\otimes\pi_T,
 \qquad \|S_T\|\le2.
\]

Define $P_K$ by induction on simplex dimension, starting with $P_K(v)=0$.
For $m>0$, once $P_K$ is defined on the proper faces of $T$, put
\[
 z_T=\Delta(T)-P_K(\partial T),\qquad P_K(T)=S_Tz_T.
\]
The induction hypothesis gives $\partial_\otimes z_T=0$, and
$(\pi_T\otimes\pi_T)z_T=0$ because $z_T$ has total degree $m>0$.
Hence $\partial_\otimes P_K(T)=z_T$. This proves
$\Delta=\partial_\otimes P_K+P_K\partial$, with
$P_K(T)\in C_\bullet(T)\otimes C_\bullet(T)$.

The cup diagonal has $m+1$ terms of coefficient $1$ on $T$.
In bidegree $(p,m-p)$, the Whitney diagonal has $(m+1)\binom{m}{p}$
nonzero face pairs, obtained by choosing their common vertex and the
other $p$ vertices of the first face. Each coefficient has magnitude
$p!(m-p)!/(m+1)!$ by Wilson's local product formula~\cite{Wilson2007}.
Thus each bidegree contributes $1$ to $\|\Delta_{\wedge}(T)\|_{\ell^1}$,
and $\|\Delta(T)\|_{\ell^1}\le2(m+1)$. Since $\partial T$ has $m+1$ faces,
the construction gives $\|P_K(T)\|_{\ell^1}\le A_m$, where
\[
 A_0=0,\qquad A_m=2(m+1)(2+A_{m-1}).
\]
Thus $A_3=496$ bounds $\|P_K(T)\|_{\ell^1}$ for every simplex of dimension at most $3$.

Finally, define $\mathcal K(a,b)(c)=\langle a\otimes b,P_K(c)\rangle$.
Pairing the chain homotopy identity with $a\otimes b$ gives
\eqref{eq:product-homotopy}, since
$\delta_\otimes(a\otimes b)=\delta a\otimes b+(-1)^pa\otimes\delta b$.
The support of $P_K(\sigma)$ on pairs of faces of $\sigma$ and its
$\ell^1$ bound give \eqref{eq:product-homotopy-bound}.
\end{proof}

\begin{proposition}\label{prop:cup-whitney-comparison}\label{prop:refinement-geometry}
For a $1$-cochain $a$ on a closed oriented triangulated $3$-manifold,
\begin{equation}\label{eq:discrete-helicity-comparison}
 |\Hc(a)-\Hw(a)|\le 496\sum_{T\in K^{(3)}}
       \max_{\substack{F\subset T\\\dim F=2}}|(\delta a)(F)|^2.
\end{equation}
On a controlled refinement, for $a_h=R_1\alpha$ with $\alpha$ smooth, we have
\begin{equation}\label{eq:smooth-cup-whitney-comparison}
 |\Hc(a_h)-\Hw(a_h)|\le Ch\|d\alpha\|_{L^\infty(M)}^2
\end{equation}
where the constant $C$ is independent of both $h$ and the vertex order.
\end{proposition}
\begin{proof}
Evaluate \eqref{eq:product-homotopy} on $[K]$ with $b=\delta a$.
The first term vanishes because $\partial[K]=0$ and the last because $\delta b=0$.
Thus $\Hc(a)-\Hw(a)=\langle\mathcal K(b,b),[K]\rangle$, and the local bound
gives \eqref{eq:discrete-helicity-comparison}.

By a slight abuse of notation, we use $C$ for positive constants that may vary from line to line
but are independent of $h$ and the vertex order.
For each triangular face $F$, Wilson's estimate gives
$\operatorname{area}(F)\le C\eta_h^2$~\cite[\S~4]{Wilson2007}. Since $b_h=R_2(d\alpha)$,
\[
 |b_h(F)|=\left|\int_F d\alpha\right|
 \le \operatorname{area}(F)\,\|d\alpha\|_{L^\infty(M)}
 \le C\eta_h^2\|d\alpha\|_{L^\infty(M)}.
\]
Squaring and taking the maximum over all faces gives the first inequality below.
For the second, each tetrahedron has volume at least $\Theta_0\eta_h^3$
by \eqref{eq:mesh-fullness}. Their interiors are disjoint and they cover $M$,
so $\#K_h^{(3)}\Theta_0\eta_h^3\le\sum_{T\in K_h^{(3)}}\operatorname{vol}(T)=\operatorname{vol}(M)$.
Together these estimates give
\begin{equation}\label{eq:local-mesh-estimates}
 \max_{F\in K_h^{(2)}}|b_h(F)|^2
 \le C\eta_h^4\|d\alpha\|_{L^\infty(M)}^2,\qquad
 \#K_h^{(3)}\le\frac{\operatorname{vol}(M)}{\Theta_0\eta_h^3}.
\end{equation}
Substitution in
\eqref{eq:discrete-helicity-comparison} gives
\eqref{eq:smooth-cup-whitney-comparison} since $\eta_h\le h$.
\end{proof}

\begin{theorem}\label{thm:helicity-convergence}
On a controlled refinement, for every smooth potential $\alpha$ and
$a_h=R_1\alpha$, both helicities converge to the smooth helicity,
\begin{equation}\label{eq:helicity-convergence}
 \lim_{h\to0}\Hw(a_h)=\lim_{h\to0}\Hc(a_h)=\Hsm(\alpha).
\end{equation}
Moreover, there is a constant $C_\alpha$, independent of $h$ and the vertex order, such that
\begin{equation}\label{eq:helicity-convergence-rate}
 |\Hw(a_h)-\Hsm(\alpha)|+|\Hc(a_h)-\Hsm(\alpha)|\le C_\alpha h.
\end{equation}
\end{theorem}
\begin{proof}
Since $R_3W_3$ is the identity on $3$-cochains,
\[
 \Hw(a_h)=\langle a_h\wedge b_h,[K_h]\rangle
         =\int_M W_3(R_1\alpha\wedge R_2d\alpha).
\]
Write $\omega_h:=W_3(R_1\alpha\wedge R_2d\alpha)-\alpha\wedge d\alpha=r_h\mu$.
We use $\mu$ for the scalar $L^2(M)$ norm. Applying Cauchy--Schwarz gives
\[
\begin{aligned}
 |\Hw(a_h)-\Hsm(\alpha)|
 &=\left|\int_M r_h\,\mu\right|\\
 &\le \left(\int_M |r_h|^2\,\mu\right)^{1/2}
       \left(\int_M 1^2\,\mu\right)^{1/2}
 = V^{1/2}\|r_h\|_{L^2(M)}.
\end{aligned}
\]
Since $\mu$ is smooth and nowhere zero on compact $M$, taking the scalar
coefficient of a $3$-form relative to $\mu$ is bounded in $L^2(M)$. Hence
\eqref{eq:wilson-product-estimate} with $(\omega,\gamma)=(\alpha,d\alpha)$ gives
$\|r_h\|_{L^2(M)}\le C_M\|\omega_h\|_{L^2(M)}\le C_\alpha h$. Proposition~\ref{prop:cup-whitney-comparison} then gives the cup
estimate for every ordering of the vertices.
\end{proof}

\begin{remark}
The same approximation also gives convergence of the kinetic energy in
\eqref{eq:whitney-energy}. Let $\mathbf{1}\colon M\to\R$ denote the constant function
$\mathbf{1}(x)=1$ for every $x\in M$. Its de Rham $0$-cochain $R_0\mathbf{1}$
takes the value $1$ at every vertex and is the unit for the discrete wedge.
For $\beta_X=\iota_X\mu$ and $b_h=R_2\beta_X$, applying
\eqref{eq:wilson-product-estimate} to $(\mathbf{1},\beta_X)$ gives $\|W_2b_h-\beta_X\|_{L^2(M)}=O(h)$. Hence
$E_{K_h}(b_h)=\tfrac12\|W_2b_h\|_{L^2(M)}^2=\mathcal E(X)+O(h)$.
The convergence of kinetic energy and Whitney helicity
shows that both sides of \eqref{eq:whitney-energy} converge under mesh refinement
to their counterparts in Arnold's energy inequality.
\end{remark}

\subsection{Convergence of discrete Lie derivative}
\label{sec:scalar-convergence}

On   functions the Lie derivative is the directional derivative
$\mathcal L_X f=df(X)$.
On each mesh $K_h$, we use the face fluxes $b_h(F)=\int_F\iota_X\mu$ to
construct $L^{\mathrm{pr}}$ as in \S\ref{sec:transport}.
We ask how the discrete Lie derivative relates to its continuous counterpart.

Let $h\to0$ and suppose the meshes satisfy Assumption~\ref{ass:volume-compatible}.
In coordinates of a fixed finite atlas, assume that each cell map $\Phi_T$
satisfies
\begin{equation}\label{eq:scalar-map-bounds}
 \|D\Phi_T\|_\infty\le C_Gh_T,\qquad
 \|(D\Phi_T)^{-1}\|_\infty\le C_Gh_T^{-1},\qquad
 \|D^2\Phi_T\|_\infty\le C_Gh_T^2.
\end{equation}
We fix the mesh regularity bound $C_G$ throughout this subsection.
It applies to every $T$ and $h$ in charts containing $T$ with uniform metric bounds.

\begin{remark}
Ern and Guermond use the scaling in \eqref{eq:scalar-map-bounds} to
analyze interpolation on nonaffine meshes~\cite{ErnGuermond2021}. To see why these powers arise
here, consider a shape-regular affine subdivision of a fixed smooth
triangulation. On a subcell of a coarse tetrahedron $T_0$, write
$\Phi_T(z)=\Phi_{T_0}(a_T+A_Tz)$. Shape regularity gives
$\|A_T\|=O(h_T)$ and $\|A_T^{-1}\|=O(h_T^{-1})$. In the second derivative,
\[
 D^2\Phi_T(z)[u,v]=D^2\Phi_{T_0}(a_T+A_Tz)[A_Tu,A_Tv],
\]
each argument contributes one factor of $h_T$. The coarse cell maps are
fixed and finite in number, so the constants are uniform even when the
subcells have different diameters.

We use these bounds to estimate the directional derivative on the
reference tetrahedron. For $f\in C^2(M)$, they give
$\|D^2(f\circ\Phi_T)\|_\infty=O(h_T^2)$. Thus affine interpolation
approximates the reference gradient to order $h_T^2$. Evaluating the
gradient error on $\nu_T$, whose magnitude is $O(h_T^{-1})$, gives the
$O(h_T)$ contribution to the scalar Lie derivative estimate below.
\end{remark}

For a $0$-cochain $c$, write $c_T=(W_0c)\circ\Phi_T$ for its affine
interpolant on $\widehat T$, as in \S\ref{sec:transport}. Recall that $L^{\mathrm{pr}}$  is defined as 
\begin{equation}\label{eq:scalar-convergence-operator}
 (L^{\mathrm{pr}}c)(v)=\sum_{T\ni v}\theta_{v,T}\,dc_T(\nu_T),
 \qquad \theta_{v,T}=\frac{V_T}{\sum_{U\ni v}V_U}.
\end{equation}

\begin{theorem}\label{thm:scalar-convergence}

For every $f\in C^2(M)$ we have
\begin{equation}\label{eq:scalar-first-order}
 \|W_0L^{\mathrm{pr}}R_0f-\mathcal L_X f\|_{L^\infty(M)}
 \le Ch\|X\|_{C^1(M)}\|f\|_{C^2(M)}.
\end{equation}
Here $C$ depends only on $M$.
\end{theorem}
\begin{proof}
Fix a tetrahedron $T$ and a point $x\in T$. Whitney interpolation gives
\[
 (W_0L^{\mathrm{pr}}R_0f)(x)
 =\sum_{v\in T^{(0)}}\lambda_{T,v}(\Phi_T^{-1}(x))\,(L^{\mathrm{pr}}R_0f)(v).
\]
For each of the four vertices $v$ of $T$, we have $d(x,v)\le h_T$.
Since $\mathcal L_X f=df(X)$, it follows that
\[
 |\mathcal L_X f(v)-\mathcal L_X f(x)|
 \le h_T\|d(\mathcal L_X f)\|_{L^\infty(M)}
 \le Ch_T\|X\|_{C^1(M)}\|f\|_{C^2(M)}.
\]
The barycentric weights are nonnegative and sum to one. It therefore
suffices to bound $|(L^{\mathrm{pr}}R_0f)(v)-\mathcal L_X f(v)|$ by
$Ch\|X\|_{C^1(M)}\|f\|_{C^2(M)}$ at the four vertices of each tetrahedron.

On $\widehat T$ write $u=f\circ\Phi_T$ and
$Y=(D\Phi_T)^{-1}(X\circ\Phi_T)$. For $z\in\widehat T$, the chain rule
and the identity $\mathcal L_X f=df(X)$ give
\[
 du_z(Y(z))
 = df_{\Phi_T(z)}\!\left(D\Phi_T(z)Y(z)\right)
 = df_{\Phi_T(z)}\!\left(X(\Phi_T(z))\right)
 = (\mathcal L_X f)(\Phi_T(z)).
\]
By the chain rule and
\eqref{eq:scalar-map-bounds},
\begin{equation}\label{eq:scalar-pullback-bounds}
 \begin{aligned}
 \|du\|_\infty&\le Ch_T\|f\|_{C^2},&
 \|D^2u\|_\infty&\le Ch_T^2\|f\|_{C^2},\\
 \|Y\|_\infty&\le Ch_T^{-1}\|X\|_{C^1},&
 \|DY\|_\infty&\le C\|X\|_{C^1}.
 \end{aligned}
\end{equation}
For the last bound, put $A=D\Phi_T$. Differentiating $A^{-1}A=I$ gives
$\|D(A^{-1})\|_\infty\le\|A^{-1}\|_\infty^2\|DA\|_\infty\le C$.
The other term in $DY$ is bounded by
$\|A^{-1}\|_\infty\|DX\|_\infty\|A\|_\infty\le C\|X\|_{C^1}$.

To compare $\nu_T$ with $Y$, take $\widehat T=[0,e_1,e_2,e_3]$. For
$i=1,2,3$, the face $\widehat F_i$ opposite $e_i$ lies in $z_i=0$, has
outward unit normal $-e_i$, and has area $1/2$. By
Assumption~\ref{ass:volume-compatible}, the pullback of $\iota_X\mu$ is
$6V_T\,\iota_Y dV_{\widehat T}$. Thus the outward flux $B_{T,i}$ and
the $i$th component of the velocity in \eqref{eq:transport-reference-field}
satisfy
\[
 B_{T,i}=-6V_T\int_{\widehat F_i}Y_i\,dA,\qquad
 (\nu_T)_i=-\frac{B_{T,i}}{3V_T}
          =\frac{1}{\operatorname{area}(\widehat F_i)}
             \int_{\widehat F_i}Y_i\,dA.
\]
Since $(\nu_T)_i$ is an average of $Y_i$ on $\widehat F_i$, we have
$|(\nu_T)_i|\le\|Y_i\|_\infty$ and hence
$|\nu_T|\le\sqrt{3}\,\|Y\|_\infty\le Ch_T^{-1}\|X\|_{C^1}$ by
\eqref{eq:scalar-pullback-bounds}.
For $z\in\widehat T$ and $w\in\widehat F_i$, the segment from $z$ to $w$
lies in $\widehat T$, so $|Y_i(w)-Y_i(z)|\le\|DY\|_\infty|w-z|$.
Averaging this difference over the face gives
\[
\begin{aligned}
 |(\nu_T)_i-Y_i(z)|
 &=\left|\frac{1}{\operatorname{area}(\widehat F_i)}
   \int_{\widehat F_i}\bigl(Y_i(w)-Y_i(z)\bigr)\,dA(w)\right|\\
 &\le \operatorname{diam}(\widehat T)\,\|DY\|_\infty.
\end{aligned}
\]
Since $\widehat T$ has fixed diameter, taking the Euclidean norm of the
three components and using \eqref{eq:scalar-pullback-bounds} yields
\begin{equation}\label{eq:scalar-reconstruction-bounds}
 |\nu_T-Y(z)|\le C\|X\|_{C^1},\qquad z\in\widehat T.
\end{equation}

Set $c=R_0f$. On $\widehat T$, the interpolant $c_T$ in
\eqref{eq:scalar-convergence-operator} is $c_T=(W_0R_0f)\circ\Phi_T=Iu$,
where $Iu$ denotes the affine interpolant of $u$. Thus $d(Iu)(\nu_T)$
enters $(L^{\mathrm{pr}}R_0f)(v)$ with weight $\theta_{v,T}$ at each
vertex $v$ of $T$. Since the weights sum to one,
\[
 (L^{\mathrm{pr}}R_0f)(v)-\mathcal L_X f(v)
 =\sum_{U\ni v}\theta_{v,U}\bigl(dc_U(\nu_U)-\mathcal L_X f(v)\bigr).
\]
We estimate the summand for $U=T$. Fix a vertex $z$ of $\widehat T$ and
write $p_0=0$ and $p_j=e_j$ for $j=1,2,3$. For $j=0,1,2,3$, define
\[
 r_j:=u(p_j)-u(z)-du_z(p_j-z).
\]
Thus $r_j$ is the error at $p_j$ in the first-order Taylor approximation
of $u$ about $z$. Taylor's theorem along the segment from $z$ to $p_j$
gives $|r_j|\le \tfrac12\|D^2u\|_\infty |p_j-z|^2$.
Since $Iu$ is affine and agrees with $u$ at $0$ and $e_j$, its derivative
in the direction $e_j$ is $u(e_j)-u(0)$. Subtracting the Taylor formulas
at $e_j$ and $0$ therefore gives
\[
 (d(Iu)-du_z)(e_j)=u(e_j)-u(0)-du_z(e_j)=r_j-r_0,
 \qquad j=1,2,3.
\]
All four vertices lie within $\operatorname{diam}(\widehat T)$ of $z$.
The bounds on $r_j$ and $r_0$ thus give
$|d(Iu)-du_z|\le C\|D^2u\|_\infty$.
At $v=\Phi_T(z)$ we obtain
\[
\begin{aligned}
 |d(Iu)(\nu_T)-\mathcal L_X f(v)|
 &\le |d(Iu)-du_z|\,|\nu_T|+|du_z|\,|\nu_T-Y(z)|\\
 &\le Ch_T\|X\|_{C^1}\|f\|_{C^2}.
\end{aligned}
\]

\end{proof}

We next prove that the primal discrete Lie derivative converges weakly to the
continuous Lie derivative on forms of every degree.

\begin{theorem}\label{thm:weak-lie-convergence}
Suppose that Assumption~\ref{ass:volume-compatible} and
\eqref{eq:scalar-map-bounds} hold. For $0\le q\le3$ and fixed smooth forms
$\omega\in\Omega^q(M)$ and $\varphi\in\Omega^{3-q}(M)$, we have
\begin{equation}\label{eq:weak-lie-convergence}
 \left|\int_M\varphi\wedge\bigl(W_qL^qR_q\omega-\mathcal L_X\omega\bigr)\right|
 \le Ch.
\end{equation}
With $C=C_{\varphi ,\omega}$ independent of $h$.
\end{theorem}

\begin{proof}
We first estimate the error in contracting a smooth $k$-form $\rho$.
The local Whitney interpolant satisfies
$\|W_kR_k\rho-\rho\|_{L^\infty(T)}\le C_\rho h_T$~\cite{Wilson2007}.
In local coordinates, each coefficient $a$ of $\Phi_T^*\rho$ is a signed
sum of products of the form
\[
 (\rho_I\circ\Phi_T)\prod_{r=1}^k\partial_{j_r}\Phi_T^{i_r},
 \qquad I=(i_1,\ldots,i_k),
\]
where $\rho_I$ is a coordinate coefficient of $\rho$. Each first derivative
of $\Phi_T$ is $O(h_T)$ by \eqref{eq:scalar-map-bounds}, so $a=O(h_T^k)$.
When we differentiate the product in a reference coordinate $z_\ell$,
the derivative of the first factor is
\[
 \partial_\ell(\rho_I\circ\Phi_T)
 =\sum_m \bigl((\partial_m\rho_I)\circ\Phi_T\bigr)\,\partial_\ell\Phi_T^m
 =O(h_T).
\]
The $k$ remaining factors then give a term of order $h_T^{k+1}$.
If the derivative falls on one of those factors instead, it replaces
$\partial_{j_r}\Phi_T^{i_r}$ by $\partial_\ell\partial_{j_r}\Phi_T^{i_r}=O(h_T^2)$.
The other $k-1$ first derivatives contribute $O(h_T^{k-1})$, so this term
is also $O(h_T^{k+1})$. Thus $\|Da\|_\infty\le C_\rho h_T^{k+1}$.
Since $\widehat T$ is convex and has fixed diameter, for any $z,z'\in\widehat T$,
\[
 |a(z)-a(z')|\le |z-z'|\,\|Da\|_\infty\le C_\rho h_T^{k+1}.
\]
Write $\widehat\rho=\Phi_T^*\rho=\sum_J a_J(z)\,dz^J$ and let $\widehat b$
be the barycenter of $\widehat T$. Define the constant form
$\rho_0=\sum_J a_J(\widehat b)\,dz^J$. Taking $z'=\widehat b$ in the
preceding estimate gives $\|\widehat\rho-\rho_0\|_\infty\le C_\rho h_T^{k+1}$.

Let $\widehat I_k$ denote Whitney interpolation on $\widehat T$. For any
smooth $k$-form $\eta$ on $\widehat T$,
\[
 \widehat I_k\eta=\sum_{F\in\widehat T^{(k)}}
    \left(\int_F\eta\right)W_k\chi_F.
\]
Each face integral is bounded by $|F|\,\|\eta\|_\infty$ . The faces and their Whitney basis forms are fixed,
so $\|\widehat I_k\eta\|_\infty\le C_k\|\eta\|_\infty$.
 Since $\widehat I_k\rho_0=\rho_0$ , it gives
\[
\begin{aligned}
 \|\widehat I_k\widehat\rho-\widehat\rho\|_\infty
 &=\|\widehat I_k(\widehat\rho-\rho_0)-(\widehat\rho-\rho_0)\|_\infty\\
 &\le (C_k+1)\|\widehat\rho-\rho_0\|_\infty
 \le C_\rho h_T^{k+1}.
\end{aligned}
\]
Since $\int_{\Phi_T(F)}\rho=\int_F\widehat\rho$, the interpolation error
on $T$ is $(\Phi_T^{-1})^*(\widehat I_k\widehat\rho-\widehat\rho)$.
This pullback applies $(D\Phi_T)^{-1}$ to each of the $k$ vector arguments.
By \eqref{eq:scalar-map-bounds},
\[
 \|W_kR_k\rho-\rho\|_{L^\infty(T)}
 \le C\|(D\Phi_T)^{-1}\|_\infty^k
       \|\widehat I_k\widehat\rho-\widehat\rho\|_\infty
 \le C_\rho h_T^{-k}h_T^{k+1}=C_\rho h_T.
\]

For $1\le k\le3$, put
$r_h^k(\rho)=\iota^k R_k\rho-R_{k-1}(\iota_X\rho)$.
On each tetrahedron $T$, write $X_T(\Phi_T(z))=D\Phi_T(z)\nu_T$.
Since $X(\Phi_T(z))=D\Phi_T(z)Y(z)$, the bounds on $D\Phi_T$ and
$\nu_T-Y$ in \eqref{eq:scalar-map-bounds} and
\eqref{eq:scalar-reconstruction-bounds} give
$\|X_T-X\|_{L^\infty(T)}\le Ch_T$. In particular, $X_T$ is uniformly bounded.
The error in the contracted form on $T$ is
\[
\begin{aligned}
 \eta_T&:=\iota_{X_T}(W_kR_k\rho)-\iota_X\rho\\
       &=\iota_{X_T}(W_kR_k\rho-\rho)+\iota_{X_T-X}\rho.
\end{aligned}
\]
The first term is bounded by $\|X_T\|_{L^\infty(T)}\|W_kR_k\rho-\rho\|_{L^\infty(T)}$
and the second by $\|X_T-X\|_{L^\infty(T)}\|\rho\|_{L^\infty(T)}$.
Thus $\|\eta_T\|_{L^\infty(T)}\le C_\rho h_T$.

Fix an oriented $(k-1)$-simplex $S$. Changing variables by $\Phi_T$ in
\eqref{eq:contraction-trace} gives the integral of $\iota_{X_T}(W_kR_k\rho)$
over $S$, restricted from $T$. Since the weights $\theta_{S,T}$ are
nonnegative and sum to one, we can subtract
$(R_{k-1}(\iota_X\rho))(S)=\int_S\iota_X\rho$ inside the average. Hence
\[
\begin{aligned}
 r_h^k(\rho)(S)&=\sum_{T\supset S}\theta_{S,T}\int_S\eta_T,\\
 |r_h^k(\rho)(S)|&\le |S|\sum_{T\supset S}\theta_{S,T}\|\eta_T\|_{L^\infty(T)}
 \le C_\rho |S|\sum_{T\supset S}\theta_{S,T}h_T\le C_\rho h|S|.
\end{aligned}
\]
Here $|S|$ is the $(k-1)$-dimensional measure of $S$. 

Now fix a tetrahedron $U$ and expand the reconstruction in its Whitney basis.
For each $(k-1)$-face $S$ of $U$, the bound on $D\Phi_U$ gives
$|S|\le Ch_U^{k-1}$. Pulling the reference basis form back by $\Phi_U^{-1}$
gives $\|W_{k-1}\chi_S\|_{L^\infty(U)}\le Ch_U^{1-k}$. The two powers of $h_U$ therefore cancel in the basis expansion, giving
\[
\begin{aligned}
 \|W_{k-1}r_h^k(\rho)\|_{L^\infty(U)}
 &\le\sum_{S\in U^{(k-1)}}|r_h^k(\rho)(S)|\,\|W_{k-1}\chi_S\|_{L^\infty(U)}\\
 &\le C_\rho h\sum_{S\in U^{(k-1)}} h_U^{k-1} h_U^{1-k}\le C_\rho h,
\end{aligned}
\]
since $U$ has only $\binom{4}{k}$ such faces. Taking the maximum over $U$
gives the same bound on $M$. Finally, adding and subtracting the Whitney
interpolant of $\iota_X\rho$ gives
\[
\begin{aligned}
 E_h^k(\rho)&:=W_{k-1}\iota^kR_k\rho-\iota_X\rho\\
 &= W_{k-1}r_h^k(\rho)+\bigl(W_{k-1}R_{k-1}(\iota_X\rho)-\iota_X\rho\bigr).
\end{aligned}
\]
The interpolation estimate applied to $\iota_X\rho$ bounds the second term,
so $\|E_h^k(\rho)\|_{L^\infty(M)}\le C_\rho h$. We take terms outside degrees $0$ through $3$ to be zero.
By Cartan's formula and the commuting identities \eqref{eq:whitney-identities},
\[
 W_q L^q R_q\omega-\mathcal L_X\omega
   =dE_h^q(\omega)+E_h^{q+1}(d\omega).
\]
For $q\ge1$, apply Stokes' theorem on each tetrahedron $T$ to the $2$-form
$\varphi\wedge E_h^q(\omega)$. Since $\deg\varphi=3-q$, the Leibniz rule gives
\[
 \int_{\partial T}\varphi\wedge E_h^q(\omega)
 =\int_T d\varphi\wedge E_h^q(\omega)
  +(-1)^{3-q}\int_T\varphi\wedge dE_h^q(\omega).
\]
On a face shared by two tetrahedra, the Whitney form
$W_{q-1}\iota^qR_q\omega$ has the same tangential trace from both sides.
The forms $\iota_X\omega$ and $\varphi$ are smooth, so the restrictions
of $\varphi\wedge E_h^q(\omega)$ to that face agree as well.
The two tetrahedra induce opposite orientations on the face, hence its
two integrals cancel. Since $M$ is closed, every face is paired in this way.
Summing over $T$ and substituting into the preceding Cartan identity gives
\[
\begin{aligned}
 \int_M\varphi\wedge\bigl(W_qL^qR_q\omega-\mathcal L_X\omega\bigr)
 &=(-1)^q\int_M d\varphi\wedge E_h^q(\omega)\\
 &\quad+\int_M\varphi\wedge E_h^{q+1}(d\omega).
\end{aligned}
\]
For $q=0$ the first term vanishes, and for $q=3$ the second vanishes.
By H\"older's inequality, the absolute value of the right-hand side is bounded,
up to a constant depending only on $M$, by
\[
\begin{aligned}
 &\|d\varphi\|_{L^1(M)}\,\|E_h^q(\omega)\|_{L^\infty(M)}
 +\|\varphi\|_{L^1(M)}\,\|E_h^{q+1}(d\omega)\|_{L^\infty(M)}\\
 &\qquad\le h\bigl(C_\omega\|d\varphi\|_{L^1(M)}
                  +C_{d\omega}\|\varphi\|_{L^1(M)}\bigr).
\end{aligned}
\]
 proving
\eqref{eq:weak-lie-convergence}.
\end{proof}

\begin{remark}
Under the hypotheses of Theorem~\ref{thm:weak-lie-convergence}, the dual
Lie derivative also converges weakly to $\mathcal L_X$ in every degree.
Here we discretize $\omega$ by integration over the oriented dual $q$-cells
and test its discrete Lie derivative against $R_{3-q}\varphi$ using the
pairing $\kappa_q^{-1}P_q$. The proof is similar, using the primal--dual
identity \eqref{eq:dual-transport-pairing}, and is left to the reader.
\end{remark}

To summarize, we have shown that both cup helicity and Whitney helicity
converge to smooth helicity at first order on the shape-regular refinements
of Definition~\ref{def:controlled-refinement}. This suggests that either
discrete helicity can be used to approximate smooth helicity in computations.
We also studied the discrete Lie derivative as a numerical approximation to
$\mathcal L_X$, proving convergence in $L^\infty(M)$ for the primal scalar
operator and weak convergence on forms of every degree.

\ifshowalgorithms
\section{Linking algorithms}\label{sec:algorithms}

The linking identity gives us a way to compute helicity directly from the
exact $2$-cochain $b=\delta a$, without recovering a potential. We focus
on Whitney helicity $\Hw(a)$ in this section. To compute cup helicity
$\Hc(a)$, we apply the same algorithms to the primal and dual 1-cycles
in \eqref{eq:cup-chains}.

To simplify the computation of linking numbers, we take $M=S^3$ throughout this section.

\subsection{Deterministic crossing scan}\label{alg:sec:03-deterministic}

Choose a point off the two 1-skeletons and identify its complement in
$S^3$ with $\R^3$ by an orientation-preserving stereographic projection chart.
The chart places both 1-cycles in $\R^3$, but their edges may still be curved. To count crossings between straight segments, we deform the curved edges into polygonal arcs. The union of the primal and dual 1-skeletons is a finite embedded graph whose piecewise smooth edges have finite total curvature. By Denne and Sullivan~\cite[\S\S~4--5]{DenneSullivan2008}, a single ambient isotopy of $\R^3$ carries this graph to a polygonal embedding. The same family of homeomorphisms, starting from the identity, thus moves both cycles together. Their supports remain disjoint throughout the deformation, and their linking number is preserved. 

Figure~\ref{fig:polygonal-linked-cycles} shows one oriented primal edge $[u,v]$ of coefficient $w$. Its polygonal representative runs from $u$ to $v$ along three successive straight segments $s_1,s_2,s_3$. Each segment carries $w$, so this edge contributes $w(s_1+s_2+s_3)$ to the primal 1-cycle. We treat dual edges in the same way. Write $c_i$ for the coefficient of primal segment $i$ and $d_j$ for that of dual segment $j$, with $1\le i\le n^{\mathrm{pr}}$ and $1\le j\le n^{\mathrm{du}}$. The segment crossings below are used to compute the linking number of the two complete 1-cycles.

\begin{figure}[H]
\centering
\definecolor{primalred}{RGB}{184,49,58}
\begin{tikzpicture}[line cap=round,line join=round,
  edge/.style={primalred,line width=1pt},
  arrow/.style={-{Stealth[length=2mm,width=1.5mm]}},
  every node/.style={font=\small}]
  \begin{scope}
    \coordinate (U) at (0,0);
    \coordinate (V) at (3.4,-.05);
    \draw[edge,arrow] (U) .. controls (.54,1.12) and (1.13,1.19) .. (1.73,1.13);
    \draw[edge] (1.73,1.13) .. controls (2.33,1.07) and (2.82,.65) .. (V);
    \fill (U) circle (1.4pt) node[below] {$u$};
    \fill (V) circle (1.4pt) node[below] {$v$};
    \node at (1.65,.50) {$w[u,v]$};
  \end{scope}
  \draw[black!65,arrow,line width=.6pt] (3.8,.51)--(4.55,.51);
  \begin{scope}[xshift=5cm]
    \coordinate (U) at (0,0);
    \coordinate (P) at (.92,.99);
    \coordinate (Q) at (2.20,1.02);
    \coordinate (V) at (3.4,-.05);
    \foreach \from/\dest/\idx in {U/P/1,P/Q/2,Q/V/3}{
      \draw[edge] (\from)--(\dest) coordinate[pos=.43] (a) coordinate[pos=.64] (b)
        node[pos=.50,below=4pt] {$w s_{\idx}$};
      \draw[edge,arrow] (a)--(b);
    }
    \fill (U) circle (1.4pt) node[below] {$u$};
    \fill (P) circle (1.2pt) node[above] {$z_1$};
    \fill (Q) circle (1.2pt) node[above] {$z_2$};
    \fill (V) circle (1.4pt) node[below] {$v$};
  \end{scope}
\end{tikzpicture}
\caption{A single oriented edge and its polygonal representative. The interior vertices $z_1,z_2$ divide the arc into three segments, each with coefficient $w$.}
\label{fig:polygonal-linked-cycles}
\end{figure}
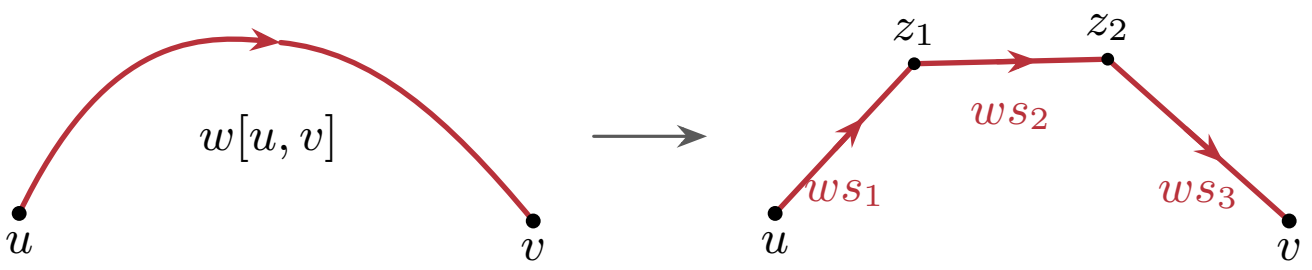

Project both polygonal cycles orthogonally onto the plane $e^\perp$ and
view the diagram from the $+e$ side, where $e$ is a unit vector. A
\emph{mixed crossing} is a crossing of a projected primal segment $i$
with a projected dual segment $j$. Choose $e$ so that the two segments
cross in their interiors at a nonzero angle, no third segment projects
through the same point, and no vertex projects onto a nonincident
segment. Let $\Cross$ be the set of pairs $(i,j)$ that give mixed
crossings.

At each mixed crossing, the segment closer to the viewer passes over
the other. Rotate the diagram until the arrow on the overpassing segment
points to the right. Following the arrow on the underpassing segment,
read the sign as follows (Figure~\ref{fig:mixed-crossing-signs}).
\[
 \sigma_{ij}=
 \begin{cases}
  +1 & \text{if the underpassing arrow crosses upward on the page},\\
  -1 & \text{if the underpassing arrow crosses downward on the page}.
 \end{cases}
\]

\begin{figure}[H]
\centering
\begin{tikzpicture}[line width=0.9pt,line cap=round,
  directed/.style={-{Stealth[length=2.6mm,width=1.9mm]}}]
  % Positive crossing: the underpassing strand is directed upwards.
  \draw (-0.54,-0.68) -- (-0.10,-0.126);
  \draw[directed] (0.10,0.126) -- (0.54,0.68);
  \draw[directed] (-0.96,0) -- (0.96,0);
  \node[anchor=west] at (1.12,0) {$\sigma_{ij}=+1$};
  \begin{scope}[xshift=4.7cm]
    % Negative crossing: reverse only the underpassing strand.
    \draw (0.54,0.68) -- (0.10,0.126);
    \draw[directed] (-0.10,-0.126) -- (-0.54,-0.68);
    \draw[directed] (-0.96,0) -- (0.96,0);
    \node[anchor=west] at (1.12,0) {$\sigma_{ij}=-1$};
  \end{scope}
\end{tikzpicture}
\caption{Signs of oriented mixed crossings viewed from the $+e$ side.}
\label{fig:mixed-crossing-signs}
\end{figure}
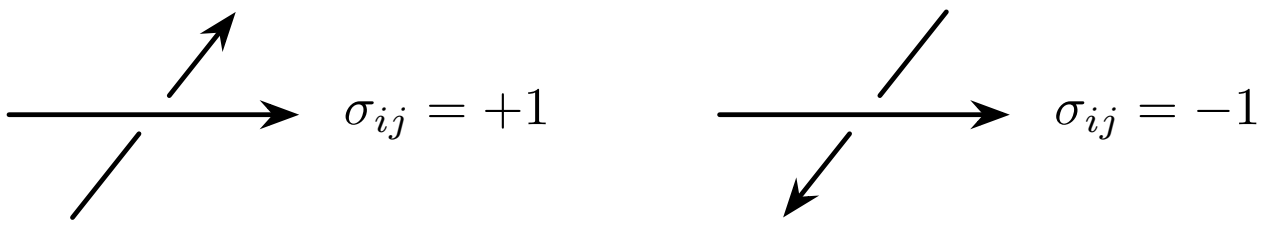

\begin{proposition}\label{alg:prop:crossing}
Under these hypotheses,
\begin{equation}\label{alg:eq:crossing}
 \Hw(a)=\Lk(c^{\mathrm{pr}},c^{\mathrm{du}})
       =\frac12\sum_{(i,j)\in\Cross}\sigma_{ij}c_i d_j.
\end{equation}
\end{proposition}
\begin{proof}
The first equality follows from Theorem~\ref{thm:whitney-linking} on $S^3$.  The formula for linking in terms of signed mixed crossings
\cite{qu2021} extends bilinearly to real cycles.
\end{proof}

Algorithm~\ref{alg:deterministic-scan} uses the same signed crossing list for
every flux on the fixed geometry.

\begin{algorithm}[H]
\caption{Deterministic crossing scan }
\label{alg:deterministic-scan}
\small
\begin{algorithmic}[1]
\Statex \textbf{Input} \enspace $K$, the oriented primal and dual segments with their parent edges,
\Statex \hspace{2.9em} a generic projection direction $e$, and the exact face flux $b=\delta a$.
\State Project all segments onto $e^\perp$ and set $\Cross\gets\varnothing$.
\For{each primal segment $i$ and each dual segment $j$}
  \If{their projected interiors cross}
    \State Add $(i,j)$ to $\Cross$ and store its sign $\sigma_{ij}$ as in Figure~\ref{fig:mixed-crossing-signs}.
  \EndIf
\EndFor
\Statex \textit{For each input flux $b$}
\Statex Orient each $T$ by $M$ and set $B_{T,w}=b(f_w)$ on the outward face opposite $w$.
\For{each oriented primal edge $[u,v]$}
  \State Assign $\displaystyle\frac{1}{12}\sum_{T\supset[u,v]}(B_{T,u}-B_{T,v})$ to $[u,v]$.
\EndFor
\For{each oriented dual edge $f^*$}
  \State Assign $b(f)$ to $f^*$, using the corresponding orientation of $f$.
\EndFor
\State $H\gets 0$.
\For{each $(i,j)\in\Cross$}
  \State Read $c_i$ and $d_j$ from the parent edges of $i$ and $j$, respectively.
  \State $H\gets H+\sigma_{ij}c_i d_j$.
\EndFor
\State \Return $H/2$.
\end{algorithmic}
\end{algorithm}

The crossing list $\Cross$ and its signs $\sigma_{ij}$ depend only on the polygonal realizations of the two 1-skeletons in the chosen chart of $S^3$ and the projection direction $e$. We therefore only need to compute and store them once.

Let $N=\sum_{q=0}^3\#K^{(q)}$ be the number of simplices in the
triangulation of $S^3$.  We count each floating-point addition, subtraction, multiplication, or division as one FLOP. Reading a stored crossing $(i,j)$ and its sign $\sigma_{ij}$ takes $O(1)$ time.
Assembling the edge coefficients costs $O(N)$ FLOPs. At each stored crossing in \eqref{alg:eq:crossing}, we multiply $c_i$ by $d_j$ and add or subtract the product according to $\sigma_{ij}$. The sum therefore costs $2|\Cross|$ FLOPs, followed by one FLOP to return $H/2$.
The vertices and edges of the two 1-skeletons together account for all $N$ simplices of $K$. Since neither graph has isolated vertices and every edge contributes at least one segment, $N=O(n^{\mathrm{pr}}+n^{\mathrm{du}})$. For the supplied regular projection, the scan tests exactly $n^{\mathrm{pr}}n^{\mathrm{du}}$ pairs. Each test uses a constant number of FLOPs and exact unit-cost geometric predicates, so preprocessing takes $O(n^{\mathrm{pr}}n^{\mathrm{du}})$ FLOPs. Each crossing contributes one entry to each of the two lists of crossings, giving a total storage bound of $O(n^{\mathrm{pr}}+n^{\mathrm{du}}+|\Cross|)$ words.

\subsection{Endpoint sampling with a fixed budget}\label{alg:sec:04-randomized}

Even with the crossing list in hand, Algorithm~\ref{alg:deterministic-scan} evaluates all $|\Cross|$ terms in \eqref{alg:eq:crossing} for each new flux $b$. In the worst case, a fixed fraction of all primal--dual segment pairs cross, so $|\Cross|$ is of order $n^{\mathrm{pr}}n^{\mathrm{du}}$. If the two segment counts are of the same order, this is quadratic in the total number of segments. When an approximation suffices, we can estimate $\Hw(a)$ from a fixed number of sampled crossings. We weight their contributions so that the estimate is unbiased.

Write $c=(c_i)_{i=1}^{n^{\mathrm{pr}}}$ and $d=(d_j)_{j=1}^{n^{\mathrm{du}}}$ for the signed flux coefficients of the oriented primal and dual segments, respectively. Each segment carries the coefficient of its parent edge. The infinity norms $\|c\|_\infty=\max_i|c_i|$ and $\|d\|_\infty=\max_j|d_j|$ give the largest flux magnitude on each side.

Let $\rho_i$ count the dual segments that cross primal segment $i$ in the projection, and let $\eta_j$ count the primal segments that cross dual segment $j$. Segment $i$ occurs in $\rho_i$ terms of the crossing sum, each proportional to $c_i$. We therefore give it sampling weight $|c_i|\rho_i$ when starting on the primal side. At equal flux magnitude, a segment with twice as many crossings is twice as likely to be chosen. At equal crossing count, selection is proportional to the flux magnitude. The dual rule uses the weight $|d_j|\eta_j$. The total weights $A_c,A_d$ and the bounds $B_c,B_d$ for a single draw are
\[
 \begin{aligned}
 A_c&=\sum_i|c_i|\rho_i,& B_c&=\tfrac12 A_c\|d\|_\infty,\\
 A_d&=\sum_j|d_j|\eta_j,& B_d&=\tfrac12 A_d\|c\|_\infty,
 \end{aligned}
 \qquad B=\min(B_c,B_d).
\]
For $m\ge1$ independent draws with replacement, let $Z_r$ denote the weighted contribution of draw $r$ in Algorithm~\ref{alg:endpoint-sampling} and let $\widehat H_m=m^{-1}\sum_{r=1}^m Z_r$. Starting on the primal side gives $|Z_r|\le B_c$, and starting on the dual side gives $|Z_r|\le B_d$. For a fixed sampling budget $m$, choosing the smaller of $B_c$ and $B_d$ gives the tighter guarantee in \eqref{alg:eq:tail}.

\begin{algorithm}[H]
\caption{Endpoint sampling with a fixed budget for Whitney helicity}
\label{alg:endpoint-sampling}
\small
\begin{algorithmic}[1]
\Statex \textbf{Input} \enspace Exact face flux $b$, the oriented segments and their parent edges,
\Statex \hspace{2.9em} the signed crossing lists for both skeletons, and an integer $m\ge 1$.
\State Assemble the segment coefficients $c_i,d_j$ from $b$ as in Algorithm~\ref{alg:deterministic-scan}.
\State Read the list lengths $\rho_i$ for primal segments and $\eta_j$ for dual segments.
\State $A_c\gets\sum_i |c_i|\rho_i$, \quad $A_d\gets\sum_j |d_j|\eta_j$.
\State $B_c\gets A_c\|d\|_\infty/2$, \quad $B_d\gets A_d\|c\|_\infty/2$.
\If{$\min(B_c,B_d)=0$}
  \State \Return $0$.
\EndIf
\State $S\gets 0$.
\If{$B_c\le B_d$}
  \State Set $p_i\gets |c_i|\rho_i/A_c$ for each primal segment $i$.
  \For{$r=1,\ldots,m$}
    \State Draw a primal segment $I$ with probabilities $(p_i)_i$.
    \State Draw a dual segment $J$ uniformly from $\{j:(I,j)\in\Cross\}$.
    \State $S\gets S+(A_c/2)\sgn(c_I)\,\sigma_{IJ}\,d_J$.
  \EndFor
\Else
  \State Set $p_j\gets |d_j|\eta_j/A_d$ for each dual segment $j$.
  \For{$r=1,\ldots,m$}
    \State Draw a dual segment $J$ with probabilities $(p_j)_j$.
    \State Draw a primal segment $I$ uniformly from $\{i:(i,J)\in\Cross\}$.
    \State $S\gets S+(A_d/2)\sgn(d_J)\,\sigma_{IJ}\,c_I$.
  \EndFor
\EndIf
\State \Return $S/m$.
\end{algorithmic}
\end{algorithm}

We now bound the probability that the absolute error of the sampling estimator $\widehat H_m$ exceeds a given tolerance.

\begin{proposition}\label{alg:prop:sampling}
If $B=0$, then $c_i d_j=0$ at every crossing $(i,j)\in\Cross$, so the
algorithm returns $\Hw(a)=0$. If $B>0$ and the input data are exact,
then for every $\varepsilon>0$,
\begin{equation}\label{alg:eq:tail}
 \E\widehat H_m=\Hw(a),\qquad
 \Prob\{|\widehat H_m-\Hw(a)|\ge\varepsilon\}
 \le2\exp\!\left(-\frac{m\varepsilon^2}{2B^2}\right).
\end{equation}
\end{proposition}
\begin{proof}
If $B_c=0$, either $\|d\|_\infty=0$ or every $c_i$ with a crossing
neighbor is zero, so all crossing products vanish.  The case $B_d=0$
is symmetric. Suppose $B>0$ and $B_c\le B_d$, so Algorithm~\ref{alg:endpoint-sampling} starts on the primal side. On draw $r$, it first selects segment $i$ with probability $|c_i|\rho_i/A_c$, then chooses uniformly among the $\rho_i$ dual segments crossing $i$. Thus, for $(i,j)\in\Cross$ with $c_i\ne0$,
\[
 \Prob\{I=i,J=j\}
 =\Prob\{I=i\}\,\Prob\{J=j\mid I=i\}
 =\frac{|c_i|\rho_i}{A_c}\,\frac{1}{\rho_i}
 =\frac{|c_i|}{A_c}.
\]
The same joint probability formula holds when $c_i=0$, since segment $i$ is then selected with probability zero. Averaging the contribution of draw $r$ over all crossings gives
\[
\begin{aligned}
 \E Z_r
 &=\sum_{(i,j)\in\Cross}\Prob\{I=i,J=j\}\,
     \frac{A_c}{2}\sgn(c_i)\,\sigma_{ij}d_j\\
 &=\sum_{(i,j)\in\Cross}\frac{|c_i|}{A_c}\,
     \frac{A_c}{2}\sgn(c_i)\,\sigma_{ij}d_j\\
 &=\frac12\sum_{(i,j)\in\Cross}\sigma_{ij}c_i d_j
 =\Hw(a)
\end{aligned}
\]
by \eqref{alg:eq:crossing}. Here $|Z_r|\le A_c\|d\|_\infty/2=B_c=B$. Starting on the dual side gives the same expectation and the bound $|Z_r|\le B_d=B$ by interchanging the roles of the two segments. Thus $\E\widehat H_m=m^{-1}\sum_{r=1}^m\E Z_r=\Hw(a)$. The draws are independent.
Hoeffding's inequality~\cite{hoeffding1963} therefore yields
\eqref{alg:eq:tail}.
\end{proof}

This provides us a guideline: to obtain absolute error at most $\varepsilon$ with probability at least
$1-\delta$, choose $m\ge\lceil 2B^2\varepsilon^{-2}\log(2/\delta)\rceil$
before sampling. When the signed terms in \eqref{alg:eq:crossing} nearly
cancel, even a small absolute error can be large relative to $|\Hw(a)|$.

Once the crossing lists have been stored, Algorithm~\ref{alg:endpoint-sampling} assembles the edge coefficients for each new flux $b$ as in \S\ref{alg:sec:03-deterministic}. It then makes one pass over the primal and dual segments to compute the sampling weights and the two bounds $B_c,B_d$, reading each crossing count from the length of the corresponding list. For $B>0$, Vose's alias method~\cite{vose1991} builds a sampling table for the chosen side in time linear in its number of segments. Together these steps cost $O(n^{\mathrm{pr}}+n^{\mathrm{du}})$ FLOPs and use $O(n^{\mathrm{pr}}+n^{\mathrm{du}})$ additional words of memory. The timings in \S\ref{alg:sec:05-experiments} include this setup for each flux.

Each draw selects a segment from the alias table and a crossing from its list in constant time. With $A_c/2$ or $A_d/2$ computed once, updating the sum takes one multiplication and one addition or subtraction per draw. The $m$ updates and the final division by $m$ therefore cost $2m+1$ FLOPs. Including setup and sampling, the total for one input flux with stored crossing lists is
\[
 O\!\left(n^{\mathrm{pr}}+n^{\mathrm{du}}+m\right)\ \text{FLOPs}.
\]

Sampling offers its main advantage on a warm start, when the crossing lists are reused for a new flux and $m\ll|\Cross|$ draws suffice for the desired accuracy.

\subsection{Numerical results}\label{alg:sec:05-experiments}

Our numerical experiments so far are limited to small examples for which the helicity can be checked independently. Whether the algorithms are practical on large 1-skeletons remains to be tested.

We test the linking algorithms on the unit cube $\Omega=[0,1]^3$ in the stereographic chart of $S^3$. We choose the potential
\[
 \alpha=\phi(-y\,dx+x\,dy+dz),\qquad
 \phi=x(1-x)y(1-y)z(1-z).
\]
It vanishes on the boundary, and $\alpha\wedge d\alpha=2\phi^2\,dx\wedge dy\wedge dz$ gives the exact helicity $1/13500$. On each mesh we integrate $\alpha$ along the edges to obtain $a_h=R_1\alpha$ and set $b_h=\delta a_h$.
On each boundary face the tangential traces of $\alpha$ and $d\alpha$ vanish, so
$a_h|_{\partial\Omega}=0$ and $b_h|_{\partial\Omega}=0$. On any triangulation of
$S^3$ extending the cube mesh, the coboundary of the zero extension of $a_h$ is
the zero extension of $b_h$. Thus no dual edge crossing $\partial\Omega$ has nonzero
flux, and exterior tetrahedra contribute zero to the primal flux on boundary edges.
Both $\cpr$ and $\cdu$ are therefore supported in $\Omega$, and the crossing formula
\eqref{alg:eq:crossing} applies. Both algorithms take $b_h$ as input and compute
crossings in the projection along $(1,\sqrt{2},\sqrt{3})$.

For the refinement test we divide $\Omega$ into $n^3$ equal cubes and split each cube into six Kuhn tetrahedra. We use $n=6,12,16,20,32$. At the four coarser resolutions we also test a perturbed mesh by displacing each interior vertex independently by a vector drawn uniformly from $[-0.1/n,0.1/n]^3$. Table~\ref{alg:fig:helicity-refinement} gives the relative error in Whitney helicity on both mesh families.

\begin{table}[H]
\centering\small
\caption{Relative error of Whitney helicity against the exact value $1/13500$.}
\label{alg:fig:helicity-refinement}
\begin{tabular}{@{}r rr rr@{}}
\toprule
& \multicolumn{2}{c}{Regular mesh} & \multicolumn{2}{c}{Perturbed mesh}\\
\cmidrule(lr){2-3}\cmidrule(l){4-5}
$n$ & $h$ & Error (\%) & $h$ & Error (\%)\\
\midrule
6  & 0.28868 & 14.4761 & 0.32419 & 14.6327\\
12 & 0.14434 & 3.7967  & 0.16574 & 3.8571\\
16 & 0.10825 & 2.1507  & 0.12431 & 2.1892\\
20 & 0.08660 & 1.3809  & 0.09979 & 1.4053\\
32 & 0.05413 & 0.5413  & ---     & ---\\
\bottomrule
\end{tabular}
\end{table}

To measure sampling error on a fixed mesh, we take $n=20$ and rescale the smooth potential to give $\Hw(a_h)=1$. We obtain the discrete reference value by integrating $W_1a_h\wedge W_2b_h$ over each tetrahedron and summing. For each sampling budget we report the root mean square error over 60 independent runs.

\begin{table}[H]
\centering\small
\caption{Absolute sampling error for the smooth field at $n=20$, with $\Hw=1$.}
\begin{tabular}{@{}rrr@{}}
\toprule
Draws $m$ & RMSE over 60 runs & $\sqrt{\operatorname{Var}(Z_r)/m}$\\
\midrule
256       & 0.39721 & 0.41328\\
1\,024    & 0.22726 & 0.20664\\
4\,096    & 0.10476 & 0.10332\\
16\,384   & 0.05952 & 0.05166\\
\bottomrule
\end{tabular}
\end{table}

We compare the crossing scan and endpoint sampling with a third method that recovers a potential from $b_h$ and evaluates its Whitney helicity. Table~\ref{alg:tab:experiments-timing} separates the one-time setup from the cost of evaluating a single flux, with $16\,384$ draws per sampling estimate. We ran these tests on an Apple M5 Pro with one BLAS thread. To check the two deterministic methods against the elementwise Whitney integral, we also use a Gaussian potential and a potential whose signed crossing contributions nearly cancel.

\begin{table}[H]
\centering\small
\caption{Setup and evaluation times on regular meshes, with $16\,384$ draws per sampling estimate.}
\label{alg:tab:experiments-timing}
\begin{tabular}{@{}rrrrrrrr@{}}
\toprule
& & & \multicolumn{2}{c}{Setup (s)} & \multicolumn{3}{c}{Time per flux (ms)}\\
\cmidrule(lr){4-5}\cmidrule(l){6-8}
$n$ & Tetrahedra & Crossings & Scan & Primitive & Scan & Primitive & Sampling\\
\midrule
6  & 1\,296   & 14\,582        & 0.076 & 0.008 & 0.025 & 0.115 & 1.94\\
12 & 10\,368  & 229\,780       & 0.772 & 0.815 & 0.172 & 3.49  & 13.90\\
16 & 24\,576  & 722\,162       & 2.19  & 4.54  & 0.460 & 13.38 & 32.65\\
20 & 48\,000  & 1\,756\,720    & 4.76  & 27.90 & 1.023 & 44.74 & 65.98\\
32 & 196\,608 & 11\,469\,528   & 27.57 & ---   & 6.539 & ---   & 274.78\\
\bottomrule
\end{tabular}
\end{table}

\begin{table}[H]
\centering\small
\caption{Maximum absolute difference between helicity evaluations.}
\begin{tabular}{@{}llr@{}}
\toprule
Test & Methods compared & Maximum difference\\
\midrule
All three fields, $n=6,12,16,20$ & Scan, primitive, Whitney integral & $1.51\times10^{-14}$\\
Smooth field, $n=32$ & Scan and Whitney integral & $3.33\times10^{-16}$\\
\bottomrule
\end{tabular}
\end{table}

To test repeated evaluations on a fixed mesh, we scale the smooth potential to unit discrete helicity and multiply it by $1+u(x-1/2)+v(y-1/2)$, drawing $u$ and $v$ independently and uniformly from $[-0.6,0.6]$ for each new flux. Both algorithms reuse the crossing lists. The sampling times in Table~\ref{alg:tab:sampling-crossover} include rebuilding the weights for each flux and drawing $65\,536$ crossings.

\begin{table}[H]
\centering\small
\caption{Median time per flux with stored crossing lists and $65\,536$ draws. At $n=56$, the median speedup over paired evaluations was $2.38\times$.}
\label{alg:tab:sampling-crossover}
\begin{tabular}{@{}rrrrrr@{}}
\toprule
$n$ & Tetrahedra & Fluxes & Scan (ms) & Sampling (ms) & Relative RMSE (\%)\\
\midrule
12 & 10\,368       & 90     & 0.147  & 1.777  & 2.68\\
20 & 48\,000       & 90     & 0.829  & 2.703  & 2.67\\
32 & 196\,608      & 90     & 4.890  & 5.411  & 2.22\\
56 & 1\,053\,696   & 4\,919 & 45.044 & 18.873 & 2.55\\
\bottomrule
\end{tabular}
\end{table}

For the largest test, we tile $\Omega$ with $30^3$ congruent blocks and use the $n=56$ mesh in each block. We generate a potential independently for each block as above, scale it by $30^{-3/2}$, and map it affinely onto the block. The potentials vanish on the block boundaries and together define a single field on $\Omega$. We use $65\,536$ draws per block. The total times in Table~\ref{alg:tab:tiled-cube} include flux generation and construction of the crossing lists.

\begin{table}[H]
\centering\small
\caption{Whitney helicity on $27\,000$ blocks. The independent reference is $1.008547495948991$.}
\label{alg:tab:tiled-cube}
\begin{tabular}{@{}lrr@{}}
\toprule
 & Deterministic scan & Endpoint sampling\\
\midrule
Computed helicity            & $1.00854750$      & $1.00874061$\\
Relative error (\%)          & $6.6\times10^{-14}$ & $0.0191$\\
Total time (min)             & 21.67             & 9.21\\
Speedup relative to scan     & $1.00\times$      & $\mathbf{2.35}\times$\\
Blockwise relative RMSE (\%) & ---               & 2.60\\
\bottomrule
\end{tabular}
\end{table}

\fi

\par
\section*{Acknowledgements}
GPT-6 Astra was used  to propose proofs for the results on Markov trajectories and convergence part. The other central mathematical ideas and the proof strategies throughout the paper are owned by the author. I also used the model for translation and to speed up the preparation of the final manuscript. The author takes full responsibility for the correctness of all results in this paper.

I am grateful to Professor Kaibo Hu of the University of Oxford for introducing me to this problem and for the many suggestions and observations he shared as the work developed. I thank Puchun Zhou for helping me shape the direction of the project and for advice on the writing. Several online discussions with Anil Hirani deepened my understanding of the problem from the perspective of discrete exterior calculus.

\begingroup
\small
\raggedright
\bibliographystyle{abbrv}
\bibliography{references}
\endgroup
\par\bigskip
\begin{samepage}
\footnotesize
\noindent\textsc{Jingqi Li, Fudan University}\par\nopagebreak
\noindent\textit{Email address:}
\href{mailto:23300180152@m.fudan.edu.cn}{\texttt{23300180152@m.fudan.edu.cn}}\par
\end{samepage}
\end{document}